\documentstyle[amssymb,12pt,a4]{lamuphys}

\makeatletter

\def\limfunc#1{\mathop{\rm #1}}%

\long\def\QQQ#1#2{%
     \long\expandafter\def\csname#1\endcsname{#2}}%
\@ifundefined{QTP}{\def\QTP#1{}}{}
\@ifundefined{Qcb}{}{}
\@ifundefined{Qct}{}{}
\@ifundefined{Qlb}{}{}
\@ifundefined{Qlt}{}{}
\long\def\QQA#1#2{}%
\def\QTR#1#2{{\csname#1\endcsname #2}}
\def\EXPAND#1[#2]#3{}%
\def\NOEXPAND#1[#2]#3{}%
\def\LaTeXparent#1{}%
\def\ChildStyles#1{}%
\def\ChildDefaults#1{}%
\def\QTagDef#1#2#3{}%
\def\QQfnmark#1{\footnotemark}

\@ifundefined{INDEX}{\def\INDEX#1#2{}{}}{}%
\@ifundefined{SUBINDEX}{\def\SUBINDEX#1#2#3{}{}{}}{}%
\def\initial#1{\bigbreak{\raggedright\large\bf #1}\kern 2\p@
   \penalty3000}%
\@ifundefined{ZZZ}{}{\makeatletter\input gnuindex.sty\makeatother\makeindex\makeatletter}%
\@ifundefined{abstract}{%
 \def\abstract{%
  \if@twocolumn
   \section*{Abstract (Not appropriate in this style!)}%
   \else \small 
   \begin{center}{\bf Abstract\vspace{-.5em}\vspace{\z@}}\end{center}%
   \quotation 
   \fi
  }%
 }{%
 }%
\@ifundefined{endabstract}{\def\endabstract
  {\if@twocolumn\else\endquotation\fi}}{}%
\@ifundefined{maketitle}{\def\maketitle#1{}}{}%
\@ifundefined{affiliation}{\def\affiliation#1{}}{}%
\@ifundefined{proof}{}{}%
\@ifundefined{endproof}{}{}%
\@ifundefined{newfield}{\def\newfield#1#2{}}{}%
\@ifundefined{chapter}{\def\chapter#1{\par(Chapter head:)#1\par }%
 \newcount\c@chapter}{}%
\@ifundefined{part}{\def\part#1{\par(Part head:)#1\par }}{}%
\@ifundefined{section}{\def\section#1{\par(Section head:)#1\par }}{}%
\@ifundefined{subsection}{\def\subsection#1%
 {\par(Subsection head:)#1\par }}{}%
\@ifundefined{subsubsection}{\def\subsubsection#1%
 {\par(Subsubsection head:)#1\par }}{}%
\@ifundefined{paragraph}{\def\paragraph#1%
 {\par(Subsubsubsection head:)#1\par }}{}%
\@ifundefined{subparagraph}{\def\subparagraph#1%
 {\par(Subsubsubsubsection head:)#1\par }}{}%
\@ifundefined{therefore}{}{}%
\@ifundefined{backepsilon}{}{}%
\@ifundefined{yen}{}{}%
\@ifundefined{registered}{%
   \def\registered{\relax\ifmmode{}\r@gistered
                    \else$\m@th\r@gistered$\fi}%
 \def\r@gistered{^{\ooalign
  {\hfil\raise.07ex\hbox{$\scriptstyle\rm\text{R}$}\hfil\crcr
  \mathhexbox20D}}}}{}%
\@ifundefined{Eth}{}{}%
\@ifundefined{eth}{}{}%
\@ifundefined{Thorn}{}{}%
\@ifundefined{thorn}{}{}%
\@ifundefined{degree}{}{}%
\def\BibTeX{{\rm B\kern-.05em{\sc i\kern-.025em b}\kern-.08em
    T\kern-.1667em\lower.7ex\hbox{E}\kern-.125emX}}%
\newdimen\theight
\def\Column{%
 \vadjust{\setbox\z@=\hbox{\scriptsize\quad\quad tcol}%
  \theight=\ht\z@\advance\theight by \dp\z@\advance\theight by \lineskip
  \kern -\theight \vbox to \theight{%
   \rightline{\rlap{\box\z@}}%
   \vss
   }%
  }%
 }%
\def\qed{%
 \ifhmode\unskip\nobreak\fi\ifmmode\ifinner\else\hskip5\p@\fi\fi
 \hbox{\hskip5\p@\vrule width4\p@ height6\p@ depth1.5\p@\hskip\p@}%
 }%
\def\miss{\hbox{\vrule height2\p@ width 2\p@ depth\z@}}%
\def\tcol#1{{\baselineskip=6\p@ \vcenter{#1}} \Column}  %
\def\newfmtname{LaTeX2e}
\def\chkcompat{%
   \if@compatibility
   \else
     \usepackage{latexsym}
   \fi
}

\ifx\fmtname\newfmtname
  \DeclareOldFontCommand{\rm}{\normalfont\rmfamily}{\mathrm}
  \DeclareOldFontCommand{\sf}{\normalfont\sffamily}{\mathsf}
  \DeclareOldFontCommand{\tt}{\normalfont\ttfamily}{\mathtt}
  \DeclareOldFontCommand{\bf}{\normalfont\bfseries}{\mathbf}
  \DeclareOldFontCommand{\it}{\normalfont\itshape}{\mathit}
  \DeclareOldFontCommand{\sl}{\normalfont\slshape}{\@nomath\sl}
  \DeclareOldFontCommand{\sc}{\normalfont\scshape}{\@nomath\sc}
  \chkcompat
\fi

\def\Greekmath#1#2#3#4{%
    \if@compatibility
        \ifnum\mathgroup=\symbold
          \mbox{\boldmath$\mathchar"#1#2#3#4$}
        \else
           \mathchar"#1#2#3#4%
        \fi 
    \else 
        \ifnum\mathgroup=5 
           \mbox{\boldmath$\mathchar"#1#2#3#4$}
        \else
           \mathchar"#1#2#3#4%
        \fi     	    
	  \fi}

\newif\ifGreekBold  \GreekBoldfalse
\let\SAVEPBF=\pbf
\def\pbf{\GreekBoldtrue\SAVEPBF}%

\@ifundefined{theorem}{\newtheorem{theorem}{Theorem}}{}
\@ifundefined{lemma}{\newtheorem{lemma}[theorem]{Lemma}}{}
\@ifundefined{corollary}{\newtheorem{corollary}[theorem]{Corollary}}{}
\@ifundefined{conjecture}{}{}
\@ifundefined{proposition}{\newtheorem{proposition}[theorem]{Proposition}}{}
\@ifundefined{axiom}{}{}
\@ifundefined{remark}{\newtheorem{remark}{Remark}}{}
\@ifundefined{example}{\newtheorem{example}{Example}}{}
\@ifundefined{exercise}{}{}
\@ifundefined{definition}{\newtheorem{definition}{Definition}}{}

\expandafter\ifx\csname ds@amstex\endcsname\relax
\else\makeatother\endinput\fi

\def\RIfM@{\relax\protect\ifmmode}
\def\text{\RIfM@\expandafter\text@\else\expandafter\mbox\fi}
\let\nfss@text\text
\def\text@#1{\mathchoice
   {\textdef@\displaystyle\f@size{#1}}%
   {\textdef@\textstyle\tf@size{\firstchoice@false #1}}%
   {\textdef@\textstyle\sf@size{\firstchoice@false #1}}%
   {\textdef@\textstyle \ssf@size{\firstchoice@false #1}}%
   \glb@settings}

\def\textdef@#1#2#3{\hbox{{%
                    \everymath{#1}%
                    \let\f@size#2\selectfont
                    #3}}}
\newif\iffirstchoice@
\firstchoice@true
\def\Let@{\relax\iffalse{\fi\let\\=\cr\iffalse}\fi}%
\def\vspace@{\def\vspace##1{\crcr\noalign{\vskip##1\relax}}}%
\def\multilimits@{\bgroup\vspace@\Let@
 \baselineskip\fontdimen10 \scriptfont\tw@
 \advance\baselineskip\fontdimen12 \scriptfont\tw@
 \lineskip\thr@@\fontdimen8 \scriptfont\thr@@
 \lineskiplimit\lineskip
 \vbox\bgroup\ialign\bgroup\hfil$\m@th\scriptstyle{##}$\hfil\crcr}%
\def\Sb{_\multilimits@}%
\def\endSb{\crcr\egroup\egroup\egroup}%
\def\Sp{^\multilimits@}%

\newdimen\ex@
\def\rightarrowfill@#1{$#1\m@th\mathord-\mkern-6mu\cleaders
 \hbox{$#1\mkern-2mu\mathord-\mkern-2mu$}\hfill
 \mkern-6mu\mathord\rightarrow$}%
\def\leftarrowfill@#1{$#1\m@th\mathord\leftarrow\mkern-6mu\cleaders
 \hbox{$#1\mkern-2mu\mathord-\mkern-2mu$}\hfill\mkern-6mu\mathord-$}%
\def\leftrightarrowfill@#1{$#1\m@th\mathord\leftarrow
\mkern-6mu\cleaders
 \hbox{$#1\mkern-2mu\mathord-\mkern-2mu$}\hfill
 \mkern-6mu\mathord\rightarrow$}%
\def\overrightarrow{\mathpalette\overrightarrow@}%
\def\overrightarrow@#1#2{\vbox{\ialign{##\crcr\rightarrowfill@#1\crcr
 \noalign{\kern-\ex@\nointerlineskip}$\m@th\hfil#1#2\hfil$\crcr}}}%

\def\overleftarrow{\mathpalette\overleftarrow@}%
\def\overleftarrow@#1#2{\vbox{\ialign{##\crcr\leftarrowfill@#1\crcr
 \noalign{\kern-\ex@\nointerlineskip}$\m@th\hfil#1#2\hfil$\crcr}}}%
\def\overleftrightarrow{\mathpalette\overleftrightarrow@}%
\def\overleftrightarrow@#1#2{\vbox{\ialign{##\crcr
   \leftrightarrowfill@#1\crcr
 \noalign{\kern-\ex@\nointerlineskip}$\m@th\hfil#1#2\hfil$\crcr}}}%
\def\underrightarrow{\mathpalette\underrightarrow@}%
\def\underrightarrow@#1#2{\vtop{\ialign{##\crcr$\m@th\hfil#1#2\hfil
  $\crcr\noalign{\nointerlineskip}\rightarrowfill@#1\crcr}}}%

\def\underleftarrow{\mathpalette\underleftarrow@}%
\def\underleftarrow@#1#2{\vtop{\ialign{##\crcr$\m@th\hfil#1#2\hfil
  $\crcr\noalign{\nointerlineskip}\leftarrowfill@#1\crcr}}}%
\def\underleftrightarrow{\mathpalette\underleftrightarrow@}%
\def\underleftrightarrow@#1#2{\vtop{\ialign{##\crcr$\m@th
  \hfil#1#2\hfil$\crcr
 \noalign{\nointerlineskip}\leftrightarrowfill@#1\crcr}}}%

\def\qopnamewl@#1{\mathop{\operator@font#1}\nlimits@}
\let\nlimits@\displaylimits
\def\setboxz@h{\setbox\z@\hbox}

\def\varlim@#1#2{\mathop{\vtop{\ialign{##\crcr
 \hfil$#1\m@th\operator@font lim$\hfil\crcr
 \noalign{\nointerlineskip}#2#1\crcr
 \noalign{\nointerlineskip\kern-\ex@}\crcr}}}}

 \def\rightarrowfill@#1{\m@th\setboxz@h{$#1-$}\ht\z@\z@
  $#1\copy\z@\mkern-6mu\cleaders
  \hbox{$#1\mkern-2mu\box\z@\mkern-2mu$}\hfill
  \mkern-6mu\mathord\rightarrow$}
\def\leftarrowfill@#1{\m@th\setboxz@h{$#1-$}\ht\z@\z@
  $#1\mathord\leftarrow\mkern-6mu\cleaders
  \hbox{$#1\mkern-2mu\copy\z@\mkern-2mu$}\hfill
  \mkern-6mu\box\z@$}

\def\projlim{\qopnamewl@{proj\,lim}}
\def\injlim{\qopnamewl@{inj\,lim}}
\def\varinjlim{\mathpalette\varlim@\rightarrowfill@}
\def\varprojlim{\mathpalette\varlim@\leftarrowfill@}
\def\varliminf{\mathpalette\varliminf@{}}
\def\varliminf@#1{\mathop{\underline{\vrule\@depth.2\ex@\@width\z@
   \hbox{$#1\m@th\operator@font lim$}}}}
\def\varlimsup{\mathpalette\varlimsup@{}}
\def\varlimsup@#1{\mathop{\overline
  {\hbox{$#1\m@th\operator@font lim$}}}}

\def\binom#1#2{{#1 \choose #2}}%
\def\QATOPD#1#2#3#4{{#3 \atopwithdelims#1#2 #4}}%
\def\stackunder#1#2{\mathrel{\mathop{#2}\limits_{#1}}}%
\begingroup \catcode `|=0 \catcode `[= 1
\catcode`]=2 \catcode `\{=12 \catcode `\}=12
\catcode`\\=12 
|gdef|@alignverbatim#1\end{align}[#1|end[align]]
|gdef|@salignverbatim#1\end{align*}[#1|end[align*]]

|gdef|@alignatverbatim#1\end{alignat}[#1|end[alignat]]
|gdef|@salignatverbatim#1\end{alignat*}[#1|end[alignat*]]

|gdef|@xalignatverbatim#1\end{xalignat}[#1|end[xalignat]]
|gdef|@sxalignatverbatim#1\end{xalignat*}[#1|end[xalignat*]]

|gdef|@gatherverbatim#1\end{gather}[#1|end[gather]]
|gdef|@sgatherverbatim#1\end{gather*}[#1|end[gather*]]

|gdef|@gatherverbatim#1\end{gather}[#1|end[gather]]
|gdef|@sgatherverbatim#1\end{gather*}[#1|end[gather*]]

|gdef|@multilineverbatim#1\end{multiline}[#1|end[multiline]]
|gdef|@smultilineverbatim#1\end{multiline*}[#1|end[multiline*]]

|gdef|@arraxverbatim#1\end{arrax}[#1|end[arrax]]
|gdef|@sarraxverbatim#1\end{arrax*}[#1|end[arrax*]]

|gdef|@tabulaxverbatim#1\end{tabulax}[#1|end[tabulax]]
|gdef|@stabulaxverbatim#1\end{tabulax*}[#1|end[tabulax*]]

|endgroup

\def\align{\@verbatim \frenchspacing\@vobeyspaces \@alignverbatim
You are using the "align" environment in a style in which it is not defined.}

\@namedef{align*}{\@verbatim\@salignverbatim
You are using the "align*" environment in a style in which it is not defined.}
\expandafter\let\csname endalign*\endcsname =\endtrivlist

\def\alignat{\@verbatim \frenchspacing\@vobeyspaces \@alignatverbatim
You are using the "alignat" environment in a style in which it is not defined.}

\@namedef{alignat*}{\@verbatim\@salignatverbatim
You are using the "alignat*" environment in a style in which it is not defined.}
\expandafter\let\csname endalignat*\endcsname =\endtrivlist

\def\xalignat{\@verbatim \frenchspacing\@vobeyspaces \@xalignatverbatim
You are using the "xalignat" environment in a style in which it is not defined.}

\@namedef{xalignat*}{\@verbatim\@sxalignatverbatim
You are using the "xalignat*" environment in a style in which it is not defined.}
\expandafter\let\csname endxalignat*\endcsname =\endtrivlist

\def\gather{\@verbatim \frenchspacing\@vobeyspaces \@gatherverbatim
You are using the "gather" environment in a style in which it is not defined.}

\@namedef{gather*}{\@verbatim\@sgatherverbatim
You are using the "gather*" environment in a style in which it is not defined.}
\expandafter\let\csname endgather*\endcsname =\endtrivlist

\def\multiline{\@verbatim \frenchspacing\@vobeyspaces \@multilineverbatim
You are using the "multiline" environment in a style in which it is not defined.}

\@namedef{multiline*}{\@verbatim\@smultilineverbatim
You are using the "multiline*" environment in a style in which it is not defined.}
\expandafter\let\csname endmultiline*\endcsname =\endtrivlist

\def\arrax{\@verbatim \frenchspacing\@vobeyspaces \@arraxverbatim
You are using a type of "array" construct that is only allowed in AmS-LaTeX.}

\def\tabulax{\@verbatim \frenchspacing\@vobeyspaces \@tabulaxverbatim
You are using a type of "tabular" construct that is only allowed in AmS-LaTeX.}

\@namedef{arrax*}{\@verbatim\@sarraxverbatim
You are using a type of "array*" construct that is only allowed in AmS-LaTeX.}
\expandafter\let\csname endarrax*\endcsname =\endtrivlist

\@namedef{tabulax*}{\@verbatim\@stabulaxverbatim
You are using a type of "tabular*" construct that is only allowed in AmS-LaTeX.}
\expandafter\let\csname endtabulax*\endcsname =\endtrivlist

\def\@@eqncr{\let\@tempa\relax
    \ifcase\@eqcnt \def\@tempa{& & &}\or \def\@tempa{& &}%
      \else \def\@tempa{&}\fi
     \@tempa
     \if@eqnsw
        \iftag@
           \@taggnum
        \else
           \@eqnnum\stepcounter{equation}\fi
     \fi
     \global\tag@false
     \global\@eqnswtrue
     \global\@eqcnt\z@\cr}

 \def\endequation{%
     \iftag@
        \addtocounter{equation}{-1} 
        \eqno \hbox{\@taggnum}
        \global\tag@false%
        $$\global\@ignoretrue
     \else
        \eqno \hbox{\@eqnnum}
        $$\global\@ignoretrue
     \fi
 } 

 \newif\iftag@ \tag@false
 
 \def\tag{\@ifnextchar*{\@tagstar}{\@tag}}
 \def\@tag#1{%
     \global\tag@true
     \global\def\@taggnum{(#1)}}
 \def\@tagstar*#1{%
     \global\tag@true
     \global\def\@taggnum{#1}%
}

\@addtoreset{equation}{section}

\@addtoreset{theorem}{section}

\@addtoreset{example}{section}

\@addtoreset{definition}{section}

\@addtoreset{remark}{section}

\makeatother
\begin{document}

\title{Quantum and Classical Integrable Systems}
\author{M.A. Semenov-Tian-Shansky}

\institute{Physique Math\'ematique, Universit\'e de Bourgogne, Dijon, France%
\and Steklov Mathematical Institute, St.Petersburg, Russia}

\maketitle

\begin{abstract}
The key concept discussed in these lectures is the relation between the
Hamiltonians of a quantum integrable system and the Casimir elements in the
underlying hidden symmetry algebra. (In typical applications the latter is
either the universal enveloping algebra of an affine Lie algebra, or its
q-deformation.) A similar relation also holds in the classical case. We
discuss different guises of this very important relation and its implication
for the description of the spectrum and the eigenfunctions of the quantum
system. Parallels between the classical and the quantum cases are thoroughly
discussed.
\end{abstract}


\section{Introduction}

The study of exactly solvable quantum mechanical models is at least as old
as Quantum Mechanics itself. Over the last fifteen years there has been a major
development aimed at a unified treatment of many examples known previously
and, more importantly, at a systematic construction of new ones. The new
method nicknamed the Quantum Inverse Scattering Method was largely created
by L. D. Faddeev and his school in St. Petersburg as a quantum counterpart of
the Classical Inverse Scattering Method, and has brought together many ideas
believed to be unrelated. [Besides the Classical Inverse Scattering Method,
one should mention the profound results of  R.Baxter in Quantum Statistical
Mechanics (which, in turn, go back to the work of L. Onsager, E. Lieb and many
others) and the seminal paper of H. Bethe on the ferromagnet model in which
the now famous Bethe Ansatz was introduced.] It also allowed to unravel
highly nontrivial algebraic structures, the Quantum Group Theory being one
of its by-products.

The origins of QISM lie in the study of concrete examples; it is designed as
a working machine which produces quantum systems together with their
spectra, the quantum integrals of motion, and their joint eigenvectors. In
the same spirit, the Classical Inverse Scattering Method (along with its
ramifications) is a similar tool to produce examples of classical integrable
systems together with their solutions. In these lectures it is virtually impossible to follow
the historic way in which this method has been developed, starting with the
famous papers of \cite{krus}, and \cite{Lax}. (A good introduction close to the
ideas of QISM may be found in (\cite{FadTakh}).) I would like to comment
only on one important turning point which gave the impetus to invent the
Quantum Inverse Scattering Method. In 1979 L. D. Faddeev exposed in his
seminar the draft paper of B. Kostant on the quantization of the Toda lattice
(\cite{kost}) which reduced the problem to the representation theory of
semisimple Lie groups. This was an indication that classical completely
integrable systems have an intimate relation to Lie groups and should
also have exactly solvable quantum counterparts. In the same talk Faddeev
introduced the now famous $RL_1L_2=L_2L_1R$ commutation relations for
quantum Lax matrices which could be extracted from R. Baxter's works on
quantum transfer matrices. The research program outlined in that talk and
implemented over the next few years was twofold: On the one hand, the
Yang-Baxter equation and the related algebra have directly led to exact
solutions of several quantum models, such as the quantum sine-Gordon
equation. On the other hand, connections with group theory and the orbits
method resulted in a systematic treatment of numerous examples in the
classical setting (\cite{AreyS}). One major problem which has remained unsettled
for more than a decade is to fill the gap between the two approaches, and in
particular, to explain the group theoretical meaning of the Bethe Ansatz.

With  hindsight we can now understand why this problem could not be
resolved immediately. First of all, while classical integrable systems are
related to ordinary Lie groups, quantum systems quite often (though not
always, cf. the discussion in section \ref{quadr}) require the full
machinery of Quantum Groups. The background took several years to prepare (%
\cite{QG}). Second, even the simplest systems such as the open Toda lattice
require a very advanced technique of representation theory (\cite{kost}, 
\cite{STSToda}). The Toda lattice is peculiar, since the underlying 'hidden
symmetry' group is finite-dimensional. All main examples are related to
(classical or quantum) affine Lie algebras. The representation theory of
affine algebras which is instrumental in the study of integrable systems has
been developed only in recent years (semi-infinite cohomologies,Wakimoto
modules, critical level representations, cf. \cite{FeiFr}, \cite{FFR}, \cite
{EFrenkel}).

The present lectures do not give a systematic overview of the Quantum
Inverse Scattering Method (several good expositions are available, cf. \cite
{BIK}, Faddeev (\cite{Fadd}, \cite{Leshouches}, \cite{bethe}), \cite{FT}, 
\cite{KulishSkl}, \cite{Skl}). Instead, I shall try to explain the parallels
between quantum and classical systems and the 'correspondence principles'
which relate the quantum and the classical case.

As already mentioned, the first key observation is that integrable systems
always have an ample hidden symmetry. (Fixing this symmetry provides some
rough classification of the associated examples. A nontrivial class of
examples is related to loop algebras or, more generally, to their
q-deformations. This is the class of examples I shall consider below.) Two
other key points are the role of (classical or quantum) R-matrices and of
the Casimir elements which give rise to the integrals of motion.

This picture appears in several different guises depending on the type of
examples in question. The simplest (so-called linear) case corresponds to
classical systems which are modelled on coadjoint orbits of Lie algebras; a
slightly more complicated group of examples is classical systems modelled
on Poisson Lie groups or their Poisson submanifolds. In the quantum setting
the difference between these two cases is deeper: while in the former case
the hidden symmetry algebra remains the same, quantization of the latter
leads to Quantum Groups. Still, for loop algebras the quantum counterpart of
the main construction is nontrivial even in the linear case; the point is
that the universal enveloping algebra of a loop algebra has a trivial center
which reappears only after central extension at the critical value of the
central charge. Thus to tackle the quantum case one needs the full machinery
of representation theory of loop algebras. (By contrast, in the classical
case one mainly deals with the evaluation representations which allow to
reduce the solution of the equations of motion to a problem in algebraic
geometry.)

The study of integrable models may be divided into two different parts. The
first one is, so to say, kinematic: it consists in the choice of appropriate
models together with their phase spaces or the algebra of observables and of
their Hamiltonians. The second one is dynamical; it consists, classically,
in the description of solutions or of the action-angle variables. The
quantum counterpart consists in the description of the spectra, the
eigenvectors, and of the various correlation functions. The algebraic scheme
proved to be very efficient in the description of kinematics. (Quantum Group
Theory may be regarded as a by-product of this kinematic problem.) The
description of spectra at the present stage of the theory remains model
dependent. The standard tool for constructing  the eigenvectors of the quantum
Hamiltonians which was the starting point of QISM is the algebraic Bethe
Ansatz. Until very recently, its interpretation in terms of the
representation theory was lacking. This problem has been finally settled by 
\cite{FFR} for an important particular model with {\em linear} commutation
relations ({\em the Gaudin model}); remarkably, their results follow the
general pattern outlined above. A similar treatment of quantum models
related to q-deformed affine algebras is also possible, although the results
in this case are still incomplete. One should be warned that much of the
'experimental material' on Quantum Integrability still resists general
explanations. I would like to mention in this respect the deep results of
E.K.Sklyanin (\cite{Skl}, \cite{separ}) relating the Bethe Ansatz to the
separation of variables; see also (\cite{HW}, \cite{Kuz}).

\subsubsection*{Acknowledgements.}

The present lectures were prepared for the CIMPA Winter School on Nonlinear
Systems which was held in Pondicherry, India, in January 1996. I am
deeply grateful to the Organizing Committee of the School and to the staff
of the Pondicherry University for their kind hospitality. I would also like to thank 
Prof. Y. Kosmann-Schwarzbach for her remarks on the draft text of these lectures.

\section{\label{gen}Generalities}

A quantum mechanical system is a triple consisting of an associative algebra
with involution ${\cal A}$ $\ $ called the algebra of observables, an
irreducible $\star $-representation $\pi $ of ${\cal A}$ in a Hilbert space $%
{\frak H}$ and a distinguished self-adjoint observable ${\cal H}$ called the
Hamiltonian. A typical question to study is the description of the spectrum
and the eigenvectors of $\pi ({\cal H})$. Quantum mechanical systems usually
appear with their classical counterparts. By definition, a classical
mechanical system is again specified by its algebra of observables ${\cal A}%
_{cl}$ which is a commutative associative algebra equipped with a Poisson
bracket (i.e., a Lie bracket which also satisfies the Leibniz rule 
\[
\left\{ a,bc\right\} =\left\{ a,b\right\} c+\left\{ a,c\right\} b; 
\]
in other words, a Poisson bracket is a derivation of ${\cal A}_{cl}$ with
respect to both its arguments), and a Hamiltonian ${\cal H}\in {\cal A}%
_{cl}. $ A commutative algebra equipped with a Poisson bracket satisfying
the Leibniz rule is called a {\em Poisson algebra}. Speaking informally, a
quantum algebra of observables ${\cal A}$ arises as a deformation of the
commutative algebra ${\cal A}_{cl}$ determined by the Poisson bracket. In
these lectures we shall not be concerned with the quantization problem in
its full generality (cf. \cite{Flato}, \cite{Weinstein}). However, it will
always be instructive to compare quantum systems with their classical
counterparts.

The algebraic language which starts with Poisson algebras makes the gap
between classical and quantum mechanics as narrow as possible; in practice,
however, we also need the dual language based on the notion of the phase
space. Roughly speaking, the phase space is the spectrum of the Poisson
algebra. Accurate definition depends on the choice of a topology in the
Poisson algebra. We shall not attempt to discuss these subtleties and shall always
assume that the underlying phase space is a smooth manifold and that the Poisson
algebra is realized as the algebra of functions on this manifold. In the
examples we have in mind, Poisson algebras always have an explicit geometric
realization of this type. The Poisson bracket itself is then as usual
represented by a bivector field on the phase space satisfying certain
differential constraints which account for the Jacobi identity. This gives
the definition of a {\em Poisson manifold} which is dual to the notion of a
Poisson algebra. Geometry of Poisson manifolds has numerous obvious parallels
with the representation theory. Recall that the algebraic version of the representation
theory is based on the study of appropriate ideals in an associative
algebra. For a Poisson algebra we have a natural notion of a {\em Poisson
ideal} (i.e., a subalgebra which is an ideal with respect to both
structures); the dual notion is that of a {\em Poisson submanifold} of a
Poisson manifold. The classical counterpart of Hilbert space representations
of an associative algebra is the restriction of functions to various Poisson
submanifolds. Poisson submanifolds are partially ordered by inclusion; {\em %
minimal} Poisson submanifolds are those for which the induced Poisson
structure is nondegenerate. (This means that the center of the Lie algebra of
functions contains only constants.) Minimal Poisson submanifolds always
carry a symplectic structure and form a stratification of the Poisson
manifold; they are called {\em symplectic leaves}. Restriction of functions
to symplectic leaves gives a classical counterpart of  the irreducible
representations of associative algebras.

Let ${\cal M}$ be a Poisson manifold, \ ${\cal H}\in C^\infty \left( {\cal M}%
\right) $. A classical system $\left( {\cal M},{\cal H}\right) $ is called 
{\em integrable }if the commutant of the Hamiltonian ${\cal H}$ in ${\cal A}%
_{cl}$ contains an abelian algebra of the maximal possible rank. (A
technical definition is provided by the well known Liouville theorem.)

Let us recall that the key idea which has started the modern age in the
study of classical integrable systems is to bring them into {\em Lax form}.
In the simplest case, the definition of a Lax representation may be given as
follows. Let $\left( {\cal A},{\cal M},{\cal H}\right) $ be a classical
mechanical system. Let ${\cal F}_t:{\cal M}\rightarrow {\cal M}$ be the
associated flow on ${\cal M}$ (defined at least locally).\ Suppose that $%
{\frak g}$ is a Lie algebra. A mapping $L:{\cal M}\rightarrow {\frak g}$ is
called a{\em \ Lax representation} of $\left( {\cal A},{\cal M},{\cal H}%
\right) $ if the following conditions are satisfied:

{\em (i) The flow} ${\cal F}_t$ {\em factorizes over} ${\frak g}$, {\em i.e.,
there exists a (local) flow} $F_t:{\frak g}\rightarrow {\frak g}$ {\em such
that the following diagram is commutative.} 
\[
\begin{array}{ccccc}
& {\cal M} & \stackrel{{\cal F}_t}{\longrightarrow } & {\cal M} &  \\ 
& \ \mid \  &  & \mid &  \\ 
^L & \downarrow &  & \downarrow & ^L \\ 
& {\frak g} & \stackrel{F_t}{\longrightarrow } & {\frak g} & 
\end{array}
\]

{\em (ii)} {\em The quotient flow} $F_t$ {\em on} ${\frak g}$ {\em is
isospectral, i.e., it is tangent to adjoint orbits in} ${\frak g}$.

Clearly, $L(x)\in {\frak g}$ for any $x\in {\cal M}$ $;$ hence we may regard 
$L$ as a 'matrix with coefficients in ${\cal A}=C^\infty \left( {\cal M}%
\right) $'$,$ i.e., as an element of ${\frak g}$ $\otimes $${\cal A};$ the
Poisson bracket on ${\cal A}$ extends to ${\frak g}$ $\otimes $${\cal A}\ $%
by linearity. Property (ii) means that there exists an element $M\in {\bf 
{\frak g}\ \otimes }$ ${\cal A}$ such that $\left\{ {\cal H},L\right\}
=\left[ L,M\right] $.

Let $\left( \rho ,V\right) $ be a (finite-dimensional) linear representation
of ${\frak g}$. Then $L_V=\rho \otimes id(L)\in EndV\otimes {\cal A}$ is a
matrix valued function on ${\cal M}$; the coefficients of its characteristic
polynomial $P(\lambda )=\det (L_V-\lambda )$ are integrals of the motion.

One may replace in the above definition a Lie algebra ${\frak g}$ with a Lie
group $G;$ in that case isospectrality means that the flow preserves
conjugacy classes in $G.$ In a more general way, the Lax operator may be a
difference or a differential operator. The difference case is of particular
importance, since the quantization of difference Lax equations is the core of
QISM and a natural source of Quantum Groups. (Cf. section 3 below)

There is no general way to find a Lax representation for a given system
(even if it is known to be completely integrable). However, there is a
systematic way to produce {\em examples} of such representations. An ample
source of such examples is provided by the general construction described in
the next section.

\subsection{Basic Theorem: Linear Case}

The basic construction outlined in this section (summarized in Theorem \ref
{AKS} below) goes back to \cite{kost} and \cite{adler} in some crucial
cases; its relation with the r-matrix method was established in (\cite{rmatr}%
).\ We shall state it using the language of symmetric algebras which
simplifies its generalization to the quantum case. Let ${\frak g}$ be a Lie
algebra over $k$ (where $k={\Bbb R}$ or ${\Bbb C}$). Let ${\cal S}({\bf 
{\frak g}})$ be the symmetric algebra of ${\frak g}$. Recall that there is a
unique Poisson bracket on ${\cal S}\left( {\frak g}\right) $ (called the 
{\em Lie-Poisson bracket}) which extends the Lie bracket on ${\frak g}$
(see, for example, (\cite{Cartier}). This Poisson bracket is more frequently
discussed from the 'spectral' point of view. Namely, let ${\frak g}^{*}$ be
the linear dual of ${\frak g;}$ the natural pairing ${\frak g\times g}%
^{*}\rightarrow k$ extends to the 'evaluation map' ${\cal S}\left( {\frak g}%
\right) \times {\frak g}^{*}\rightarrow k$ which induces a canonical
isomorphism of ${\cal S}\left( {\frak g}\right) \ $with the space of
polynomials $P\left( {\frak g}^{*}\right) ;$ thus ${\frak g}^{*}$ is a
linear Poisson manifold and ${\cal S}\left( {\frak g}\right) \simeq $ $%
P\left( {\frak g}^{*}\right) $ is the corresponding algebra of observables.
Linear functions on ${\frak g}^{*}$ form a subspace in $P\left( {\frak g}%
^{*}\right)$ which may be identified with ${\frak g.}$ The Lie-Poisson
bracket on ${\frak g}^{*}$ is uniquely characterized by the following
properties:

\begin{enumerate}
\item  {\em The Poisson bracket of linear functions on }${\frak g}^{*}${\em \ is
again a linear function.}

\item  {\em The restriction of the Poisson bracket to }${\frak g}\subset
P\left( {\frak g}^{*}\right) ${\em \ coincides with the Lie bracket in} $%
{\frak g.}$
\end{enumerate}

Besides the Lie-Poisson structure on ${\cal S}({\frak g})$ we shall need its
Hopf structure. Further on we shall deal with other more complicated
examples of Hopf algebras, so it is probably worth recalling the general
definitions (though we shall not use them in full generality until section
4). Recall that a Hopf algebra is a set $\left( A,m,\Delta ,\epsilon
,S\right) $ consisting of an associative algebra $A$ over $k$\ with
multiplication $m:A\otimes A\rightarrow A\ $and unit element $1,$ the
coproduct $\Delta :A\rightarrow A\otimes A,$ the counit $\epsilon
:A\rightarrow k$ and the antipode $S:A\rightarrow A$ which satisfy the
following axioms

\begin{itemize}
\item  $m:A\otimes A\rightarrow A,\Delta :A\rightarrow A\otimes
A,i:k\rightarrow A:\alpha \longmapsto \alpha \cdot 1,\epsilon :A\rightarrow
k\ ${\em are homomorphisms of algebras.}

\item  {\em The following diagrams are commutative:}
\end{itemize}

\begin{equation}
\begin{array}{ccccccccc}
A & \otimes & A & \otimes & A & \stackrel{id\otimes m}{\rightarrow } & A & 
\otimes & A \\ 
&  & \mid &  &  &  &  & \mid &  \\ 
& ^{m\otimes id} & \downarrow &  &  &  &  & \downarrow &  \\ 
& A & \otimes & A &  & \rightarrow &  & A & 
\end{array}
,\; 
\begin{array}{ccccccccc}
& A &  & \rightarrow &  & A & \otimes & A &  \\ 
& \mid &  &  &  &  & \mid &  &  \\ 
^\Delta & \downarrow &  &  &  &  & \downarrow & ^{\Delta \otimes id} &  \\ 
A & \otimes & A & \stackrel{id\otimes \Delta }{\rightarrow } & A & \otimes & 
A & \otimes & A
\end{array}
,
\end{equation}
(these diagrams express associativity of the product $m$ and the
coassociativity of the coproduct $\Delta ,$ respectively),

\begin{equation}
\begin{array}{ccccccc}
A & \otimes & A & \stackrel{id\otimes i}{\longleftarrow } & A & \otimes & k
\\ 
& \mid &  &  &  & \mid &  \\ 
^m & \downarrow &  &  &  & \downarrow &  \\ 
& A &  & \stackrel{id}{\longleftarrow } &  & A & 
\end{array}
,\; 
\begin{array}{ccccccc}
A & \otimes & A & \stackrel{i\otimes id}{\longleftarrow } & A & \otimes & k
\\ 
& \mid &  &  &  & \mid &  \\ 
^m & \downarrow &  &  &  & \downarrow &  \\ 
& A &  & \stackrel{id}{\longleftarrow } &  & A & 
\end{array}
,
\end{equation}

\begin{equation}
\begin{array}{ccccccc}
A & \otimes & A & \stackrel{id\otimes \epsilon }{\longrightarrow } & A & 
\otimes & k \\ 
_\Delta & \uparrow &  &  &  & \mid &  \\ 
& \mid &  &  &  & \downarrow &  \\ 
& A &  & \stackrel{id}{\longrightarrow } &  & A & 
\end{array}
,\; 
\begin{array}{ccccccc}
A & \otimes & A & \stackrel{\epsilon \otimes id}{\longrightarrow } & A & 
\otimes & k \\ 
_\Delta & \uparrow &  &  &  & \mid &  \\ 
& \mid &  &  &  & \downarrow &  \\ 
& A &  & \stackrel{id}{\longrightarrow } &  & A & 
\end{array}
\end{equation}
(these diagrams express, respectively, the properties of the unit element $%
1\in A$ and of the counit $\epsilon \in A^{*}$ ).

\begin{itemize}
\item  {\em The antipode }$S${\em \ is an antihomomorphism of algebras and
the following diagrams are commutative: } 
\begin{equation}
\begin{array}{ccccccc}
A & \stackrel{\Delta }{\longrightarrow } & A\otimes A & \stackrel{id\otimes S%
}{\longrightarrow } & A\otimes A & \stackrel{m}{\longrightarrow } & A \\ 
&  &  &  &  &  &  \\ 
&  & \searrow ^\epsilon &  & \;^i\nearrow &  &  \\ 
&  &  & k &  &  & 
\end{array}
,
\end{equation}
\end{itemize}

\begin{equation}
\begin{array}{ccccccc}
A & \stackrel{\Delta }{\longrightarrow } & A\otimes A & \stackrel{S\otimes id%
}{\longrightarrow } & A\otimes A & \stackrel{m}{\longrightarrow } & A \\ 
& \; &  &  &  &  &  \\ 
&  & \searrow ^\epsilon &  & \;\nearrow _i &  &  \\ 
&  &  & k &  &  & 
\end{array}
.
\end{equation}

If $\left( A,m,1,\Delta ,\epsilon ,S\right) $ is a Hopf algebra, its linear
dual $A^{*}$ is also a Hopf algebra; moreover, the coupling $A\otimes
A^{*}\rightarrow k$ interchanges the roles of product and coproduct. Thus we
have 
\[
\left\langle m\left( a\otimes b\right) ,\varphi \right\rangle =\left\langle
a\otimes b,\Delta \varphi \right\rangle . 
\]

The Hopf structure on the symmetric algebra ${\cal S}\left( {\frak g}\right) 
$ is determined by the structure of the additive group on ${\frak g}^{*}.$
Namely, let ${\cal S}\left( {\frak g}\right) \times {\frak g}^{*}\rightarrow
k:\left( a,X\right) \longmapsto a\left( X\right) $ be the evaluation map$;$
since ${\cal S}\left( {\frak g}\right) \otimes {\cal S}\left( {\frak g}%
\right) \simeq {\cal S}\left( {\frak g}\oplus {\frak g}\right) \simeq
P\left( {\frak g}^{*}\oplus {\frak g}^{*}\right) ,$ this map extends to $%
{\cal S}\left( {\frak g}\right) \otimes {\cal S}\left( {\frak g}\right) .$
We have 
\begin{equation}
\Delta a\left( X,Y\right) =a\left( X+Y\right) ,\epsilon \left( a\right)
=a\left( 0\right) ,S\left( a\right) \left( X\right) =a\left( -X\right)
\end{equation}
(we shall sometimes write $S\left( a\right) =a^{\prime }$ for brevity).

For future reference we recall also that the coproducts in the universal
enveloping algebra ${\cal U}\left( {\frak g}\right) $ and in ${\cal S}\left( 
{\frak g}\right) $ coincide (in other words, ${\cal U}\left( {\frak g}%
\right) $ is canonically isomorphic to ${\cal S}\left( {\frak g}\right) $ as
a coalgebra).

Let us now return to the discussion of the Poisson structure on ${\cal S}%
\left( {\frak g}\right) ,$ i.e., of the Lie-Poisson bracket of the Lie
algebra ${\frak g.}$ It is well known that the Lie-Poisson bracket is always
degenerate; its symplectic leaves are precisely the coadjoint orbits of the
corresponding Lie group.

\begin{proposition}
The center of ${\cal S}({\frak g})$ (regarded as a Poisson algebra)
coincides with the subalgebra $I={\cal S}({\frak g})^{{\frak g}}$ of $ad%
{\frak g}$-invariants in ${\cal S}({\frak g})$ (called Casimir elements).
\end{proposition}

Restrictions of Casimir elements to symplectic leaves are constants; so
fixing the values of Casimirs provides a rough classification of symplectic
leaves, although it is not true in general that different symplectic leaves
are separated by the Casimirs. (In the semisimple case there are enough
Casimirs to separate generic orbits.)

Fix an element ${\cal H}\in {\cal S}\left( {\frak g}\right) ;$ it defines a
derivation \ $D_{{\cal H}}$ of ${\cal S}\left( {\frak g}\right) ,$%
\[
D_{{\cal H}}\varphi =\left\{ {\cal H},\varphi \right\} .
\]
By duality, this derivation determines a (local) Hamiltonian flow on ${\frak %
g}^{*}.$ Since the Lie-Poisson bracket on ${\frak g}^{*}$ is degenerate,
this flow splits into a family of independent flows which are confined to
Poisson submanifolds in ${\frak g}^{*}.$ Thus a Hamiltonian ${\cal H}$
defines Hamiltonian flows on the coadjoint orbits in ${\frak g}^{*}.$ In the
special case when ${\cal H}\in I$ all these flows are trivial. However, one
may still use Casimir elements to produce nontrivial equations of motion if
the Poisson structure on ${\cal S}({\frak g})$ is properly modified, and
this is the way in which Lax equations associated with ${\frak g}$ do arise.
The formal definition is as follows.

Let $r\in End\ {\frak g}${\bf \ }be a linear operator; we shall say that $r$
is a {\em classical r-matrix} if the {\em r-bracket} 
\begin{equation}
\left[ X,Y\right] _r=\frac 12\left( \left[ rX,Y\right] +\left[ X,rY\right]
\right) ,X,Y\in {\frak g},  \label{rbr}
\end{equation}
satisfies the Jacobi identity. In that case we get two different Lie
brackets on the same underlying linear space ${\frak g.}$ We shall assume
that $r$ satisfies the following stronger condition called {\em the} {\em %
modified classical Yang-Baxter identity:} 
\begin{equation}
\left[ rX,rY\right] -r\left( \left[ rX,Y\right] +\left[ X,rY\right] \right)
=-\left[ X,Y\right] ,X,Y\in {\bf {\frak g}.}  \label{cybe}
\end{equation}
{\bf \ }

\begin{proposition}
Identity (\ref{cybe}) implies the Jacobi identity for (\ref{rbr}).
\end{proposition}

\begin{remark}
Identity (\ref{cybe}) is of course only a sufficient condition; however,
we prefer to impose it from the very beginning, since it assures a very
important factorization property (see below).
\end{remark}

Let $\left( {\frak g},r\right) $ be a Lie algebra equipped with a classical
r-matrix $r\in End{\frak g}$ satisfying $\left( \ref{cybe}\right) .$ Let $%
{\frak g}_r$ be the corresponding Lie algebra (with the same underlying
linear space). Put $r_{\pm }=\frac 12(r\pm id);$ then (\ref{cybe}) implies
that $r_{\pm }:{\frak g}_r\rightarrow {\frak g}$ are Lie algebra
homomorphisms. Let us extend them to Poisson algebra morphisms ${\cal S}$($%
{\frak g}_r)\rightarrow {\cal S}$(${\frak g})$ which we shall denote by the
same letters. These morphisms also agree with the standard Hopf structure on 
${\cal S}$(${\frak g})$, ${\cal S}$(${\bf {\frak g}}_r).$ Define the action 
\[
{\cal S}({\frak g}_r)\otimes {\cal S}({\frak g})\rightarrow {\cal S}({\frak g%
}) 
\]
by setting 
\begin{equation}
x\cdot y=\sum r_{+}(x_i^{(1)})\;y\;r_{-}(x_i^{(2)})^{\prime },x\in {\cal S}(%
{\frak g}_r),y\in {\cal S}({\frak g}),  \label{act}
\end{equation}
where $\Delta x=\sum x_i^{(1)}\otimes x_i^{(2)}$ is the coproduct and $%
a\mapsto a^{\prime }$ is the antipode map.

\begin{theorem}
\label{AKS} (i) ${\cal S}$(${\frak g})$ is a free graded ${\cal S}$(${\frak g%
}_r)$-module generated by $1\in {\cal S}$(${\frak g}).$ (ii) Let $i_r$: $%
{\cal S}$(${\frak g})$ $\rightarrow {\cal S}$(${\bf {\frak g}}_r)$ be the
induced isomorphism of graded linear spaces; its restriction to $I={\cal S}(%
{\frak g}){\bf ^{{\frak g}}}$ is a morphism of Poisson algebras. (iii)
Assume, moreover, that ${\frak g}$ is equipped with a nondegenerate
invariant bilinear form; the induced mapping ${\frak g}_r^{*}\rightarrow 
{\frak g}$ defines a Lax representation for all Hamiltonians ${\cal H}=i_r(%
\widehat{{\cal H}}),\;\widehat{{\cal H}}\in $ ${\cal S}({\frak g}){\bf ^{%
{\frak g}}}$.
\end{theorem}

The geometric meaning of theorem \ref{AKS} is very simple. In fact, it
becomes more transparent if one uses the language of Poisson manifolds
rather than the dual language of Poisson algebras. (However, it is this more
cumbersome language that may be generalized to the quantum case.) Since $%
{\frak g}\simeq {\frak g}_r$ as linear spaces, the polynomial algebras $%
P\left( {\frak g}^{*}\right) $ and $P\left( {\frak g}_r^{*}\right) $ are
also isomorphic as graded linear spaces (though not as Poisson algebras).
The subalgebra $I\subset P\left( {\frak g}^{*}\right) $ remains commutative
with respect to the Lie-Poisson bracket of ${\frak g}_r;$ moreover, the
flows in ${\frak g}^{*}$ which correspond to the Hamiltonians ${\cal H}\in I$
are tangent to {\em two} systems of coadjoint orbits in ${\frak g}^{*},$ the
orbits of ${\frak g}$ and of ${\frak g}_r.$ The latter property is trivial,
since {\em all} Hamiltonian flows in ${\frak g}_r^{*}$ are tangent to
coadjoint orbits; the former is an immediate corollary of the fact that $I=%
{\cal S}({\frak g}){\bf ^{{\frak g}}\simeq }P\left( {\frak g}^{*}\right) ^{%
{\frak g}}$ is the subalgebra of coadjoint invariants. (As we mentioned, in
the semisimple case generic orbits are level surfaces of the Casimirs.)
Finally, if ${\frak g}$ carries an invariant inner product, the coadjoint
and the adjoint representations of ${\frak g}$ are equivalent and $I$
consists precisely of spectral invariants.

Let ${\frak u}\subset {\frak g}_r$ be an ideal, ${\frak s}={\frak g}_r/%
{\frak u}$ the quotient algebra, $p:{\cal S}\left( {\frak g}_r\right)
\rightarrow {\cal S}\left( {\frak s}\right) $ the canonical projection;
restricting $p $ to the subalgebra $I_r=i_r\left( {\cal S}({\frak g}){\bf ^{%
{\frak g}}}\right) $ we get a Poisson commutative subalgebra in ${\cal S}%
\left( {\frak s}\right) ;$ we shall say that the corresponding elements of $%
{\cal S}\left( {\frak s}\right) $ are obtained by {\em specialization.}

The most common examples of classical r-matrices are associated with {\em %
decompositions of Lie algebras. } Let ${\frak g}$ be a Lie algebra with a
nondegenerate invariant inner product, ${\frak g}_{+},{\frak g}_{-}\subset 
{\frak g}$ its two Lie subalgebras such that ${\frak g}={\frak g}_{+}\dot +%
{\frak g}_{-}$ as a linear space. Let $P_{{\frak +}},P_{{\frak -}}$ be the
projection operators associated with this decomposition; the operator 
\begin{equation}
r=P_{{\frak +}}-P_{{\frak -}}  \label{standard}
\end{equation}
satisfies (\ref{cybe}); moreover, the Lie algebra ${\frak g}_r$ splits into
two parts, ${\frak g}_r\simeq {\frak g}_{+}\oplus {\frak g}_{-}$.\ Let us
choose ${\cal A}={\cal S}({\frak g}_{+})$ as our algebra of observables. Let 
${\frak g}_{+}^{\perp }{\frak ,g}_{-}^{\perp }$ be the orthogonal
complements of ${\frak g_{+}},{\frak g}_{-}$ in ${\frak g.}$ Then ${\frak g=g%
}_{+}^{\perp }{\frak \oplus g}_{-}^{\perp }$ and we may identify ${\frak g}%
_{+}^{*}\simeq {\frak g}_{-}^{\perp },{\frak g}_{-}^{*}\simeq {\frak g}%
_{+}^{\perp }.$ Let $L\in {\frak g}_{-}^{\perp }\otimes {\frak g}_{+}\subset 
{\frak g}_{-}^{\perp }\otimes {\cal S}({\frak g}_{+})$ be the canonical
element (that is, $L$ $=\sum e_i\otimes e^i,$ where $\left\{ e^i\right\} $
is a basis in ${\frak g}_{+}$ and $\left\{ e_i\right\} $ the dual basis in $%
{\frak g}_{+}^{*}\simeq {\frak g}_{-}^{\bot }$); using the canonical
isomorphism ${\frak \ g}_{-}^{\bot }\otimes {\frak g}_{+}\simeq Hom\left( 
{\frak g}_{+}^{*},{\frak g}_{-}^{\bot }\right) $, we may regard $L$ as an
embedding ${\frak g}_{+}^{*}\hookrightarrow {\frak g}_{-}^{\perp }$. Let $P:%
{\cal S}({\frak g})\rightarrow {\cal S}({\frak g}_{-})$ be the projection
onto ${\cal S}({\frak g}_{-})$ in the decomposition 
\[
{\cal S}({\frak g})={\cal S}({\frak g}_{-})\oplus {\frak g}_{+}{\cal S}(%
{\frak g}){\bf .} 
\]

\begin{corollary}
\label{CorAKS}(i) The restriction of $P$ to the subalgebra $I={\cal S}(%
{\frak g}){\bf ^{{\frak g}}}$ of Casimir elements is a Poisson algebra
homomorphism. (ii) $L$ defines a Lax representation for all Hamiltonian
equations of motion defined by the Hamiltonians $\overline{{\cal H}}=P({\cal %
H}),\;{\cal H}\in I\ $. (iii) The corresponding Hamiltonian flows on ${\frak %
g}_{+}^{*}\simeq {\frak g}_{-}^{\perp }{\bf \subset }{\frak g}$ preserve
intersections of coadjoint orbits of ${\frak g}_{+}$ with the adjoint orbits
of ${\frak g}$.
\end{corollary}

One sometimes calls $L$ {\em the universal Lax operator;} restricting the
mapping $L:{\frak g}_{+}^{*}\hookrightarrow {\frak g}_{-}^{\perp }$ to
various Poisson submanifolds in ${\frak g}_{+}^{*}\ $we get Lax
representations for particular systems.

The situation is especially simple when ${\frak g}_{+},{\frak g}_{-}\subset 
{\frak g}$ are isotropic subspaces with respect to the inner product in $%
{\frak g;}$ then ${\frak g}_{+}=$ ${\frak g}_{+}^{\perp },{\frak g}_{-}=%
{\frak g}_{-}^{\perp },{\frak g}_{+}\simeq {\frak g}_{-}^{*}.$ In that case $%
\left( {\frak g}_{+}{\frak ,g_{-},g}\right) $ is referred to as a {\em Manin
triple}. This case is of particular importance, since then both ${\frak g}\ $
and its Lie subalgebras ${\frak g}_{+}{\frak ,g}_{-}$ carry an additional
structure of {\em Lie bialgebra}.

\begin{definition}
Let ${\frak a}$ be a Lie algebra, ${\frak a}^{*}$ its dual; assume that $%
{\frak a}^{*}$ is equipped with a Lie bracket $\left[ ,\right] _{*}:{\frak a}%
^{*}\wedge {\frak a}^{*}\rightarrow {\frak a}^{*}.$ The brackets on ${\frak a%
}$ and ${\frak a}^{*}$ are said to be compatible if the dual map $\delta :%
{\frak a}\rightarrow {\frak a}\wedge {\frak a}$ is a 1-cocycle on ${\frak a}$
(with values in the wedge square of the adjoint module). A pair $\left( 
{\frak a,a}^{*}\right) $ with consistent Lie brackets is called a Lie
bialgebra.
\end{definition}

One can prove that this definition is actually symmetric with respect to $%
{\frak a,a}^{*};$ in other words, if $\left( {\frak a,a}^{*}\right) $ is a
Lie bialgebra, so is $\left( {\frak a}^{*}{\frak ,a}\right) .$

\begin{proposition}
Let $\left( {\frak g}_{+}{\frak ,g_{-},g}\right) $ be a Manin triple.
Identify ${\frak g}^{*}$ with ${\frak g}$ by means of the inner product and
equip it with the r-bracket associated with $r=P_{{\frak +}}-P_{{\frak -}};$
then (i) $\left( {\frak g,g}^{*}\right) $ is a Lie bialgebra. (ii) $\left( 
{\frak g}_{+}{\frak ,g}_{-}\right) $ is a Lie sub-bialgebra of $\left( 
{\frak g,g}^{*}\right) .$ (iii) Conversely, if $\left( {\frak a},{\frak a}%
^{*}\right) $ is a Lie bialgebra there exists a unique Lie algebra ${\frak d}%
=d\left( {\frak a},{\frak a}^{*}\right) $ called the double of $\left( 
{\frak a},{\frak a}^{*}\right) $ such that
\end{proposition}

\begin{itemize}
\item  ${\frak d}={\frak a}\dot +{\frak a}^{*}${\em \ as a linear space},

\item  ${\frak a,a}^{*}\subset {\frak d}$ {\em are Lie subalgebras},

\item  {\em The canonical bilinear form on }${\frak d}$ i{\em nduced by the
natural pairing between} ${\frak a}$ {\em and} ${\frak a}^{*}${\em \ is} $ad%
{\frak d}$ {\em -invariant}.
\end{itemize}

In a more general way, let $r\in End{\frak g}$ be a classical r-matrix on a
Lie algebra ${\frak g}$ which is equipped with a fixed inner product.
Identify ${\frak g}_r$ with the dual of ${\frak g}$ by means of the inner
product. Assume that $r$ is skew and satisfies the modified Yang-Baxter
identity. Then $\left( {\frak g},{\frak g}^{*}\right) $ is a Lie bialgebra;
it is usually called{\em \ a factorizable Lie bialgebra. }Let us once more
list the properties of factorizable Lie bialgebras for future reference:

\begin{itemize}
\item  ${\frak g}$ {\em is equipped with a nondegenerate inner product.}

\item  ${\frak g}^{*}$ {\em is identified with} ${\frak g}$ {\em as a linear
space by means of the inner product; the Lie bracket in} ${\frak g}^{*}$ 
{\em is given by} 
\[
\left[ X,Y\right] _{*}=\frac 12\left( \left[ rX,Y\right] +\left[ X,rY\right]
\right) ,X,Y\in {\frak g}, 
\]

\item  $r_{\pm }=\frac 12\left( r\pm id\right) $ {\em are Lie algebra
homomorphisms }${\frak g}^{*}\rightarrow {\frak g}$ ; {\em moreover,} $%
r_{+}^{*}=-r_{-},${\em \ and} $r_{+}-r_{-}=id.$
\end{itemize}

Lie bialgebras are particularly important, since the associated Lie-Poisson
structures may be extended to Lie groups (cf. section 3).

\subsubsection{2.1.1. Hamiltonian Reduction.}

Theorem \ref{AKS} admits a very useful global version which survives
quantization. To formulate this theorem we shall first recall some basic
facts about Hamiltonian reduction.

Let $G$ be a Lie group with Lie algebra ${\frak g.}$ Let $T^{*}G$ be the
cotangent bundle of $G$ equipped with the canonical symplectic structure.
Let $B\subset G$ be a Lie subgroup with Lie algebra ${\frak b.}$ The action
of $B$ on $G$ by right translations extends canonically to a Hamiltonian
action $B\times T^{*}G$ $\rightarrow T^{*}G.$ This means that the vector
fields on \ $T^{*}G$ generated by the corresponding infinitesimal action of $%
{\frak b}$ are globally Hamiltonian; moreover, the mapping 
\[
{\frak b}\rightarrow C^\infty \left( T^{*}G\right) :X\longmapsto h_X 
\]
which assigns to each $X\in {\frak b}$ the Hamiltonian of the corresponding
vector field is a homomorphism of Lie algebras. In particular, this map is
linear and hence for each $X\in{\frak g}, x\in T^{*}G$%
\[
h_X(x)=\left\langle X,\mu \left( x\right) \right\rangle , 
\]
where $\mu :T^{*}G\rightarrow {\frak b}^{*}$ is the so-called {\em moment
mapping.}${\frak \ }$Choose a trivialization $T^{*}G{\frak \simeq }G{\frak %
\times g}^{*}$ by means of left translations; then the action of $B$ is
given by 
\[
b:\left( g,\xi \right) \longmapsto \left( \ gb^{-1},\left( Ad^{*}b\right)
^{-1}\xi \right) . 
\]
The corresponding moment map $\mu :T^{*}G\rightarrow {\frak b}^{*}{\frak \ }$%
is 
\[
\mu :\left( g,\xi \right) \longmapsto -\xi \mid _{{\frak b}}, 
\]
where $\xi \mid _{{\frak b}}$ means the restriction of $\xi \in {\frak g}%
^{*} $ to ${\frak b}\subset {\frak g.}$ The moment mapping is {\em %
equivariant}; in other words, the following diagram is commutative
\[
\begin{tabular}{lllll}
$B$ & $\times $ & $T^{*}G$ & $\rightarrow $ & $T^{*}G$ \\ 
&  & $\downarrow $ &  & $\downarrow $ \\ 
$B$ & $\times $ & ${\frak b}^{*}$ & $\rightarrow $ & ${\frak b}^{*}$%
\end{tabular}
, 
\]
where $B\times {\frak b}^{*}\rightarrow {\frak b}^{*}$ is the ordinary
coadjoint action.

Recall that the reduction procedure consists of two steps:

\begin{itemize}
\item  {\em Fix the value of the moment map} $\mu $ {\em and consider the
level surface} 
\[
{\cal M}_F=\left\{ x\in T^{*}G;\mu \left( x\right) =F\right\} ,F\in {\frak b}%
^{*} 
\]
Speaking in physical terms, we impose on $T^{*}G$ linear constraints which
fix the moment $\mu :T^{*}G\rightarrow {\frak b}^{*}.$

\item  {\em Take the quotient of }${\cal M}_F$ {\em over the action of the
stabilizer of }$F$ {\em in} $B.$
\end{itemize}

(Due to the equivariance of $\mu ,$ the stabilizer $B^F\subset B$ preserves $%
{\cal M}_F$ .)

The resulting manifold $\overline{{\cal M}}_F$ is symplectic; it is called
the {\em reduced phase space} obtained by reduction over $F\in {\frak b}%
^{*}. $ A particularly simple case is the reduction over $0\in {\frak b}%
^{*}; $ in that case the stabilizer coincides with \ $B$ itself. We shall
denote the reduced space ${\cal M}_0/B$ by $T^{*}G//B.$

In a slightly more general fashion, we may start with an arbitrary
Hamiltonian $B$-space ${\cal M}$ and perform the reduction over the zero
value of the moment map getting the reduced space ${\cal M}//B.$

Let ${\sf S}_G$ be the category of Hamiltonian $G$-spaces; recall that (up
to a covering) homogeneous Hamiltonian $G$-spaces are precisely the
coadjoint orbits of $G.$ We shall define a functor (called {\em symplectic
induction,} cf. \cite{Kazhdan}), 
\[
{\rm Ind}_B^G:{\sf S}_B\rightsquigarrow {\sf S}_G,
\]
which associates to each Hamiltonian $B$-space a Hamiltonian $G$-space.
Namely, we put 
\[
{\rm Ind}_B^G\left( {\cal M}\right) =T^{*}G\times {\cal M}//B, 
\]
where the manifold $T^{*}G\times {\cal M}$ is equipped with the symplectic
structure which is the {\em difference }of symplectic forms on $T^{*}G$ and
on ${\cal M}$ and $B$ acts on $T^{*}G\times {\cal M}$ diagonally (hence the
moment map on $T^{*}G\times {\cal M}$ is the {\em difference} of the moment
maps on $T^{*}G$ and on ${\cal M}$). The structure of the Hamiltonian $G$%
-space is induced on $T^{*}G\times {\cal M}//B$ by the action of $G$ on $%
T^{*}G$ by left translations.

The most simple case of symplectic induction is when the coadjoint orbit
consists of a single point. In representation theory, this corresponds to
 representations induced by 1-dimensional representations of a subgroup.
Let now $F\in {\frak b}^{*}$ be an arbitrary element; let ${\cal O}_F$ be
its coadjoint orbit, ${\cal K}_F=$ ${\rm Ind}_B^G\left( {\cal O}_F\right) .$
It is sometimes helpful to have a different construction of ${\cal K}_F$ (in
much the same way as in representation theory it is helpful to have
different realizations of a given representation). In many case it is
possible to induce ${\cal K}_F\ $starting from a smaller subgroup $L\subset
B $ and its single-point orbit (a similar trick in representation theory
means that we wish to induce a given representation from a 1-dimensional
one). The proper choice of $L$ is the {\em Lagrangian subgroup} $\
L_F\subset B$ subordinate to $F$ (whenever it exists). Recall that a Lie
subalgebra${\frak \ l}_F\subset {\frak b}$ is called a {\em Lagrangian
subalgebra} subordinate to $F\in {\frak b}^{*}$ if $ad^{*}{\frak l}_F\cdot
F\subset T_F{\cal O}_F$ is a Lagrangian subspace of the tangent space to the
coadjoint orbit of $F\ $(with respect to the Kirillov form on $T_F{\cal O}_F$
). This condition implies that $F\mid _{\left[ {\frak l}_F,{\frak l}%
_F\right] }=0;$ in other words, $F$ defines a character of ${\frak l}_F%
{\frak .}$ This is of course tantamount to saying that $F$ is a single-point
orbit of ${\frak l}_F.$ To define the corresponding Lagrangian subgroup $%
L_F\subset B$ let us observe that ${\frak l}_F$ always coincides with its
own normalizer in ${\frak b.}$ Let us define $L_F\subset B$ as the
normalizer of ${\frak l}_F$ in $B$ (this definition fixes the group of
components of $L_F$ which may be not necessarily connected).

\begin{proposition}
\label{induction}Let $L_F\subset B$ be a Lagrangian subgroup subordinate to $%
F;$ then $T^{*}G\times {\cal O}_F//B\simeq T^{*}G\times \left\{ F\right\}
//L_F.$
\end{proposition}

In the next section we shall discuss an example of a coadjoint orbit with
real polarization.

Let us now return to the setting of theorem \ref{AKS}. Consider the action $%
G\times G\times G\rightarrow G$ by left and right translations: 
\begin{equation}
\left( h,h^{\prime }\right) :x\longmapsto hxh^{\prime -1};
\end{equation}
this action may be lifted to the Hamiltonian action $G\times G\times
T^{*}G\rightarrow T^{*}G.$  Let us choose again a trivialization $T^{*}G\simeq
G\times {\frak g}^{*}$ by means of left translations; then 
\begin{equation}
\left( h,h^{\prime }\right) :\left( g,\xi \right) \longmapsto \left(
hgh^{\prime -1},\left( Ad^{*}h^{\prime }\right) ^{-1}\xi \right) .
\label{left}
\end{equation}
The corresponding moment map $\mu :T^{*}G\rightarrow {\frak g}^{*}{\frak %
\oplus g}^{*}$ is 
\begin{equation}
\mu :\left( g,\xi \right) \longmapsto \left( \xi ,-\left( Ad^{*}g\right)
^{-1}\xi \right) .  \label{moment}
\end{equation}
Using the left trivialization of $T^{*}G$ we may extend polynomial functions 
$\varphi \in $ $P\left( {\frak g}^{*}\right) \simeq {\cal S}\left( {\frak g}%
\right) $ to left-invariant functions on $T^{*}G$ which are polynomial on
the fiber. The Casimir elements ${\cal H}\in {\cal S}\left( {\frak g}\right)
^{{\frak g}}$ give rise to {\em bi-invariant} functions; thus the
corresponding Hamiltonian systems admit reduction with respect to any
subgroup $S\subset $ $G\times G$ .

\begin{lemma}
Let $S\subset G\times G$ be a Lie subgroup, ${\frak s}\subset {\frak g}%
\oplus {\frak g}$ its Lie subalgebra, $p:{\frak g}^{*}\oplus {\frak g}^{*}$ $%
\rightarrow {\frak s}^{*}$ the canonical projection, $S\times
T^{*}G\rightarrow T^{*}G$ the restriction to $H$ of the action $\left( \ref
{left}\right) .$ This action is Hamiltonian and the corresponding moment map
is $\mu _S=p\circ \mu .$
\end{lemma}

Now let us assume again that $\left( {\frak g},r\right) $ is a Lie algebra
equipped with a classical r-matrix $r\in End\ {\frak g}$ satisfying $\left( 
\ref{cybe}\right) .$ Combining the Lie algebra homomorphisms $r_{\pm }:%
{\frak g}_r\rightarrow {\frak g}$ we get an\ embedding $i_r:{\frak g}%
_r\rightarrow {\frak g}\oplus {\frak g};$ we may identify the Lie group
which corresponds to ${\frak g}_r$ with the subgroup $G_r\subset G\times G$
which corresponds to the Lie subalgebra ${\frak g}_r\subset {\frak g}\oplus 
{\frak g}.$ Fix $F\in {\frak g}_r^{*}$ and let ${\cal O}_F\subset {\frak g}%
_r^{*}$ be its coadjoint orbit. We want to apply to ${\cal O}_F$ the
`symplectic induction' procedure. Since half of $G_r$ is acting by left
translations and the other half by right translations, the action of $G$ on
the reduced space will be destroyed. However, the Casimir elements survive
reduction.

\begin{theorem}
\label{red}(i) The reduced symplectic manifold $T^{*}G\times {\cal O}_F//G_r$
is canonically isomorphic to ${\cal O}_F.$ (ii) Extend Casimir elements $%
{\cal H}\in {\cal S}\left( {\frak g}\right) ^{{\frak g}}\simeq P\left( 
{\frak g}^{*}\right) $ to $G\times G$ -invariant functions on $T^{*}G;$
these functions admit reduction with respect to $G_r;$ the reduced
Hamiltonians coincide with those described in theorem \ref{AKS}, namely, $%
{\cal H}_{red}=i_r\left( {\cal H}\right) \mid _{{\cal O}_F};$ the quotient
Hamiltonian flow on ${\cal O}_F$ which corresponds to ${\cal H}\in {\cal S}%
\left( {\frak g}\right) ^{{\frak g}}$ is described by the Lax equation.
\end{theorem}

Theorem \ref{red} allows to get explicit formulae for the solutions of Lax
equations in terms of the factorization problem in $G.$ In these lectures we
are mainly interested in its role for quantization.

The construction of Lax equations based on theorem \ref{AKS} may be applied
to finite-dimensional semisimple Lie algebras (\cite{kost}). However, its
really important applications are connected with loop algebras, which
possess sufficiently many Casimirs. We shall describe two examples both of
which have interesting quantum counterparts: the {\em open Toda lattice}
(this is the example first treated in (\cite{kost})) and the so-called {\em %
generalized Gaudin model }(its quantum counterpart was originally proposed
in (\cite{Gaud})).

\subsection{ Two Examples}

\subsubsection{\label{ClassToda} 2.2.1. Generalized Toda Lattice.}

The Toda Hamiltonian describes transverse oscillations of a 1-dimensional
cyclic molecule with the nearest-neighbour interaction around the
equilibrium configuration. Thus the phase space is ${\cal M}={\Bbb R}^{2n}$
with the canonical Poisson bracket; the Hamiltonian is 
\[
{\cal H}=\frac 12\sum_{i=1}^np_i^2+\sum_{i=1}^{n-1}\exp \left(
q_i-q_{i+1}\right) +\exp \left( q_n-q_1\right) . 
\]
Removing the last term in the potential energy gives the {\em open Toda
lattice} whose behaviour is qualitatively different (all potentials are
repulsive, and the Hamiltonian describes the scattering of particles in the
conical valley 
\[
C_{+}=\left\{ q=\left( q_1,...,q_n\right) \in {\Bbb R}^n;q_i-q_{i+1}\leq
0,i=1,...,n-1\right\} 
\]
with steep walls). The {\em generalized open Toda lattice} may be associated
to an arbitrary semi\-simple Lie algebra. (Generalized {\em periodic} Toda
lattices correspond to affine Lie algebras.) I shall briefly recall the
corresponding construction (cf. \cite{kost}, \cite{GW}, \cite{AreyS}).

Let ${\frak g}$ be a real split semisimple Lie algebra, $\sigma :{\frak g}%
\rightarrow {\frak g}$ a Cartan involution, ${\frak g}={\frak k}\oplus 
{\frak p}$ the corresponding Cartan decomposition (i.e., $\sigma =id$ on $%
{\frak k},$ $\sigma =-id$ on ${\frak p}$). Fix a split Cartan subalgebra $%
{\frak a}\subset {\frak p}$ ; let $\Delta \subset {\frak a}^{*}$ be the root
system of $\left( {\frak g,a}\right) .$ For $\alpha \in \Delta $ let 
\begin{equation}
{\frak g}_\alpha =\left\{ X\in {\frak g};adH\cdot X=\alpha \left( H\right)
X,H\in {\frak a}\right\}  \label{root}
\end{equation}
be the corresponding root space. Fix an order in the system of roots and let 
$\Delta _{+}\subset \Delta $ be the set of positive roots, $P\subset $ $%
\Delta _{+}$ the corresponding set of simple roots. Put ${\frak n}=\oplus
_{\alpha \in \Delta _{+}}{\frak g}_\alpha .$ Let ${\frak b}={\frak a}+{\frak %
n}$ ; recall that ${\frak b}$ is a maximal solvable subalgebra in ${\frak g}$
(Borel subalgebra). We have 
\begin{equation}
{\frak g}={\frak k}\dot +{\frak a}\dot +{\frak n}  \label{iwasa}
\end{equation}
(the Iwasawa decomposition). Equip ${\frak g}$ with the standard inner
product (the Killing form). Clearly, we have ${\frak k}^{\perp }={\frak p,b}%
^{\perp }={\frak n;}$ thus the dual of ${\frak b}$ is modelled on ${\frak p}$,
and we may regard ${\frak p}$ as a ${\frak b}$ -module with respect to the
coadjoint representation. Let $G$ be a connected split semisimple Lie group
with finite center which corresponds to ${\frak g}$, $G=ANK$ its Iwasawa
decomposition which corresponds to the decomposition $\left( \ref{iwasa}%
\right) $ of its Lie algebra.

\begin{remark}
The Iwasawa decomposition usually does {\em not} give rise to a Lie
bialgebra structure on ${\frak g};$ indeed, the subalgebras ${\frak k}$ and $%
{\frak b}$ are not isotropic (unless ${\frak g}$ is a complex Lie algebra
considered as an algebra over ${\Bbb R}$). In other words, the classical
r-matrix associated with the decomposition (\ref{iwasa}) is not skew. We
shall briefly recall the construction of the so-called standard
skew-symmetric classical r-matrix on a semisimple Lie algebra in Section
2.2.3 below.
\end{remark}

The Borel subalgebra admits a decreasing filtration ${\frak b}\supset {\frak %
b}^{(1)}\supset {\frak b}^{(2)}\supset ...,$ where ${\frak b}^{(1)}=\left[ 
{\frak b,b}\right] ={\frak n},{\frak b}^{(2)}=\left[ {\frak n},{\frak n}%
\right] ,...,{\frak b}^{(k)}=\left[ {\frak n,b}^{(k-1)}\right] .$ By
duality, there is an increasing filtration of ${\frak p}\simeq {\frak b}^{*}$
by $ad^{*}{\frak b}$ -invariant subspaces. To describe it let us introduce
some notation. For $\alpha \in \Delta $ let 
\begin{equation}
\alpha =\sum_{\alpha _i\in P}m_i\alpha _i  \label{ht}
\end{equation}
be its decomposition with respect to the simple roots; put $d\left( \alpha
\right) =\sum_{\alpha _i\in P}m_i.$ Put 
\begin{equation}
{\frak d}_p=\bigoplus_{\left\{ \alpha \in \Delta ;d\left( \alpha \right)
=p\right\} }{\frak g}_\alpha ,p\neq 0,{\frak d}_0={\frak a,n}^{\left(
k\right) }=\bigoplus_{p\geq k}{\frak d}_p.  \label{diag}
\end{equation}
Define the mapping $s:$ ${\frak g}\rightarrow {\frak p}:X\longmapsto \frac
12\left( id-\sigma \right) X$, and let ${\frak p}_k=\oplus _{0\leq p\leq
k}s\left( {\frak d}_p\right) ,k=0,1,...$ It is easy to see that ${\frak a}=%
{\frak p}_0\subset {\frak p}_1\subset ...$ are invariant with respect to the
coadjoint action of ${\frak b;}$ moreover, ${\frak p}_k$ is the annihilator
of ${\frak b}^{\left( k+1\right) }$ in ${\frak p}.$ Thus ${\frak p}_k$ is
set in duality with ${\frak b}/{\frak b}^{\left( k+1\right) }.$ Put ${\frak %
s=b}/{\frak b}^{\left( 2\right) }$ and let $L$ $\in {\frak p}_1\otimes 
{\frak s\subset p}_1\otimes {\cal S}\left( {\frak s}\right) $ be the
corresponding canonical element. Let $I_{{\frak s}}$ be the specialization
to ${\frak s}$ of the subalgebra $I={\cal S}\left( {\frak g}\right) ^{{\frak %
g}}$ of Casimir elements. An example of a Hamiltonian $H\in I_{{\frak s}}$
is constructed from the Killing form: 
\begin{equation}
H=\frac 12\left( L,L\right) ,  \label{square}
\end{equation}
in obvious notation. To describe all others we may use the following trick.
Fix a faithful representation $\left( \rho ,V\right) $ of ${\frak g}$ and
put $L_V=\left( \rho \otimes id\right) L.$ ($\rho ,V)$ is usually called the%
{\em \ auxiliary linear representation. }

\begin{proposition}
\label{aux}(i) $I_{{\frak s}}$ is generated by $H_k=tr_VL_V^k\in {\cal S}%
\left( {\frak s}\right) ,k=1,2,...$ (ii) $L$ defines a Lax representation
for any Hamiltonian equation on ${\cal S}\left( {\frak s}\right) $ with the
Hamiltonian $H\in I_{{\frak s}}.$
\end{proposition}

We have ${\frak s}=$ ${\frak a}\ltimes {\frak u,}$ where ${\frak u}={\frak n/%
}\left[ {\frak n},{\frak n}\right] ;$ the corresponding Lie group is $%
S=A\ltimes U,$ $U=N/N^{\prime }.$ Choose $e_\alpha \in {\frak g}_\alpha
,\alpha \in P,$ in such a way that $\left( e_\alpha ,e_{-\alpha }\right) =1,$
and put 
\begin{equation}
{\cal O}_T=\left\{ p+\sum_{\alpha \in P}b_\alpha \left( e_\alpha +e_{-\alpha
}\right) ;p\in {\frak a,\;}b_\alpha >0\right\} ;  \label{orbtoda}
\end{equation}
by construction, ${\cal O}_T\subset {\frak p}_1.$ It is convenient to
introduce the para\-metri\-za\-tion $b_\alpha =\exp \alpha \left( q\right) ,$
$q\in {\frak a.}$

\begin{proposition}
(i) ${\cal O}_T$ is an open coadjoint orbit of $S$ in ${\frak p}_1\simeq 
{\frak s}^{*}$ (ii) Restriction to ${\cal O}_T$ of the Hamiltonian $\left( 
\ref{square}\right) $ is the generalized Toda Hamiltonian 
\[
H_T=\frac 12\left( p,p\right) +\sum_{\alpha \in P}\exp 2\alpha \left(
q\right) .
\]
\end{proposition}

The Toda orbit has the natural structure of a polarized symplectic manifold,
that is, ${\cal O}_T$ admits an $S$-invariant fibering whose fibers are
Lagrangian submanifolds; the description of this Lagrangian structure is a
standard part of the geometric quantization program (\cite{Kirillov}). More
precisely, put $f=\sum_{\alpha \in P}\left( e_\alpha +e_{-\alpha }\right) $;
let us choose $f$ as a marked point on ${\cal O}_T.$

\begin{proposition}
\label{lagr} (i) We have $S/U$ $\simeq A;$ the isomorphism ${\cal O}%
_T=S\cdot f$ induces a Lagrangian fibering ${\cal O}_T\rightarrow A$ with
fiber $U.$ (ii) ${\cal O}_T$ is isomorphic to $T^{*}A$ as a polarized
symplectic manifold. (iii) The Lagrangian subalgebra ${\frak l}_f\subset 
{\frak s}$ subordinate to $f$ is ${\frak l}_f={\frak u.}$
\end{proposition}

\begin{remark}
We may of course regard ${\cal O}_T$ as a coadjoint orbit of the group $B=AN$
itself; indeed, there is an obvious projection $B\rightarrow S$ and its
kernel lies in the stabilizer of $f.$ The 'big' Lagrangian subgroup $\hat
L_f\subset B$ subordinate to $f$ coincides with $N.$
\end{remark}

An alternative description of the open Toda lattice which provides   more
detailed information on its behavior and survives quantization is based on
Hamiltonian reduction. The following result is a version of theorem \ref{red}%
. Let us consider again the Hamiltonian action of $G\times G$ on $T^{*}G$.
To get the generalized Toda lattice we shall restrict this action to the
subgroup $G_r=K\times B\subset G\times G.$ We may regard ${\cal O}_T$ as a $%
G_r$-orbit, the action of $K$ being trivial.

\begin{theorem}
\label{todared} (i) $T^{*}G\times {\cal O}_T//$ $G_r$ is isomorphic to $%
{\cal O}_T$ as a symplectic manifold. (ii) Let ${\cal H}_2$ be the quadratic
Casimir element which corresponds to the Killing form on ${\frak g.}$ The
Toda flow is the reduction of the Hamiltonian flow on $T^{*}G$ generated by $%
{\cal H}_2.$
\end{theorem}

Theorem \ref{todared} leads to explicit formulae for the trajectories of the
Toda lattice in terms of the Iwasawa decomposition of matrices (\cite{OP}, 
\cite{GW}, \cite{AreyS}).

Since ${\cal O}_T$ admits a real polarization, we may state the following
version of proposition \ref{induction}. Observe that if we regard ${\cal O}%
_T $ as a $G_r$-orbit, the Lagrangian subgroup subordinate to $f\in {\cal O}%
_T$ is $K\times N.$

\begin{theorem}
\label{qtoda2red} (i) $T^{*}G\times $ $\left\{ f\right\} //K\times N$ is
isomorphic to ${\cal O}_T$ as a symplectic manifold. (ii) The standard
Lagrangian polarization on $T^{*}G{\cal \ }$by the fibers of projection $%
T^{*}G\rightarrow G$ gives rise to the Lagrangian polarization on ${\cal O}_T
$ described in proposition \ref{lagr}.
\end{theorem}

It is easy to see that the Toda flow is again reproduced as the reduction of
the Hamiltonian flow on $T^{*}G$ generated by the quadratic Hamiltonian $%
{\cal H}_2.$

As discussed in Section 4.1.2, theorem \ref{qtoda2red} gives a very simple
hint to the solution of the quantization problem.

\subsubsection{2.2.2. More examples: standard r-matrices on semisimple Lie
algebras and on loop algebras.}

As already noted, the classical r-matrix associated with the Iwasawa
decomposition of a simple Lie algebra ${\frak g}$ is not skew and hence does
not give rise to a Lie bialgebra structure on ${\frak g.}$ While this is not
a disadvantage for the study of the Toda lattice, the bialgebra structure is
of course basic in the study of the q-deformed case. Let us briefly recall
the so-called {\em standard} Lie bialgebra structure on ${\frak g.}$ We
shall explicitly describe the corresponding Manin triple.

Let ${\frak g}$ be a complex simple Lie algebra, ${\frak a}\subset {\frak g}$
its Cartan subalgebra, $\Delta \subset {\frak a}^{*}$ the root system of $%
\left( {\frak g,a}\right) ,{\frak b}_{+}={\frak a}+{\frak n}_{+}$ a positive
Borel subalgebra which corresponds to some choice of order in $\Delta ,$ $%
{\frak b}_{-}$ the opposite Borel subalgebra. For $\alpha \in \Delta $ let $%
e_\alpha $ be the corresponding root space vector normalized in such a way
that $\left( e_\alpha ,e_{-\alpha }\right) =1.$ We may identify ${\frak a}$
with the quotient algebra ${\frak b}_{\pm }/{\frak n}_{\pm };$ let $\pi
_{\pm }:{\frak b}_{\pm }\rightarrow {\frak a}$ be the corresponding
canonical projection. Put ${\frak d}={\frak g}\oplus {\frak g}$ (direct sum
of two copies); we equip ${\frak d}$ with the inner product which is the 
{\em difference} of the Killing forms on the first and the second copy. Let $%
{\frak g}^\delta \subset {\frak d}$ be the diagonal subalgebra, 
\[
{\frak g}^{*}=\left\{ \left( X_{+},X_{-}\right) \in {\frak b}_{+}\oplus 
{\frak b}_{-};\pi _{+}\left( X_{+}\right) =-\pi _{-}\left( X_{-}\right)
\right\} . 
\]
(Mind the important minus sign in this definition!)

\begin{proposition}
(i) $\left( {\frak d},{\frak g}^\delta ,{\frak g}^{*}\right) $is a Manin
triple. (ii) The corresponding Lie bialgebra structure on ${\frak g}$ is
associated with the classical r-matrix 
\begin{equation}
r=\sum_{\alpha \in \Delta _{+}}e_\alpha \wedge e_{-\alpha }.  \label{r-stand}
\end{equation}
(iii) The Lie bialgebra $\left( {\frak g}^\delta ,{\frak g}^{*}\right) $ is
factorizable.
\end{proposition}

The construction described above admits a straightforward generalization for
loop algebras. The associated r-matrix is called {\em trigonometric}. For
future reference we shall recall its definition as well.

Let  ${\frak g}$ be a complex semisimple Lie algebra, let $L{\frak g}=%
{\frak g}\otimes {\Bbb C}\left[ z,z^{-1}\right] $ be the loop algebra of $%
{\frak g}$ consisting of rational functions with values in ${\frak g}$ which
are regular on ${\Bbb C}P(1)\smallsetminus \left\{ 0,\infty \right\} .$ Let $%
L{\frak g}_0,L{\frak g}_\infty $ be its local completions at 0 and $\infty $;
 by definition, $L{\frak g}_0,L{\frak g}_\infty $ consist of
formal Laurent series in local parameters $z,$ $z^{-1},$ respectively. Put $%
\widehat{{\frak d}}=L{\frak g}_0\oplus L{\frak g}_\infty $ (direct sum of
Lie algebras). Clearly, the diagonal embedding $L{\frak g}\hookrightarrow 
\widehat{{\frak d}}:X\longmapsto \left( X,X\right) $ is a homomorphism of
Lie algebras. Let 
\[
\begin{array}{c}
L{\frak g}_{\pm }={\frak g}\otimes {\Bbb C}[\left[ z^{\pm 1}\right] ], \\ 
\widehat{{\frak b}}_{+}=\left\{ X\in L{\frak g}_{+};X\left( 0\right) \in 
{\frak b}_{+}\right\} , \\ 
\widehat{{\frak b}}_{-}=\left\{ X\in L{\frak g}_{-};X\left( \infty \right)
\in {\frak b}_{-}\right\}
\end{array}
\]
(where $X\left( 0\right) ,X\left( \infty \right) $ denote the constant term
of the formal series). Combining the 'evaluation at 0' (respectively, at $%
\infty )$ with the projections $\pi _{\pm }:{\frak b}_{\pm }\rightarrow 
{\frak a}$ we get two canonical projection maps $\widehat{\pi }_{\pm }:%
\widehat{{\frak b}}_{\pm }\rightarrow {\frak a.}$ Put 
\[
\left( L{\frak g}\right) ^{*}=\left\{ \left( X_{+},X_{-}\right) \in \widehat{%
{\frak b}}_{+}\oplus \widehat{{\frak b}}_{-};\widehat{\pi }_{+}\left(
X_{+}\right) =-\widehat{\pi }_{-}\left( X_{-}\right) \right\} \subset 
\widehat{{\frak d}}. 
\]
Let us set $\left( L{\frak g}\right) ^{*}$ and $L{\frak g}$ in duality by
means of the bilinear form on $\widehat{{\frak d}}$%
\begin{equation}
\begin{array}{l}
\left\langle \left( X_{+},X_{-}\right) ,\left( Y_{+},Y_{-}\right)
\right\rangle = \\ 
\qquad Res_{z=0}\left( X_{+}(z),Y_{+}(z)\right) dz/z+Res_{z=\infty }\left(
X_{-}(z),Y_{-}(z)\right) dz/z.
\end{array}
\label{dzoverz}
\end{equation}

\begin{proposition}
$\left( \widehat{{\frak d}},L{\frak g},\left( L{\frak g}\right) ^{*}\right) $%
is a Manin triple.
\end{proposition}

\begin{remark}
\hfill Properly \hfill speaking, \hfill with \hfill our \hfill choice \hfill %
of \hfill a \hfill completion\hfill for\hfill $\left( L{\frak g}\right) ^{*}$
\hfill the \hfill Lie \\ bialgebra\hfill $\left( L{\frak g},\left( L{\frak g}%
\right) ^{*}\right) $\hfill is \hfill{\em not} \hfill factorizable;\hfill
indeed, \hfill our \hfill definition\hfill requires\hfill that\hfill the%
\hfill Lie\hfill algebra\\ and \hfill its \hfill dual\hfill should \hfill be %
\hfill isomorphic \hfill as \hfill linear\hfill spaces.\hfill However, %
\hfill a\hfill slightly\hfill weaker \hfill assertion\\ still\hfill holds %
\hfill true: \hfill the\hfill dual \hfill Lie \hfill algebra \hfill$\left( L%
{\frak g}\right) ^{*}$ \hfill contains\hfill an \hfill{\em open \hfill dense%
\hfill subalgebra} \hfill$\left( L{\frak g}\right) ^{\circ }=$\\ $\left\{
\left( X_{+},X_{-}\right) \in \left( L{\frak g}\right) ^{*};X_{\pm }\in L%
{\frak g}\right\} $ such that the mappings $\left( L{\frak g}\right) ^{\circ
}\rightarrow L{\frak g}:\left( X_{+},X_{-}\right) \longmapsto X_{\pm }$ are
Lie algebra homomorphisms and $\left( L{\frak g}\right) ^{\circ }$ is
isomorphic to $L{\frak g}$ as a linear space. This subtlety is quite typical
of infinite-dimensional Lie algebras.
\end{remark}

Let us finally say a few words on the trigonometric r-matrix itself. Define
the cobracket $\delta :L{\frak g\rightarrow }L{\frak g}\otimes L{\frak g}$
by 
\begin{equation}
\left\langle \delta X,Y\otimes Z\right\rangle =\left\langle X,\left[
Y,Z\right] _{L{\frak g}^{*}}\right\rangle ,X\in L{\frak g},Y,Z\in \left( L%
{\frak g}\right) ^{*},  \label{trig}
\end{equation}
i.e. as the dual of the Lie bracket in $\left( L{\frak g}\right) ^{*}$ with
respect to the pairing (\ref{dzoverz}). The element $\delta (X)\in L{\frak g}%
\otimes L{\frak g}$ may be regarded as a Laurent polynomial in two variables
with values in ${\frak g};$ the cobracket (\ref{trig}) is then given by 
\[
\delta (X)\left( z,w\right) =\left[ r_{trig}\left( \frac zw\right)
,X_1\left( z\right) +X_2\left( w\right) ,\right] 
\]
where $r_{trig}\left( z/w\right) $ is a rational function in $z/w$ with
values in ${\frak g}\otimes {\frak g.}$ (We use the dummy indices $1,2,...$
to denote different copies of linear spaces; in other words, 
\[
X_1\left( z\right) :=X\left( z\right) \otimes 1,X_2\left( w\right)
:=1\otimes X\left( w\right) ; 
\]
below we shall frequently use this abridged tensor notation.) The explicit
expression for $r_{trig}$ is given by 
\begin{equation}
r_{trig}\left( x\right) =t\cdot \frac{x+1}{x-1}+r,  \label{rtrig}
\end{equation}
where $t\in {\frak g}\otimes {\frak g}$ is the canonical element ({\em %
tensor Casimir}) which represents the inner product $\left\langle
,\right\rangle $ in ${\frak g}$ and $r\in {\frak g}\wedge {\frak g}$ is the
standard classical r-matrix in ${\frak g}$ given by (\ref{r-stand}).

\begin{remark}
Note that the expression for $r_{trig}$ which may be found in \cite{beldr}
is different from (\ref{rtrig}), since these authors use a different grading
in the loop algebra (the so-called principal grading, as compared to the
standard grading which we use in these lectures.)
\end{remark}

\subsubsection{\label{ClassGaudin}2.2.3. Generalized Gaudin Model.}

We shall return to the bialgebras described in the previous section when we
shall discuss the q-deformed algebras and the quadratic Poisson Lie groups.
In order to deal with our next example of an integrable system we shall use
another important Lie bialgebra associated with the Lie algebra of rational
functions. The example to be discussed is a model of spin-spin interaction
(the so-called {\em Gaudin model}) which may be deduced as a limiting case
of lattice spin models (such as the Heisenberg XXX model) for small values
of the coupling constant (\cite{Gaud}). The original Gaudin model is related
to the ${\frak sl}\left( 2\right) $ or ${\frak su}\left( 2\right) $
algebras; below we shall discuss its generalization to arbitrary simple Lie
algebras. The generalized Gaudin Hamiltonians also arise naturally in the
semiclassical approximation to the Knizhnik-Zamolodchikov equations (\cite
{RV}). The underlying r-matrix associated with the Gaudin model is the so
called {\em rational r-matrix} which is described below.

Let ${\frak g}$ be a complex simple Lie algebra. Fix a finite set $D=\left\{
z_1,...,z_N\right\} \subset {\Bbb C}\subset {\Bbb C}P_1$ and let ${\frak g}%
(D)$ be the algebra of rational functions on ${\Bbb C}P_1$ with values in%
{\bf \ }${\frak g}$ which vanish at infinity and are regular outside $D$.
Let ${\frak g}_{z_i}={\frak g}\otimes {\Bbb C}((z-z_i))$ be the {\em %
localization} of ${\frak g}(D)$ at $\ z_i\in D.$ [As usual, we denote by $%
{\Bbb C}((z))$ the algebra of formal Laurent series in $z,$ and by ${\Bbb C}%
\left[ \left[ z\right] \right] $ its subalgebra consisting of formal Taylor
series.] Put ${\frak g}_D=\oplus _{z_i\in D}{\frak g}_{z_i}.$ There is a
natural embedding ${\frak g}(D)\hookrightarrow {\frak g}_D$ which assigns to
a rational function the set of its Laurent expansions at each point $z_i\in
D $. Put ${\frak g}_{z_i}^{+}={\frak g}\otimes {\Bbb C}\left[ \left[
z-z_i\right] \right] ,\ {\frak g}_D^{+}=\oplus _{z_i\in D}\;{\frak g}_{z_i}.$
Then 
\begin{equation}
{\frak g}_D={\frak g}_D^{+}\ \dot +\ {\frak g}(D)  \label{split}
\end{equation}
as a linear space. [This assertion has a very simple meaning. Fix an element 
$X=\left( X_i\right) _{z_i\in D}\in {\frak g}_D;$ truncating the formal
series $X_i$, we get a finite set of Laurent polynomials which may be
regarded as principal parts of a rational function. Let $X^0$ be the
rational function defined by these principal parts; it is unique up to a
normalization constant which may be fixed by a condition at infinity. Expand 
$X^0$ in its Laurent series at $z=z_i;$ by construction, $X_i-X_i^0\in 
{\frak g}_{z_i}^{+}$ ; thus $P:X\longmapsto X_0$ is a projection operator
from ${\frak g}_D$ onto ${\frak g}(D)$ parallel to ${\frak g}_D^{+}.$ In
other words, the direct sum decomposition (\ref{split}) is equivalent to the
existence of a rational function with prescribed principal parts; this is
the assertion of the well known Mittag-Leffler theorem in complex analysis. ]

Fix an inner product on ${\frak g}$ and extend it to ${\frak g}_D$ by
setting 
\begin{equation}
\left\langle X,Y\right\rangle =\sum_{z_i\in D}Res\left( X_i,Y_i\right) dz,\
X=\left( X_i\right) _{z_i\in D},Y=\left( Y_i\right) _{z_i\in D}.  \label{res}
\end{equation}
(One may notice that the key difference from the formula (\ref{dzoverz}) in
the previous section is in the choice of the differential $dz$ instead of $%
dz/z.)$ Both subspaces ${\frak g}(D),{\frak g}_D^{+}$ are isotropic with
respect to the inner product (\ref{res}) which sets them into duality. Thus $%
\left( {\frak g}_D,{\frak g}_D^{+},{\frak g}(D)\right) \ $ is a Manin
triple.\ The linear space ${\frak g}(D)$ may be regarded as a ${\frak g}%
_D^{+}-$module with respect to the coadjoint representation. The action of $%
\ {\frak g}_D^{+}$ preserves the natural ${\Bbb Z}^N$-filtration of ${\frak g%
}(D)$ by the order of poles. In particular, rational functions with simple
poles form an invariant subspace ${\frak g}(D)_1\subset {\frak g}(D).$ Let $%
{\frak g}_{z_i}^{++}{\frak \ }\subset {\frak g}_{z_i}^{+}$ be the subalgebra
consisting of formal series without constant terms. Put ${\frak g}%
_D^{++}=\oplus _{z_i\in D}\;{\frak g}_{z_i}^{++}$; clearly, ${\frak g}%
_D^{++} $ is an ideal in ${\frak g}_D^{+}$ and its action in ${\frak g}(D)_1$
is trivial.\ The quotient algebra ${\frak g}_D^{+}/{\frak g}_D^{++}$ is
isomorphic to ${\frak g}^N\ =\oplus _{z_i\in D}\;{\frak g}$. The inner
product (\ref{res}) sets the linear spaces ${\frak g}(D)_1$ and ${\frak g}^N$ in
duality.\ Let $L(z)\in {\frak g}(D)_1\otimes {\frak g}^N$ be the canonical
element; we shall regard $L(z)$ as a matrix-valued rational function with
coefficients in ${\frak g}^N\subset S({\frak g}^N).$ Fix a faithful linear
representation $(\rho ,V)$ of ${\frak g}$ ({\em the auxiliary linear
representation}); it extends canonically to a representation $(\rho ,V(z))$
of the Lie algebra ${\frak g}(z)={\frak g}\otimes {\Bbb C}(z)$ in the space $%
V(z)=V\otimes {\Bbb C}(z)$. (In section 2.2.1 we have already used this
auxiliary linear representation for a similar purpose, cf. proposition \ref
{aux} .) Put 
\begin{equation}
L_V(z)=(\rho \otimes id)L(z).  \label{Lax}
\end{equation}
The matrix coefficients of $L_V(z)$ generate the algebra of observables $%
{\cal S}({\frak g}^N).$ The Poisson bracket relations in this algebra have a
nice expression in 'tensor form', the brackets of the matrix coefficients of 
$L_V\ \ $ forming a matrix in $EndV\otimes EndV$ with coefficients in $\ \ 
{\Bbb C}(u,v)\otimes {\cal S}({\frak g}^N);$ the corresponding formula,
suggested for the first time by \cite{rmat}, was the starting point of the
whole theory of classical r-matrices. To describe it let us first introduce
the {\em rational r-matrix,} 
\begin{equation}
r_V(u,v)=\frac{t_V}{u-v},\;  \label{rat}
\end{equation}
Here $t$ is the {\em tensor Casimir} of ${\frak g}$ which corresponds to the
inner product in ${\frak g}$ (i.e., $t=\sum e_a\otimes e_a,$ where$\ \left\{
e_a\right\} $ is an orthonormal basis in ${\frak g}$) and $t_V=\rho
_V\otimes \rho _V\left( t\right) .$ Notice that $r_V(u,v)$ is essentially
the Cauchy kernel solving the Mittag-Leffler problem on ${\Bbb C}P_1$ with
which we started.

\begin{proposition}
The matrix coefficients of $L_V(u)$ satisfy the Poisson bracket relations 
\begin{equation}
\left\{ L_V(u)\otimes _{,}L_V(v)\right\} =\left[ r_V(u,v),L_V(u)\otimes
1+1\otimes L_V(v)\right] .  \label{pbr}
\end{equation}
\end{proposition}

We shall now specialize corollary \ref{CorAKS} to the present setting; it
shows that in some way spectral invariants of the Lax operator $L(z)$ may be
regarded as 'radial parts' of the Casimir elements.

\begin{remark}
One should use some caution, as Casimir elements do not lie in the
symmetric algebra ${\cal S}({\frak g}_D)$ itself but rather in its
appropriate local completion; however, their projections to ${\cal S}({\frak %
g}^N)$ are well defined. More precisely, for each $z_i\in D$ consider the
projective limit
\[
\tilde {{\cal S}}\left( {\frak g}_{z_i}\right) =\stackunder{n\rightarrow
+\infty }{\stackrel{\longleftarrow }{\lim }}{\cal S}\left( {\frak g}%
_{z_i}\right) /{\cal S}\left( {\frak g}_{z_i}\right) \left( {\frak g}\otimes
z_i^n{\Bbb C}\left[ \left[ z_i\right] \right] \right) 
\]
and put $\tilde {{\cal S}}({\frak g}_D)=\otimes _i\tilde {{\cal S}}\left( 
{\frak g}_{z_i}\right) .$ Passing to the dual language let us observe that
the inner product on ${\frak g}_D$ induces the evaluation map ${\frak g}%
(D)\times {\cal S}({\frak g}_D)\rightarrow {\Bbb C}$ ; the completion is
chosen in such a way that this map makes sense for $\tilde {{\cal S}}({\frak %
g}_D).$ To relate Casimir elements lying in $\tilde {{\cal S}}({\frak g}_D)$
to the more conventional spectral invariants, let us notice that if $L(z)\in 
{\frak g}(D)$ is a 'Lax matrix', its spectral invariants, e.g., the
coefficients of its characteristic polynomial 
\begin{equation}
P\left( z,\lambda \right) \ =\det \left( L\left( z\right) -\lambda \right)
=\sum \sigma _k\left( z\right) \lambda ^k,  \label{det}
\end{equation}
are rational functions in $z.$ To get numerical invariants we may expand $%
\sigma _k\left( z\right) $ in a local parameter (at any point $\zeta \in {\Bbb %
C}P_1$ ) and take any coefficient of this expansion. The resulting
functionals are well defined as polynomial mappings from ${\frak g}(D)$ into 
${\Bbb C}$ . It is easy to see that any such functional is obtained by
applying the evaluation map to an appropriate Casimir element $\zeta \in $ $%
\tilde {{\cal S}}({\frak g}_D)^{{\frak g}_D}.$
\end{remark}

Clearly, we have ${\frak g}_D={\frak g}(D)\ \dot +\ {\frak g}^N\ \dot +\ 
{\frak g}_D^{++}$, and hence 
\[
{\cal S}({\frak g}_D)={\cal S}({\frak g}^N)\oplus ({\frak g}(D)\;{\cal S}(%
{\frak g}_D)+{\cal S}({\frak g}_D)\;{\frak g}_D^{++}). 
\]
Let $\gamma :{\cal S}({\frak g}_D)\rightarrow {\cal S}({\frak g}^N)$ be the
projection map associated with this decomposition.

\begin{proposition}
\label{Harish}The restriction of $\gamma $ to the subalgebra $\widetilde{%
{\cal S}}({\frak g}_D)^{{\frak g}_D}$ of Casimir elements is a morphism of
Poisson algebras; under the natural pairing ${\cal S}({\frak g}^N)\times 
{\frak g}(D)_1\rightarrow {\Bbb C}$ induced by the inner product (\ref{res})
the restricted Casimirs coincide with the spectral invariants of $L(z).$
\end{proposition}

The mapping $\gamma $ is an analogue of the {\em Harish-Chandra homomorphism}
(\cite{dixmier}); its definition may be extended to the quantum case as well.

\begin{corollary}
\label{inv} Spectral invariants of $L(z)$ are in involution with respect to
the standard Lie-Poisson bracket on ${\frak g}^N;L(z)$ defines a Lax
representation for any of these invariants (regarded as a Hamiltonian on ($%
{\frak g}^N)^{*}\simeq {\frak g}^N.$)
\end{corollary}

The generalized Gaudin Hamiltonians are, by definition, the quadratic
Hamiltonians contained in this family; one may take, e.g., 
\begin{equation}
{\cal H}_V(u)=\frac 12tr_V\;L(u)^2.  \label{gaudc}
\end{equation}
(Physically, they describe, e.g., the bilinear interaction of several 'magnetic
momenta'.)

Corollary \ref{inv} does not mention the 'global' Lie algebra ${\frak g}_D$
(and may in fact be proved by elementary means). However, of the three Lie
algebras involved, it is probably the most important one, as it is
responsible for the dynamics of the system. It is this 'global' Lie algebra that may be
called the{\em \ hidden symmetry algebra}${\em .}$

Setting $P\left( z,\lambda \right) =0$ in (\ref{det}), we get an affine
algebraic curve (more precisely, a family of curves parametrized by the
values of commuting Hamiltonians); this switches on the powerful
algebraic-geometric machinery and makes the complete integrability of the
generalized Gaudin Hamiltonians almost immediate (\cite{AreyS}).

\begin{remark}
\ In the exposition above we kept the divisor $D$ fixed. A more invariant
way is of course to use ad\`elic language.\ Thus one may replace the algebra 
${\frak g}_D$ with the global algebra{\bf \ }${\frak g}${\bf $_A$} defined
over the ring{\bf \ A} of ad\`eles of ${\Bbb C}(z);$ fixing a divisor
amounts to fixing a Poisson subspace inside the Lie algebra ${\frak g}(z)$
of rational functions (which is canonically identified with the dual space
of the subalgebra ${\frak g}_{{\bf A}}^{+}).$ The use of ad\`elic freedom
allowing us to add new points to the divisor is essential in the treatment of
the quantum Gaudin model.
\end{remark}

\begin{remark}
Further generalizations consist in replacing the rational r-matrix with more
complicated ones. One is of course tempted to repeat the construction above,
replacing ${\Bbb C}P_1$ with an arbitrary algebraic curve.\ There are
obvious obstructions which come from the Mittag-Leffler theorem for curves:
the subalgebra ${\frak g}(C)$ of rational functions on the curve $C$ does
not admit a complement in ${\frak g}{\bf _A}$ which is again a Lie
subalgebra; for elliptic curves (and ${\frak g}=sl(n))$ this problem may be
solved (\cite{beldr}, \cite{AreyS}) by considering quasiperiodic functions on 
$C,$ and in this way we recover elliptic r-matrices (which were originally
the first example of classical r-matrices ever considered (\cite{rmat}))!
\end{remark}

\section{\label{quadr}Quadratic Case}

The integrable systems which we have considered so far are modelled on Poisson
submanifolds of linear spaces, namely on the dual spaces of appropriate Lie
algebras equipped with the Lie-Poisson bracket. As we shall see,
quantization of such systems (whenever possible) still leaves us in the
realm of ordinary (i.e., non-quantum) Lie groups and algebras. Let us note
in passing that, in the context of the Integrability phenomena, the name `Quantum
Groups' leads to some confusion: one is tempted to believe that quantization
of integrable systems implies that the corresponding `hidden symmetry
groups' also become quantum. As a matter of fact, the real reason to
introduce Quantum groups is {\em different}. (This is clear from the fact
that the Planck constant and the deformation parameter $q$ which enters the
definition of Quantum groups are independent of each other!)  The point is
that we are usually interested not in {\em individual} integrable systems,
but rather in ${\em families}$ of such systems with an arbitrary number of
'particles'. (This is particularly natural for problems arising in Quantum
Statistical Mechanics.) A good approximation is obtained if one assumes that the phase
space of a multiparticle system is the direct product of phase spaces for
single particles. However, the definition of {\em observables} for
multiparticle systems is not straightforward and tacitly assumes the
existence of a Hopf structure (or some of its substitutes) on the algebra $%
{\cal A}$ of observables, i.e. of a map $\Delta ^{\left( N\right) }:{\cal A}%
\rightarrow \otimes ^N{\cal A}$ which embeds the algebra of observables of a
single particle into that of the multiparticle system and thus allows us to
speak of individual particles inside the complex system. For systems which
are modelled on dual spaces of Lie algebras the underlying Hopf structure is
that of ${\cal S}\left( {\frak g}\right) $ or ${\cal U}\left( {\frak g}%
\right) $ which is derived from the additive structure on ${\frak g}^{*}$
(cf. section \ref{gen}). There is an important class of integrable systems (%
{\em difference Lax equations}) for which the natural Hopf structure is
different; it is derived from the multiplicative structure in a nonabelian
Lie group. Let me briefly recall the corresponding construction.

As already mentioned, difference Lax equations arise when the 'auxiliary
linear problem' which underlies the construction is associated with a
difference equation. A typical example is a first order difference system, 
\begin{equation}
\psi _{n+1}=L_n\psi _n,n\in {\Bbb Z}{\bf .}  \label{lax}
\end{equation}
The discrete variable $n\in {\Bbb Z}$ labels the points of an infinite
1-dimensional lattice; if we impose the periodic boundary condition $%
L_n=L_{n+N},$ it is replaced with a finite periodic lattice parametrized by $%
{\Bbb Z}{\bf /}N{\Bbb Z}{\bf .}$ It is natural to assume that the Lax matrices $%
L_n$ belong to a matrix Lie group $G.$ Difference Lax equations arise as
compatibility conditions for the linear system, 
\begin{equation}
\begin{array}{c}
\psi _{n+1}=L_n\psi _n, \\ 
\partial _t\psi _n=A_n\psi _n,
\end{array}
\label{diff}
\end{equation}
with an appropriately chosen $A_n.$ They have the form of {\em finite-difference
 zero-curvature equations}, 
\begin{equation}
\partial _tL_n=L_nA_{n+1}-A_nL_n.  \label{curv}
\end{equation}
Let $\Psi $ be the fundamental solution of the difference system (\ref{lax})
normalized by $\Psi _0=I.$ The value of $\Psi $ at $n=N$ is called the {\em %
monodromy matrix} for the periodic problem and is denoted by $M_L.$ Clearly, 
\[
\begin{array}{c}
\Psi _n=\overleftarrow{\prod }_{0\leq k<n}L_k, \\ 
M_L=\overleftarrow{\prod }_{0\leq k<N}L_k.
\end{array}
\]
Spectral invariants of the difference operator associated with the periodic
difference system (\ref{lax}) are the eigenvalues of the monodromy matrix.
Thus one expects 
\[
I_s=trM_L^s,s=1,2,..., 
\]
to be the integrals of motion for any Lax equation associated with the
linear problem (\ref{lax}). This will hold if the monodromy matrix itself
evolves isospectrally, i.e. , 
\begin{equation}
\partial _tM_L=\left[ M_L,A(M_L)\right]  \label{nov}
\end{equation}
for some matrix $A(T_L).$ Let us now discuss possible choices of the Poisson
structure on the phase space which will turn (\ref{curv}) into Hamiltonian
equations of motion. Observe first of all that the phase space of our system
is the product ${\Bbb G}=\prod^NG;$ if we assume that the variables
corresponding to different copies of $G$ are independent of each other, we
may equip ${\Bbb G}$ with the product Poisson structure. The monodromy may
be regarded as a mapping $M:{\Bbb G}\rightarrow G$ which assigns to a
sequence $(L_0,...,L_{N-1})$ the ordered product $M_L=\overleftarrow{\prod }%
L_k.$ Let ${\bf F}_t:{\Bbb G}\rightarrow {\Bbb G}$ be the (local) dynamical
flow associated with equation (\ref{curv}) and $F_t:G\rightarrow G$ the
corresponding flow associated with equation (\ref{nov}) for the monodromy.
We expect the following diagram to be commutative: 
\[
\begin{array}{ccccc}
& {\Bbb G} & \stackrel{{\bf F}_t}{\rightarrow } & {\Bbb G} &  \\ 
& \mid &  & \mid &  \\ 
^M & \downarrow &  & \downarrow & ^M \\ 
& G & \stackrel{F_t}{\rightarrow } & G & 
\end{array}
. 
\]
In other words, the dynamics in ${\Bbb G}$ factorizes over $G.$ It is
natural to demand that all maps in this diagram should be Poissonian; in
particular, the monodromy map $M:{\Bbb G}\rightarrow G$ should be compatible
with the product structure on ${\Bbb G}=\prod^NG$. By induction, this
property is reduced to the following one:

\vspace{0.5cm}

{\em Multiplication} $m:G\times G\rightarrow G$ {\em is a Poisson mapping.}

\vspace{0.5cm}

This is precisely the axiom introduced by \cite{PGr} as a definition of
Poisson Lie groups. In brief, one can say that the multiplicativity property
of the Poisson bracket on $G$ matches perfectly with the kinematics of
one-dimensional lattice systems.

\subsection{\label{quadr abstr}. Abstract Case: Poisson Lie Groups and
Factorizable Lie Bialgebras}

To construct lattice Lax system one may start with an arbitrary Manin
triple, or, still more generally, with a factorizable Lie bialgebra. Let us
briefly describe this construction which is parallel to the linear case but
involves the geometry of Poisson Lie groups. Although main applications are
connected with infinite-dimensional groups (e.g., loop groups) it is more
instructive to start with the finite-dimensional case.

Let $\left( {\frak g,g}^{*}\right) $ be a finite-dimensional factorizable
Lie bialgebra. Recall that in this case $r_{\pm }=\frac 12\left( r\pm
id\right) $ are Lie algebra homomorphisms from ${\frak g}^{*}\simeq {\frak g}%
_r$ into ${\frak g.}$ Let $G,G^{*}$ be the connected simply connected Lie groups
 which correspond to ${\frak %
g,g}^{*},$ respectively.  Let us assume that $G$ is a linear group; let $\left( \rho
,V\right) $ be its faithful linear representation. We extend $r_{\pm }$ to
Lie group homomorphisms from $G^{*}$ into $G$ which we denote by the same
letters; we shall also write $h_{\pm }=r_{\pm }\left( h\right) ,$ $h\in
G^{*}.$ Both $G$ and $G^{*}$ carry a natural Poisson structure which makes
them Poisson Lie groups. The description of the Poisson bracket on $G$ is
particularly simple. For $\varphi \in F\left( G\right) $ let $D\varphi
,D^{\prime }\varphi \in {\frak g}^{*}$ be its left and right differentials
defined by 
\begin{eqnarray}
&&\left\langle D\varphi \left( L\right) ,X\right\rangle =\frac d{dt}\varphi
\left( \exp tX\cdot L\right) \mid _{t=0}, \\
&&\left\langle D^{\prime }\varphi \left( L\right) ,X\right\rangle =\frac
d{dt}\varphi \left( L\exp tX\right) \mid _{t=0},L\in G,X\in {\frak g}. 
\nonumber
\end{eqnarray}
The{\em \ Sklyanin bracket} on $G$ is defined by 
\begin{equation}
\left\{ \varphi ,\psi \right\} =\left\langle r,D^{\prime }\varphi \wedge
D^{\prime }\psi \right\rangle -\left\langle r,D\varphi \wedge D\psi
\right\rangle .  \label{ASkl}
\end{equation}
An equivalent description which makes sense when $G$ is a linear group: For $L\in G$ put 
$L_V=\rho _V\left( L\right) ;$ the matrix coefficients of $L_V$ generate the
affine ring ${\cal F}\left[ G\right] $ of $G.$ It is convenient not to fix a
representation $\left( \rho _V,V\right) $ but to consider all
finite-dimensional representations of $G$ simultaneously. The reason is that
the Hamiltonians which are used to produce integrable systems are spectral
invariants of the Lax matrix. To get a complete set of invariants one may
consider either $tr_VL_V^n$ for a {\em fixed} representation and all $n$ $%
\geq 1$ or, alternatively, $tr_VL_V$ for {\em all} linear representations.
The second version is more convenient in the quantum case. If $\left( \rho
_W,W\right) $ is another linear representation, we get from the definition $%
\left( \ref{ASkl}\right) :$%
\begin{equation}
\left\{ L_V\stackunder{,}{\otimes }L_W^{\ }\right\} =\left[
r_{VW},L_V\otimes L_W\right] ,\text{ where }r_{VW}=\left( \rho _V\otimes
\rho _W\right) \left( r\right) .  \label{skl}
\end{equation}
To avoid excessive use of tensor product signs one usually puts 
\[
L_V^1=L_V\otimes I,L_W^2=I\otimes L_W; 
\]
formula (\ref{skl}) is then condensed to 
\begin{equation}
\left\{ L_V^1,L_W^2\right\} =\left[ r_{VW}^{},L_V^1L_W^2\right] .
\label{skl12}
\end{equation}
Let $I\subset {\cal F}\left[ G\right] $ be the subalgebra generated by $%
tr_VL_V,$ $V\in {\sf Rep\ }G.$ Formula $\left( \ref{skl12}\right) $
immediately implies the following assertion.

\begin{proposition}
$I\subset {\cal F}\left[ G\right] $ is a commutative subalgebra with respect
to the Poisson bracket $\left( \ref{ASkl}\right) .$
\end{proposition}

The dual Poisson bracket on $G^{*}$ may be described in similar terms.
Extend $r_{\pm }$ to Lie group homomorphisms $r_{\pm }:G^{*}\rightarrow G;$
it is easy to see that $i_r:G^{*}\rightarrow G\times G:T\longmapsto \left(
r_{+}\left( T\right) ,r_{-}\left( T\right) \right) $ is an injection; hence $%
G^{*}$ may be identified with a subgroup in $G\times G.$ We shall write $%
r_{\pm }\left( T\right) =T_{\pm }$ for short. Let again $\left( \rho
_V,V\right) ,\left( \rho _W,W\right) $ be two linear representation of $G;$
put $T_V^{\pm }=\rho _V\circ T_{\pm }\ ,T_W^{\pm }=\rho _W\circ T_{\pm }\ $ $%
T\in G^{*}.$ Matrix coefficients of $T_V^{\pm },$  $\rho _V\in {\sf Rep}G,$
generate the affine ring ${\cal F}\left[ G^{*}\right] $ of $G^{*}.$

The Poisson bracket relations for the matrix coefficients of $T_V^{\pm },$ $%
T_W^{\pm }$ are given by 
\begin{eqnarray}
\left\{ T_V^{\pm 1},T_W^{\pm 2}\right\} &=&\left[ r_{VW},T_V^{\pm 1}T_W^{\pm
2}\right] ,  \label{dual} \\
\left\{ T_V^{+1},T_W^{-2}\right\} &=&\left[
r_{VW}+t_{VW},T_V^{+1}T_W^{-2}\right] ;  \nonumber
\end{eqnarray}
here $t\in {\frak g}\otimes {\frak g}$ is the tensor Casimir which
corresponds to the inner product on ${\frak g}$ and $t_{VW}=\rho _V\otimes
\rho _W\left( t\right) .$

Consider the action 
\begin{equation}
G^{*}\times G\rightarrow G:\left( h,x\right) \longmapsto h_{+}xh_{-}^{-1}.
\label{g act}
\end{equation}

\begin{lemma}
\label{A}$G^{*}$ acts locally freely on an open cell in $G$ containing unit
element.
\end{lemma}

In other words, almost all elements $x\in G$ admit a decomposition $%
x=h_{+}h_{-}^{-1},$ where $h_{\pm }=r_{\pm }\left( h\right) $ for some $h\in
G^{*}.$ (Moreover, this decomposition is unique if both $h$ and $%
(h_{+},h_{-})$ are sufficiently close to the unit element.) This explains
the term 'factorizable group'.

Lemma \ref{A} allows us to push forward the Poisson bracket from $G^{*} $ to $%
G. $ More precisely, we have the following result.

\begin{proposition}
There exists a unique Poisson structure $\left\{ ,\right\} ^{*}$ on $G$ such
that 
\[
F:G^{*}\rightarrow G:\left( x_{+},x_{-}\right) \longmapsto x_{+}x_{-}^{-1}
\]
is a Poisson mapping; it is given by
\begin{equation}
\begin{array}{l}
\left\{ \varphi ,\psi \right\} ^{*}=\left\langle r,D\varphi \wedge D\psi
\right\rangle +\left\langle r,D^{\prime }\varphi \wedge D^{\prime }\psi
\right\rangle  \\ 
\qquad -\left\langle 2r_{+},D\varphi \wedge D^{\prime }\psi \right\rangle
-\left\langle 2r_{-},D^{\prime }\varphi \wedge D\psi \right\rangle
,\;\varphi ,\psi \in C^\infty \left( G\right) .
\end{array}
\label{DDual}
\end{equation}
Equivalent description of the bracket $\left\{ ,\right\} ^{*}$ for linear
groups: let again $L_V=\rho _V\left( L\right) ,L_W=\rho _W\left( L\right) ;$
then
\begin{equation}
\begin{array}{l}
\left\{ L_V^1,L_W^2\right\} ^{*}=  \\
\ r_{VW}L_V^1L_W^2+L_V^1L_W^2r_{VW}  \\
\ -L_V^1r_{VW}L_W^2-L_W^2r_{VW}L_V^1-L_V^1t_{VW}L_W^2+L_W^2t_{VW}L_V^1.
\end{array}
  \label{dual12}
\end{equation}
\end{proposition}

\begin{remark}
(i) A priori the bracket $\left\{ ,\right\} ^{*}$ is defined only on the big
cell in $G$ (i.e., on the image of $G^{*});$ however, formulae (\ref{DDual},%
\ref{dual12}) show that it extends smoothly to the whole $G.$ (ii) One may
notice that the choice of $F$ is actually rigid: this is essentially the
only combination of $x_{+},x_{-}$ such that the bracket for $x=F(x_{+},x_{-})
$ is expressed in terms of $x$ (not of its factors).
\end{remark}

Recall that $I\subset {\cal F}\left[ G\right] $ consists of spectral
invariants; thus we are again in the setting of the generalized
Kostant-Adler theorem: there are two Poisson structures on the same
underlying space and the Casimirs of the former are in involution with
respect to the latter.

\begin{theorem}
Hamiltonians ${\cal H}\in I$ generate Lax equations on $G$ with respect to
the Sklyanin bracket $.$
\end{theorem}

Geometrically, this means that the Hamiltonian flow generated by ${\cal H}%
\in I$ preserves two systems of symplectic leaves in $G,$ namely, the
symplectic leaves of $\left\{ ,\right\} $ and $\left\{ ,\right\} ^{*}.$ The
latter coincide with the {\em conjugacy classes} in $G$ (for a proof see 
jjcite{CIMPA}). In order to include lattice Lax equations into
our general framework we have to introduce the notion of {\em twisting}
which is discussed in the next section.

\begin{remark}
The Hamiltonian reduction picture discussed in section 2.1.1 may be fully
generalized to the present case as well (\cite{CIMPA}). To define the
symplectic induction functor we now need the {\em nonabelian moment map }%
(see \cite{varenna} for a review). In the present exposition we shall not
describe this construction.
\end{remark}

\subsection{Duality theory for Poisson Lie groups and twisted spectral
invariants}

Our key observation so far has been  the duality between the Hamiltonians of
integrable systems on a Lie group and the Casimir functions of its Poisson
dual. In section 4.2 we shall see that a similar relation holds for
quantized universal enveloping algebras. In order to keep the parallel
between the two cases as close as possible we must introduce into the
picture still another ingredient:{\em \ the twisted Poisson structure on the
dual group. }This notion will also allow us to include lattice Lax equations
into our general framework. Roughly speaking, twisting is possible whenever
there is an automorphism $\tau \in Aut{\frak g}$ which preserves the
r-matrix (i.e., $\tau \otimes \tau \cdot r=r).$ Outer automorphisms are
particularly interesting. However, even in the case of inner automorphisms
the situation does not become completely trivial. The natural explanation of
twisting requires the full duality theory which involves the notions of the
double, the Heisenberg double and the twisted double of a Poisson Lie group (%
\cite{STSdouble}). In the present exposition we shall content ourselves by
presenting the final formulae for the twisted Poisson brackets. The
structure of these formulae is fairly uniform; we list them starting with
the case of finite-dimensional simple Lie algebras (where all automorphisms
are inner) and then pass to the lattice case (which accounts for the
treatment of Lax equations on the lattice) and to the case of affine Lie
algebras. In the quantum case which will be considered in the next section
twisting will play an important role in the description of quantum Casimirs.

\subsubsection{3.2.1. Finite dimensional simple Lie groups.}

We begin with a finite dimensional example and return to the setting of
section 2.2.2. Let again ${\frak g}$ be a complex simple Lie algebra, ${\frak %
b}_{+}={\frak a}+{\frak n}_{+}$ its Borel subalgebra, ${\frak b}_{-}$ the
opposite Borel subalgebra,  and $\pi _{\pm }:{\frak b}_{\pm }\rightarrow {\frak a}
$ the canonical projection, 
\[
{\frak g}^{*}=\left\{ \left( X_{+},X_{-}\right) \in {\frak b}_{+}\oplus 
{\frak b}_{-}\subset {\frak g}\oplus {\frak g};\pi _{+}\left( X_{+}\right)
=-\pi _{-}\left( X_{-}\right) \right\} . 
\]
Let $B_{\pm }\subset G$ be the corresponding Borel subgroups,  and $\pi _{\pm
}:B_{\pm }\rightarrow A$ the canonical projections, 
\begin{equation}
G^{*}=\left\{ \left( x_{+},x_{-}\right) \in B_{+}\times B_{-};\pi _{+}\left(
x_{+}\right) =\pi _{-}\left( x_{-}\right) ^{-1}\right\} .  \label{Gdual}
\end{equation}
Define $r_d\in \bigwedge^2$ $\left( {\frak g}\oplus {\frak g}\right) $ by 
\begin{equation}
r_d=\left( 
\begin{array}{cc}
r & -r-t \\
-r+t & r
\end{array}
\right) ,  \label{rdouble}
\end{equation}
where $t\in {\frak g}\otimes {\frak g}$ is the tensor Casimir. (One may
notice that $r_d$ equips ${\frak g}\oplus {\frak g}$ with the structure of a
Lie bialgebra which coincides, up to an isomorphism, with that of the {\em %
Drinfeld double} of ${\frak g.}$) Our next assertion specializes the results
stated in the previous section. Let us define the Poisson bracket on $G^{*}$
by the formula 
\begin{equation}
\left\{ \varphi ,\psi \right\} =\left\langle \left\langle r_d,D^{\prime
}\varphi \wedge D^{\prime }\psi \right\rangle \right\rangle -\left\langle
\left\langle r_d,D\varphi \wedge D\psi \right\rangle \right\rangle ,
\label{Gdualbrack}
\end{equation}
where $D\varphi =\binom{D_{+}\varphi }{D_{-}\varphi },$ $D^{\prime }\varphi =%
\binom{D_{+}^{\prime }\varphi }{D_{-}^{\prime }\varphi }$ are left and right
differentials of $\varphi $ with respect to its two arguments and $%
\left\langle \left\langle ,\right\rangle \right\rangle $ is the natural
coupling between $\bigwedge^2$ $\left( {\frak g}\oplus {\frak g}\right) $
and its dual.

\begin{proposition}
(i) $G^{*}$ is the dual Poisson Lie group of $G$ which corresponds to the
standard r-matrix described in (\ref{r-stand}). (ii) The mapping $%
G^{*}\rightarrow G:\left( x_{+},x_{-}\right) \longmapsto x_{+}x_{-}^{-1}$ is
a bijection of $G^{*}$ onto the 'big Schubert cell' in $G.$ (iii) The dual
Poisson bracket on $G^{*}$ extends smoothly from the big cell to the whole
of $G;$ it is explicitly given by (\ref{DDual}) and its Casimirs are central
functions on $G.$
\end{proposition}

The next definition prepares the introduction of {\em twisting. }

\begin{definition}
Let $\left( {\frak g},{\frak g}_r\right) $ be a factorizable Lie bialgebra
with the classical r-matrix $r\in {\frak g}\otimes {\frak g}.$ An
automorphism of $\left( {\frak g},{\frak g}_r\right) $ is an automorphism $%
\varphi \in Aut$ ${\frak g}$ which preserves the r-matrix and the inner
product on ${\frak g}$.
\end{definition}

We have the following simple result.

\begin{proposition}
\label{f}Let ${\frak g}$ be a finite-dimensional simple Lie algebra with the
standard classical r-matrix. Then $Aut\left( {\frak g},{\frak g}_r\right) $%
 coincides with the Cartan subgroup $A\subset G$.
\end{proposition}

The action of $A$ on $\wedge ^2{\frak g}$ is the restriction to $A\subset G$
of the wedge square of the standard adjoint action.

\begin{remark}
If ${\frak g}$ is only semisimple, $Aut\left( {\frak g},{\frak g}_r\right) $%
 may have a nontrivial group of components.
\end{remark}

The Sklyanin bracket on $G$ is invariant with respect to the action $A\times
G\rightarrow G$ by right translations. By contrast, the dual bracket is {\em %
not} invariant, and in this way we get a family of Poisson structures on $%
G^{*}$ with different Casimir functions. More precisely, define the action $%
A\times G^{*}\rightarrow G^{*}$ by 
\begin{equation}
h:\left( x_{+},x_{-}\right) \longmapsto \left( hx_{+},h^{-1}x_{-}\right)
\label{ac}
\end{equation}
(Note that this action is compatible with (\ref{Gdual}); the flip $%
h\longmapsto h^{-1}$ is possible, since $A$ is abelian.) We denote by $%
\lambda _h$ the contragredient action of $h$ on ${\cal F}\left[ G^{*}\right]
:$ $\lambda _h\varphi \left( x_{+},x_{-}\right) =\varphi \left(
h^{-1}x_{+},hx_{-}\right) .$ Define the twisted Poisson bracket $\left\{
,\right\} _h^{*}$ on $G^{*}$ by 
\[
\left\{ \varphi ,\psi \right\} _h^{*}=\lambda _h^{-1}\left\{ \lambda
_h\varphi ,\lambda _h\psi \right\} . 
\]
Explicitly we get 
\begin{equation}
\left\{ \varphi ,\psi \right\} _h^{*}=\left\langle \left\langle
r_d^hD^{\prime }\varphi ,D^{\prime }\psi \right\rangle \right\rangle
-\left\langle \left\langle r_dD\varphi ,D\psi \right\rangle \right\rangle ,
\label{gdualbrackh}
\end{equation}
where 
\[
r_d^h=\bigwedge^2\left( Ad\ h\oplus Ad\ h^{-1}\right) \cdot r_d=\left( 
\begin{array}{cc}
r & \ \left( -r-t\right) ^h \\ 
\left( -r+t\right) \ ^{h^{-1}} & r
\end{array}
\right) , 
\]
or, in matrix form, 
\begin{eqnarray}
&&\left\{ T_V^{\pm 1},T_W^{\pm 2}\right\} _h^{*}=\left[ r_{VW},T_V^{\pm
1}T_W^{\pm 2}\right] ,  \label{hdual} \\
&&\left\{ T_V^{+1},T_W^{-2}\right\}
_h^{*}=(r_{VW}^h+t_{VW}^h)T_V^{+1}T_W^{-2}-T_V^{+1}T_W^{-2}(r_{VW}+t_{VW}), 
\nonumber
\end{eqnarray}
where $r_{VW}^h=Ad\ h\otimes Ad\ h^{-1}\cdot r_{VW}=Ad\ h^2\otimes id\cdot
r_{VW}$ and similarly for $t_{VW}^h.$

Let us denote by $I_h\left( G\right) \ $the set of{\em \ twisted spectral
invariants}; by definition, $\phi \in I_h\left( G\right) \ $ if 
\[
\phi \left( gx\right) =\phi \left( xgh\right) \text{for any }x,g\in G. 
\]

\begin{example}
Let $\left( \rho ,V\right) $ be a linear representation of $G;$ then $\phi
:x\longmapsto tr_V\rho (hx)$ is a twisted spectral invariant.
\end{example}

\begin{proposition}
(i) There exists a unique Poisson structure on $G$ such that the{\em \ twisted
factorization map} $F_h:G^{*}\rightarrow G:\left( x_{+},x_{-}\right)
\longmapsto h^2x_{+}x_{-}$ is a Poisson mapping. (ii) $F_h$ transforms the
Casimir functions of $\left\{ ,\right\} _h^{*}$ into the set of twisted
spectral invariants $I_{h^2}\left( G\right) .$ (ii) Twisted spectral
invariants also commute with respect to the Sklyanin bracket on $G$ and
generate on $G$ generalized Lax equations 
\[
\frac{dL}{dt}=AL-LB,\;A=Ad\ h^2\cdot B.
\]
.
\end{proposition}

It is instructive to write an explicit formula for the bracket on $G$
which is the push-forward of $\left\{ ,\right\} _h^{*}$ with respect to the
factorization map. We get 
\begin{equation}
\begin{array}{l}
\left\{ \varphi ,\psi \right\} _h^{*}=\left\langle r,D\varphi \wedge D\psi
\right\rangle +\left\langle r,D^{\prime }\varphi \wedge D^{\prime }\psi
\right\rangle \\ 
\qquad -\left\langle r_{}^h,D\varphi \wedge D^{\prime }\psi \right\rangle
-\left\langle r^{h^{-1}},D^{\prime }\varphi \wedge D\psi \right\rangle \\ 
-\left\langle t_{}^h,D\varphi \wedge D^{\prime }\psi \right\rangle
+\left\langle t^{h^{-1}},D^{\prime }\varphi \wedge D\psi \right\rangle
,\;\varphi ,\psi \in C^\infty \left( G\right) ,
\end{array}
\label{Dh}
\end{equation}
where $t$ is the tensor Casimir, and we put $r^h=Ad\ h\otimes \cdot Ad\
h^{-1}\cdot r$ and $t^h=Ad\ h\otimes Ad\ h^{-1j}\cdot t.$ This bracket again
extends smoothly to the whole of $G$ (cf. (\ref{DDual}). Equivalently, 
\begin{eqnarray}
\left\{ L_V^1,L_W^2\right\} _h^{*} &=&r_{VW}L_V^1L_V^2+L_V^1L_W^2r_{VW} 
\nonumber \\
&&-L_V^1r_{VW}^hL_V^2-L_W^2r_{VW}^{h^{-1}}L_V^1  \label{dualh12} \\
&&\ -L_V^1t_{VW}^hL_W^2+L_W^2t_{VW}^{h^{-1}}L_V^1.  \nonumber
\end{eqnarray}

The twisted bracket $\left\{ ,\right\} _h^{*}$ on $G^{*}$ is {\em not}
multiplicative with respect to the group structure on $G^{*}$; to explain
its relation to the Poisson group theory we shall need some more work.

Let ${\frak a\ltimes g}^{*}$ be the semidirect product of Lie algebras which
corresponds to action (\ref{ac}); in other words, we set 
\[
\left[ H,\left( X_{+},X_{-}\right) \right] =\left( \left[ H,X_{+}\right]
,-\left[ H,X_{-}\right] \right) ,H\in {\frak a},X_{\pm }\in {\frak b}_{\pm
}. 
\]

Fix a basis $\left\{ H_i\right\} $ of ${\frak a}$ and let $\left\{
H^i\right\} $ be the dual basis of ${\frak a}^{*}$. We define a (trivial)
2-cocycle on ${\frak g}$ with values in ${\frak a}^{*}$ by 
\begin{equation}
\omega \left( X,Y\right) =\sum_i\left( H_i,\left[ X,Y\right] \right) H^i.
\label{triv}
\end{equation}
Let $\widehat{{\frak g}}={\frak g\oplus a}^{*}$ be the central extension of $%
{\frak g}$ by ${\frak a}^{*}$ which corresponds to this cocycle.

\begin{proposition}
$\left( {\frak a\ltimes g}^{*},\widehat{{\frak g}}\right) $ is a Lie
bialgebra.
\end{proposition}

The Poisson Lie group which corresponds to ${\frak a\ltimes g}^{*}$ is the
semidirect product $A\ltimes G^{*}$ associated with action (\ref{ac}).
It is easy to see that the variable $h\in A$ is central with respect to this
bracket; thus the bracket on $A\ltimes G^{*}$ splits into a family of
brackets on $G^{*}$ parametrized by $h\in A.$ This is precisely the family $%
\left\{ ,\right\} _h^{*}.$

\subsubsection{3.2.2.Twisting on the lattice.}

One may notice that the previous example is in fact trivial, since the
twisted bracket differs from the original one by a change of variables%
\footnote{%
However, even in this case twisted spectral invariants give rise to a {\em %
different} set of Hamiltonians, as compared with the non-twisted case.}.
This is of course closely related to the fact that the central extension
associated with the cocycle (\ref{triv}) is also trivial. Our next example
is more interesting; it is adapted to the study of Lax equations on the
lattice. Put $\Gamma ={\Bbb Z}/N{\Bbb Z}.$ Let again ${\frak g}$ be a finite
dimensional simple Lie algebra equipped with a standard r-matrix. Let ${\sf g%
}={\frak g}^\Gamma $ be the Lie algebra of functions on $\Gamma $ with
values in ${\frak g};$ obviously, ${\sf g}$ inherits from ${\frak g}$ the
natural Lie bialgebra structure; the corresponding r-matrix ${\sf r}\in
\bigwedge^2{\sf g}$ is given by ${\sf r}_{mn}=\delta _{mn}r,m,n\in \Gamma .$
We define ${\sf r}_d\in \bigwedge^2({\sf g}\oplus {\sf g})$ by 
\begin{equation}
{\sf r}_d=\left( 
\begin{array}{cc}
{\sf r} & -{\sf r}-{\sf t} \\ 
-{\sf r}+{\sf t} & {\sf r}
\end{array}
\right) ,  \label{ldouble}
\end{equation}
where ${\sf t\in g\otimes g}$ is the tensor Casimir which corresponds to the
standard inner product on ${\sf g},$%
\[
\left\langle \left\langle X,Y\right\rangle \right\rangle =\sum_n\left\langle
X_n,Y_n\right\rangle . 
\]
Let $\tau $ be the cyclic permutation on ${\sf g,}$ $\left( \tau X\right)
_n=X_{n+1\limfunc{mod}N}.$

\begin{lemma}
$\tau \in Aut\left( {\sf g},{\sf g}^{*}\right) .$
\end{lemma}

Let us extend the action of $\tau $ to ${\sf g\oplus g}$ in the following
way: 
\[
\tau \cdot \left( X,Y\right) =\left( \tau X,Y\right) . 
\]
(By analogy with (\ref{ac}), one might rotate the second copy of ${\sf g}$
in the opposite sense, but then the total shift in formulae below would be
by two units, which is less natural on the lattice.) Put${\sf \ }$%
\begin{eqnarray*}
{\sf r}_d^\tau &=&\bigwedge^2\left( \tau \oplus id\right) \cdot {\sf r}%
_d=\left( 
\begin{array}{cc}
{\sf r} & -{\sf r}^\tau -{\sf t}^{{\sf \tau }} \\ 
-{\sf r}^{\tau ^{-1}}+{\sf t}^{{\sf \tau }^{-1}} & {\sf r}
\end{array}
\right) , \\
{\sf r}^\tau &=&\tau \otimes id\cdot {\sf r},{\sf r}^{\tau ^{-1}}=\tau
^{-1}\otimes id\cdot {\sf r=}id\otimes \tau \cdot {\sf r.}
\end{eqnarray*}
Put ${\Bbb G}=G^\Gamma ,$ ${\Bbb G}^{*}=G^{*\ \Gamma }.$ The action of $\tau $
extends to ${\Bbb G}$ and to ${\Bbb G}^{*}$ in an obvious way. As in
section 3.2.1 we may embed ${\Bbb G}^{*}$ into the direct product ${\Bbb G}%
\times {\Bbb G}.$ By definition, the twisted Poisson bracket on ${\Bbb G}%
^{*}$ is given by 
\[
\left\{ \varphi ,\psi \right\} _\tau ^{*}=\left\langle \left\langle {\sf r}%
_d^\tau ,D\varphi \wedge D\psi \right\rangle \right\rangle -\left\langle
\left\langle {\sf r}_dD^{\prime }\varphi \wedge D^{\prime }\psi
\right\rangle \right\rangle , 
\]
where $D\varphi ,D^{\prime }\varphi \in {\sf g\oplus g}$ are the
two-component left and right differentials of $\varphi .$

We may push forward the Poisson bracket from ${\Bbb G}^{*}$ to ${\Bbb G}$
using the factorization map. More precisely:

\begin{proposition}
There exists a unique Poisson structure on ${\Bbb G}$ (which we shall still
denote by $\left\{ ,\right\} _\tau ^{*}$ ) such that the {\em twisted
factorization map} ${\Bbb G}^{*}\rightarrow {\Bbb G}:\left(
x_{+},x_{-}\right) $ $\longmapsto x_{+}^\tau x_{-}^{-1}$ becomes a Poisson
mapping. This Poisson structure is given by
\end{proposition}

\begin{equation}
\begin{array}{l}
\left\{ \varphi ,\psi \right\} _\tau ^{*}=\left\langle r,D\varphi \wedge
D\psi \right\rangle +\left\langle r,D^{\prime }\varphi \wedge D^{\prime
}\psi \right\rangle \\ 
\qquad -\left\langle r^{\tau ^{}},D\varphi \wedge D^{\prime }\psi
\right\rangle -\left\langle r^{\tau ^{-1}},D^{\prime }\varphi \wedge D\psi
\right\rangle \\ 
-\left\langle t^{\tau ^{}},D\varphi \wedge D^{\prime }\psi \right\rangle
+\left\langle t^{\tau ^{-1}},D^{\prime }\varphi \wedge D\psi \right\rangle
,\;\varphi ,\psi \in C^\infty \left( {\Bbb G}\right) ,
\end{array}
\label{Dtau}
\end{equation}

We shall denote the Lie group ${\Bbb G}$ equipped with the Poisson structure
(\ref{Dtau}) by ${\Bbb G}_\tau ,$ for short. Equivalent description of the
bracket $\left\{ ,\right\} _\tau ^{*}$ for linear groups: let again $%
L_V=\rho _V\left( L\right) ;$ then 
\begin{eqnarray}
\left\{ L_V^1,L_V^2\right\} ^{*} &=&rL_V^1L_V^2+L_V^1L_V^2r-L_V^1r^\tau
L_V^2-L_V^2r^{\tau ^{-1}}L_V^1  \label{dtau12} \\
&&\ \ -L_V^1t^\tau L_V^2-L_V^2t^{\tau ^{-1}}L_V^1  \nonumber
\end{eqnarray}
Another obvious Poisson structure on ${\Bbb G}$ is the Sklyanin bracket
which corresponds to the Lie bialgebra $\left( {\sf g},{\sf g}^{*}\right) .$
Let us consider the following action of ${\Bbb G}$ on ${\Bbb G}_\tau $
called the (lattice) gauge action: 
\[
x:T\longmapsto x^\tau Tx^{-1}. 
\]

\begin{theorem}
(i) The gauge action ${\Bbb G}\times {\Bbb G}_\tau \rightarrow {\Bbb G}_\tau 
$ is a Poisson group action. (ii) The Casimir functions of $\left\{
,\right\} _\tau ^{*}$ coincide with the gauge invariants. (iii) Let $M:{\Bbb %
G}_\tau \rightarrow G$ be the monodromy map, 
\[
M:T=\left( T_1,T_2,...,T_N\right) \longmapsto T_N\cdot \cdot T_1;
\]
we equip the target space with Poisson structure (\ref{DDual}). Then $M$
is a Poisson map and the Casimirs of ${\Bbb G}_\tau $ coincide with the
spectral invariants of the monodromy. (iv) Spectral invariants of the
monodromy are in involution with respect to the Sklyanin bracket on ${\Bbb G}
$ and generate lattice Lax equations on ${\Bbb G}$.
\end{theorem}

Parallels with theorem \ref{AKS} are completely obvious.

\subsubsection{3.2.3. Twisting for loop algebras.}

Our last example of twisting is based on the outer automorphisms of loop
algebras. We shall use the notation and definitions introduced in section
2.2.2. Let ${\frak g}$ be a simple Lie algebra, $L{\frak g}$ the
corresponding loop algebra with the standard trigonometric r-matrix,  and let $L%
{\frak g}_r=\left( L{\frak g}\right) ^{*}$ be the dual algebra. The Lie group
which corresponds to $L{\frak g}$ is the group $LG$ of polynomial loops with
values in $G;$ the elements of the dual Poisson group may be identified with
pairs of formal series $\left( T^{+}(z),T^{-}(z^{-1})\right) ,$%
\[
T^{\pm }(z)=\sum_{n\geq 0}T^{\pm }\left[ \pm n\right] z^{\pm n},T^{\pm
}\left[ 0\right] \in B_{\pm },\pi _{+}\left( T^{+}\left[ 0\right] \right)
\cdot \pi _{-}\left( T^{-}\left[ 0\right] \right) =1 
\]
(here $\pi _{\pm }:B_{\pm }\rightarrow H\ $are canonical projections). The
group $\left( LG\right) ^{*}$ is pro-algebraic and its affine ring is
generated by the matrix coefficients of $T^{\pm }\left[ n\right] .$

The next proposition describes the Poisson structure on the affine ring of $%
\left( LG\right) ^{*}.$ Let $t\in {\frak g}\otimes {\frak g}$ be the
canonical element (tensor Casimir) which represents the inner product $%
\left\langle ,\right\rangle $ in ${\frak g.}$ Let $\delta \left( x\right) $
be the Dirac delta, 
\[
\delta \left( x\right) =\sum_{n=-\infty }^\infty x^n 
\]
(observe that $t\cdot \delta \left( z/w\right) $ represents the kernel of
the identity operator acting in $L{\frak g).}$

\begin{proposition}
The Poisson bracket on $\left( LG\right) ^{*}$ which corresponds to the
standard Lie bialgebra structure on $L{\frak g}$ is given by 
\begin{eqnarray}
\left\{ T_1^{\pm }(z),T_2^{\pm }(w)\right\} =\left[ r\left( \frac zw\right)
,T_1^{\pm }(z)T_2^{\pm }(w)\right] ,  \label{t12} \\
\left\{ T_1^{+}(z),T_2^{-}(w)\right\} =\left[ r\left( \frac zw\right)
,T_1^{+}(z)T_2^{-}(w)\right] +\delta \left( \frac zw\right) \left[
t,T_1^{+}(z)T_2^{-}(w)\right] ,  \nonumber
\end{eqnarray}
where $r\left( x\right) $ is the trigonometric r-matrix.
\end{proposition}

We pass to the description of the twisted structure on $\left( LG\right)
^{*}.$

\begin{lemma}
$Aut\left( L{\frak g},L{\frak g}_r\right) $ is the extended Cartan subgroup $%
\hat H=H\times {\Bbb C}^{\times }$
\end{lemma}

By definition, $p\in {\Bbb C}^{\times }$ acts on $L{\frak g}$ by
sending $X\left( z\right) \longmapsto X\left( pz\right) ;$ we shall denote
the corresponding dilation operator by $D_p.$ The action of $H$ gives
nothing new, as compared with our first example (section 2.2.1). To describe
the twisting of the Poisson structure on $\left( LG\right) ^{*}$ by the
action of ${\Bbb C}^{\times }$ we start with the r-matrix 
\[
r_d(z/w)=\left( 
\begin{array}{cc}
r(z/w) & -r(z/w)-t\delta (z/w) \\ 
-r(z/w)+t\delta (z/w) & r(z/w)
\end{array}
\right) 
\]
(as usual, $r_d(z/w)$ does not belong to the $\bigwedge^2\left( L{\frak g}%
\oplus L{\frak g}\right) $, but is well defined as a kernel of a linear
operator acting in $L{\frak g}\oplus L{\frak g).}$ For $p\in {\Bbb C}%
^{\times }$ we define the twisted r-matrix by 
\[
r_d(z/w)^p=\bigwedge^2(D_p\oplus D_p^{-1})\cdot r_d(z/w)=\left( 
\begin{array}{cc}
r(z/w) & -r(pz/w)-t\delta (pz/w) \\ 
-r(z/pw)+t\delta (z/pw) & r(z/w)
\end{array}
\right) . 
\]

\begin{definition}
The twisted Poisson bracket on $\left( LG\right) ^{*}$ is given by 
\begin{eqnarray}
\left\{ T_1^{\pm }(z),T_2^{\pm }(w)\right\} _p^{*}=\left[ r\left( \frac
zw\right) ,T_1^{\pm }(z)T_2^{\pm }(w)\right] ,  \nonumber \\
\left\{ T_1^{+}(z),T_2^{-}(w)\right\} _p^{*}=r\left( \frac{pz}w\right)
T_1^{+}(z)T_2^{-}(w)-T_1^{+}(z)T_2^{-}(w)r\left( \frac zw\right) 
\label{tp12} \\
+\delta \left( \frac{pz}w\right) t\
T_1^{+}(z)T_2^{-}(w)-T_1^{+}(z)T_2^{-}(w)\delta \left( \frac zw\right) t 
\nonumber
\end{eqnarray}
\[
\]
\end{definition}

Define the new generating function (`full current') by 
\begin{equation}
T\left( z\right) =T^{+}(pz)T^{-}(p^{-1}z)^{-1};  \label{cur}
\end{equation}
clearly, the coefficients of $T\left( z\right) $ generate the affine ring of 
$\left( LG\right) ^{*}.$

\begin{proposition}
The current $T\left( z\right) $ satisfies the following Poisson bracket
relations: 
\begin{eqnarray}
\left\{ T_1\left( z\right) ,T_2\left( w\right) \right\} _p^{*}=  \nonumber \\
\ \ \ r\left( \frac zw\right) T_1\left( z\right) T_2\left( w\right)
+T_1\left( z\right) T_2\left( w\right) r\left( \frac zw\right)   \nonumber \\
\ \ \ -T_1\left( z\right) r\left( \frac{p^2z}w\right) T_2\left( w\right)
-T_2\left( w\right) r\left( \frac z{p^2w}\right) T_1\left( z\right) 
\label{Dp} \\
\ \ -T_1\left( z\right) \delta \left( \frac{p^2z}w\right) \ t\ T_2\left(
w\right) +T_2\left( w\right) \delta \left( \frac z{p^2w}\right) \ t\
T_1\left( z\right) .  \nonumber
\end{eqnarray}
\end{proposition}

The existence of twisting in the affine case is due to the existence of the
non-trivial central extension of the loop algebra. (This should be compared
with the situation in our first example, where the central extension is
trivial because all automorphisms of simple Lie algebras are inner.)  More
precisely, we have the following result. Let $\omega $ be the 2-cocycle on $L%
{\frak g}$ defined by 
\begin{equation}
\omega \left( X,Y\right) =Res_{z=0}\left\langle X\left( z\right) ,z\partial
_zY\left( z\right) \right\rangle dz/z.
\end{equation}
Let $\widehat{L{\frak g}}$ be the central extension of $L{\frak g}$
determined by this cocycle,  and $\widehat{L{\frak g}}^{*}$ the semidirect
product of $\left( L{\frak g}\right) ^{*}$ and ${\Bbb C}^{\times }.$

\begin{proposition}
$\left( \widehat{L{\frak g}},\widehat{L{\frak g}}^{*}\right) $ is a Lie
bialgebra.
\end{proposition}

The Poisson Lie group which corresponds to $\widehat{L{\frak g}}^{*}$ is the
semidirect product ${\Bbb C}^{\times }\ltimes (LG)^{*};$ it is easy to see
that the affine coordinate on ${\Bbb C}^{\times }$ is central with respect
to the Poisson structure on ${\Bbb C}^{\times }\ltimes (LG)^{*};$ fixing its
value, we get the family of Poisson brackets $\left\{ ,\right\} _p^{*}$ on $%
(LG)^{*}.$

In section 4.2 we shall discuss a similar phenomenon for quantum systems
associated with quasitriangular Hopf algebras.

\subsection{$\label{Sklyanin}$Sklyanin Bracket on ${\bf G}\left( z\right) $}

In our previous example we have described the Poisson structure associated
with the loop algebra, namely, the bracket on the dual group $\left(
LG\right) ^{*};$ applications to lattice systems involve the Poisson
structure on the loop group itself. As a matter of fact, this was the first
example of a Poisson Lie group which appeared prior to their formal
definition (\cite{rmat}). To put it more precisely, let ${\frak g}$ be a
complex semisimple Lie algebra, $G$ a linear complex Lie group with the Lie
algebra ${\frak g.}$ Let $G(z)$ be the associated loop group, i.e., the group
of rational functions with values in $G.$ We may regard $G(z)$  as an affine
algebraic group defined over ${\Bbb C}\ (z)$ . Let us again fix a faithful
representation $\left( \rho ,V\right) $ of $G;$ the affine ring ${\cal A}%
_{aff}$ of $G(z)$ is generated by `tautological functions' which assign to $%
L\in G\left( z\right) $ the matrix coefficients of $L_V\left( z\right) \in
EndV\otimes {\Bbb C}(z).$ The Poisson bracket on the affine ring of $G(z)$
is defined by the following formula (\cite{rmat}), 
\begin{equation}
\left\{ L_V^1(u),L_V^2(v)\right\} =\left[ r_V(u,v),L_V^1(u)L_V^2(v)\right].
\label{affskl}
\end{equation}
(By an abuse of notation, we do not distinguish between the generators of $%
{\cal A}_{aff}$ and their values at $L\in G\left( z\right) $ in the l.h.s.
of this formula.)  Bracket (\ref{affskl}) defines the structure
of a Poisson-Lie group on $G(z)$. Equivalent formulation:

{\em Let }$\Delta :{\cal A}_{aff}\rightarrow {\cal A}_{aff}\otimes {\cal A}%
_{aff}${\em \ be the coproduct in }${\cal A}_{aff}${\em \ induced by the
group multiplication in }$G\left( z\right) :$%
\[
\Delta \left( L_V\left( z\right) \right) =L_V\left( z\right) \dot \otimes
L_V\left( z\right) , 
\]
{\em or, in a less condensed notation, } 
\begin{equation}
\Delta \left( L_V^{ij}\left( z\right) \right) =\sum_kL_V^{ik}\left( z\right)
\otimes L_V^{kj}\left( z\right) .
\end{equation}
{\em Then }$\Delta ${\em \ is a morphism of Poisson algebras (in other
words, } 
\[
\left\{ \Delta \varphi ,\Delta \psi \right\} =\Delta \left\{ \varphi ,\psi
\right\} 
\]
{\em for any }$\varphi ,\psi \in {\cal A}_{aff})${\em .}

The Lax matrix $L_V\left( u\right) $ in this context is basically a
tautological mapping which assigns to an element $L\in G(z)$ the matrix $%
L_V\left( u\right) \in End$ $V\otimes {\Bbb C}\left( u\right) .$
Alternatively, for any Poisson submanifold ${\cal M\subset }$ $G(z)$ the Lax
matrix may be regarded as an embedding map ${\cal M}\hookrightarrow
G(z)\rightarrow End$ $V\otimes {\Bbb C}\left( u\right) .$ An explicit
description of all Poisson submanifolds in $G(z)$ follows from the theory of 
{\em dressing transformations} (\cite{CIMPA}); for our present goals it
suffices to know that rational functions with a prescribed divisor of poles
form a finite-dimensional Poisson submanifold in $G(z)$; moreover, since the
Poisson structure is multiplicative, the product of Poisson submanifolds is
again a Poisson submanifold. Generic Poisson submanifolds correspond to
functions with simple poles which may be written in multiplicative form 
\begin{equation}
L_V(z)=\stackrel{\curvearrowleft }{\prod_i}(I-\frac{X_i}{z-z_i}).
\label{monc}
\end{equation}
In other words, the description of Poisson submanifolds is tantamount to the
choice of an {\em Ansatz} for the Lax matrix; the free parameters in this
Ansatz (e.g., the residues $X_i$) become dynamical variables. Spectral
invariants of 'Lax matrices' of this type may be used to generate completely
integrable lattice Lax equations. More precisely, to get such system one may
proceed as follows:

\begin{itemize}
\item  Pick a Poisson submanifold ${\cal M}\subset G(z).$

\item  Consider the product space ${\cal M}^N\subset G(z)^{{\Bbb \ Z}/N{\Bbb %
Z}}$ and the monodromy map\ $m:G(z)\times ...\times G(z)\rightarrow G(z)$
(more precisely, its restriction to ${\cal M}^N).$

\item  Choose any central function ${\cal H}$ on $G(z)$ as a Hamiltonian;
its pullback to ${\cal M}^N$ defines a lattice Lax equation.
\end{itemize}

Under some mild assumptions on the choice of ${\cal M},$ ${\cal H}$ defines
a completely integrable system on $m\left( {\cal M}^N\right) $ ; its
pullback to ${\cal M}^N$ still remains completely integrable; this may be
regarded as the main content of the Inverse Scattering Method (in this
slightly simplified setting).

The study of Lax systems on the lattice thus breaks naturally into two
parts: {\it (a) Solve the Lax equation for the monodromy. (b) Lift the
solutions back to $G(z)^N.$} The second stage (by no means trivial) is the
inverse problem, {\em \ sensu strictu.}

\section{\label{Quant}Quantization}

\subsection{\label{Quant lin}Linear case}

Speaking informally, quantization consists in replacing Poisson bracket
relations with commutation relations in an associative algebra. We shall
have to distinguish between linear and quadratic cases (as explained in the
previous section, this difference stems from the different Hopf structures
on the algebras of observables). The linear case is easier, since we remain
in the conventional setting of Lie groups and Lie algebras. Quantization of
quadratic Poisson bracket relations $\left( \ref{skl}\right) $ leads to {\em %
quasitriangular Hopf algebras}.

The standard way to quantize the algebra ${\cal \ A}_{cl}={\cal S}({\frak g}%
) $ (equipped, as usual, with the Lie-Poisson bracket of ${\frak g}$ ) is to
replace it with the universal enveloping algebra ${\cal U}({\frak g}).$ Let
me briefly recall the corresponding construction.

Let ${\cal U}({\frak g})=\bigcup_{n\geq 0}{\cal U}_n$ be the canonical
filtration of ${\cal U}({\frak g}),$  and ${\cal S}({\frak g})=\oplus _{n\geq 0}%
{\cal S}_n$ the canonical grading of ${\cal S}({\frak g})$. Clearly, ${\cal S%
}_n\simeq {\cal U}_n/{\cal U}_{n-1}.$ Let $gr_n:{\cal U}_n\rightarrow {\cal S%
}_n$ be the canonical projection. If $u\in {\cal U}_n,v\in {\cal U}_m$,  their
commutator $\left[ u,v\right] \in {\cal U}_{n+m-1}.$ It is easy to see that $%
gr_{n+m-1}\left( \left[ u,v\right] \right) =\left\{ gr_nu,gr_mv\right\} $ is
precisely the Lie-Poisson bracket; hence ${\cal U}({\frak g})$ is a
quantization of ${\cal S}({\frak g}).$

Alternatively, let $\hbar $ be a formal parameter. Put ${\frak g}_\hbar =%
{\frak g}\otimes _kk\left[ \left[ \hbar \right] \right] $; we rescale the
commutator in ${\frak g}_\hbar $ by putting $\left[ X,Y\right] _\hbar =\hbar
\left[ X,Y\right] .$ Let ${\cal U}\left( {\frak g}_\hbar \right) $be the
universal enveloping algebra of ${\frak g}_\hbar ,$  and $\hbar \,{\cal U}\left( 
{\frak g}_\hbar \right) $  be its maximal ideal consisting of formal series in $%
\hbar $ with zero constant term; the quotient algebra ${\cal U}\left( {\frak %
g}_\hbar \right) /\hbar \,{\cal U}\left( {\frak g}_\hbar \right) $ is
canonically isomorphic to ${\cal S}\left( {\frak g}\right) .$ Let $p:{\cal U}%
\left( {\frak g}_\hbar \right) \rightarrow {\cal S}\left( {\frak g}\right) $
be the canonical projection; then $\left\{ p\left( x\right) ,p\left(
y\right) \right\} =p\left( \hbar ^{-1}\left( xy-yx\right) \right) $ for any $%
x,y\in {\cal U}\left( {\frak g}_\hbar \right) $. (In the sequel we prefer to
set $\hbar =1$ and so the first definition will be more convenient.)

Let $r\in End{\frak g}$ be a classical r-matrix. As we have seen, the phase spaces
of Lax systems associated with $\left( {\frak g},r\right) $ are coadjoint
orbits of ${\frak g}_r;$ the Hamiltonians are obtained from the Casimir
elements $H\in {\cal S}\left( {\frak g}\right) ^{{\frak g}}$ via the
isomorphism ${\cal S}\left( {\frak g}_r\right) \rightarrow {\cal S}\left( 
{\frak g}\right) $ described in theorem \ref{AKS}. The quantum counterpart
of ${\cal S}\left( {\frak g}\right) ^{{\frak g}}$ is the center ${\cal %
Z\subset U}\left( {\frak g}\right) $ of the universal enveloping algebra.
Moreover, under favorable conditions there is a nice correspondence between
coadjoint orbits of a Lie algebra and unitary representations of the
corresponding Lie group. Thus we may establish the following heuristic
vocabulary which describes the correspondence between classical and quantum
systems: 
\[
\begin{tabular}{|c|c|}
\hline
{\bf Linear Classical Case} & {\bf Linear Quantum case} \\ \hline
${\cal S}({\frak g})$ & ${\cal U}({\frak g})$ \\ \hline
{\em Classical r-matrices} & {\em Classical r-matrices} \\ \hline
{\em Casimir functions} $\left( {\cal S}\left( {\frak g}\right) \right) ^{%
{\frak g}}\subset {\cal S}\left( {\frak g}\right) $ & {\em Casimir operators}
${\cal Z}\subset {\cal U}\left( {\frak g}\right) $ \\ \hline
${\cal S}({\frak g}_r)$ & ${\cal U}({\frak g}_r)$ \\ \hline
{\em Coadjoint orbits in} ${\frak g}_r^{*}$ & {\em Irreducible }${\cal U}%
\left( {\frak g}_r\right) ${\em -modules} \\ \hline
\end{tabular}
\]

Notice that, in the linear case, classical r-matrices are used in the quantum case
as well!

The following theorem is an exact analogue of theorem \ref{AKS}

Let $\left( {\frak g},r\right) $ be a Lie algebra equipped with a classical
r-matrix $r\in End{\frak g}$ satisfying $\left( \ref{cybe}\right) .$ Let $%
{\frak g}_r$ be the corresponding Lie algebra (with the same underlying
linear space). Put $r_{\pm }=\frac 12(r\pm id);$ then (\ref{cybe}) implies
that $r_{\pm }:{\frak g}_r\rightarrow {\frak g}$ are Lie algebra
homomorphisms. Extend them to homomorphisms ${\cal U}({\frak g}%
_r)\rightarrow {\cal U}({\frak g})$ which we denote by the same letters.
These morphisms also agree with the standard Hopf structure on ${\cal U}(%
{\frak g})$, ${\cal U}({\frak g}_r).$ Define the action 
\begin{equation}
{\cal U}({\frak g}_r)\otimes {\cal U}({\frak g})\rightarrow {\cal U}\left( 
{\frak g}\right)  \nonumber
\end{equation}
by setting 
\begin{equation}
x\cdot y=\sum r_{+}(x_i^{(1)})\;y\;r_{-}(x_i^{(2)})^{\prime },x\in {\cal U}(%
{\frak g}_r),y\in {\cal U}({\frak g}),  \label{qact}
\end{equation}
where $\Delta x=\sum x_i^{(1)}\otimes x_i^{(2)}$ is the coproduct and $%
a\mapsto a^{\prime }$ is the antipode map.

\begin{theorem}
\label{K} (i) ${\cal U(}{\frak g})$ is a free filtered Hopf ${\cal U(}{\frak %
g}_r)$-module generated by $1\in {\cal U}({\frak g}).$ (ii) Let $i:{\cal U(%
}{\frak g})$ $\rightarrow {\cal U(}{\frak g}_r)$ be the induced isomorphism
of filtered linear spaces; its restriction to ${\cal Z}\subset {\cal U}%
\left( {\frak g}\right) $ is an algebra homomorphism; in particular, $%
i\left( {\cal Z}\right) \subset {\cal U}({\frak g}_r)$ is commutative.
\end{theorem}

Let ${\frak u}\subset {\frak g}_r$ be an ideal, ${\frak s}={\frak g}_r/%
{\frak u}$ the quotient algebra,  and $p:{\cal U}\left( {\frak g}_r\right)
\rightarrow {\cal U}\left( {\frak s}\right) $ be the canonical projection;
restricting $p$ to the subalgebra $I_r=i_r\left( {\cal Z}\right) $ we get a
commutative subalgebra in ${\cal U}\left( {\frak s}\right) ;$ we shall say
that the corresponding elements of ${\cal U}\left( {\frak s}\right) $ are
obtained by {\em specialization.}

In particular, let ${\frak g}={\frak g}_{+}{\frak \dot +g}_{-}$ be a
splitting of ${\frak g}$ into a linear sum of two its Lie subalgebras; let $%
P:{\cal U}({\frak g})\rightarrow {\cal U}({\frak g}_{-})$ be the projection
onto ${\cal U}({\frak g}_{-})$ in the decomposition 
\[
{\cal U}({\frak g})={\cal U}({\frak g}_{-})\oplus {\frak g}_{+}{\cal U}(%
{\frak g}){\bf .} 
\]

\begin{corollary}
\label{CorK} The restriction of $P$ to the subalgebra ${\cal Z}\subset {\cal %
U}({\frak g})$ of Casimir elements is an algebra homomorphism.
\end{corollary}

\subsubsection{4.1.1. Quantum reduction.}

The quantum analogue of symplectic induction discussed in section 2.1.1 is
the ordinary induction. For completeness we recall the standard definitions
(cf., for example, \cite{dixmier}). Let ${\frak g}$ be a Lie algebra. Let us
denote by ${\sf Rep}_{{\frak g}}$ the category of ${\cal U}\left( {\frak g}%
\right) $ -modules. Let ${\frak b}\subset {\frak g}$ be a Lie subalgebra. In
complete analogy with section 2.2.1 we construct the induction functor $Ind_{%
{\frak b}}^{{\frak g}}:{\sf Rep}_{{\frak b}}\rightsquigarrow {\sf Rep}_{%
{\frak g}}$ which associates to each ${\cal U}\left( {\frak b}\right) $%
-module an ${\cal U}\left( {\frak g}\right) $-module. Namely, we put $Ind_{%
{\frak b}}^{{\frak g}}\left( V\right) =$ ${\cal U}\left( {\frak g}\right)
\bigotimes_{{\cal U}\left( {\frak b}\right) }V,$ where ${\cal U}\left( 
{\frak g}\right) $ is regarded as a right ${\cal U}\left( {\frak b}\right) $%
-module. [By definition, ${\cal U}\left( {\frak g}\right) \bigotimes_{{\cal U%
}\left( {\frak b}\right) }V$ is the quotient of ${\cal U}\left( {\frak g}%
\right) \bigotimes_{{\Bbb C}}V$ over its submodule generated by $\left(
u\otimes v\right) \cdot X,$ $u\in {\cal U}\left( {\frak g}\right) ,v\in V,$ $%
X\in {\cal U}\left( {\frak b}\right) ;$ the (right) action of ${\cal U}%
\left( {\frak b}\right) $ on ${\cal U}\left( {\frak g}\right) \bigotimes_{%
{\Bbb C}}V$ is defined by 
\[
\left( u\otimes v\right) \cdot X=\sum_iuS\left( X^{\left( i\right) }\right)
\otimes X_{\left( i\right) }v, 
\]
where $\Delta X=\sum_iX^{\left( i\right) }\otimes X_{\left( i\right) }$ is
the coproduct in ${\cal U}\left( {\frak b}\right) $ and $S$ is its antipode.
In physical terms passing to the quotient means that we impose 'constraints
on the wave functions'; these constraints express invariance of wave
functions with respect to the diagonal action of ${\cal U}\left( {\frak b}%
\right) .$] The structure of a ${\cal U}\left( {\frak g}\right) $-module in $%
Ind_{{\frak b}}^{{\frak g}}\left( V\right) $ is induced by the left action
of ${\cal U}\left( {\frak g}\right) $ on itself.

Fix a point $F\in {\frak b}^{*}$ and assume that ${\frak l}_F$ is a
Lagrangian subalgebra subordinate to $F.$ In that case $F$ defines a
character of ${\frak l}_F;$ let $V_F$ be the corresponding 1-dimensional $%
{\cal U}\left( {\frak l}_F\right) $-module. Put $V={\cal U}\left( {\frak b}%
\right) \bigotimes_{{\cal U}\left( {\frak l}_F\right) }V_F.$ $V$ is a
natural ${\cal U}\left( {\frak b}\right) $-module associated with the
coadjoint orbit of $F.$

\begin{proposition}
$Ind_{{\frak b}}^{{\frak g}}\left( V\right) \simeq Ind_{{\frak l}%
_F}^{{\frak g}}\left( V_F\right) .$
\end{proposition}

Informally, we may say that an ${\cal U}({\frak g)}$-module associated with
a coadjoint orbit admitting a Lagrangian polarization may be induced from a
1-dimensional module.

Let us now return to the setting of theorem \ref{K}. Let  $\left( 
{\frak g},r\right) $ be a Lie algebra equipped with a classical r-matrix $%
r\in End\ {\frak g}$ satisfying $\left( \ref{cybe}\right) .$ We regard $%
{\cal U}\left( {\frak g}\right) $ as a ${\cal U}\left( {\frak g}_r\right) $%
-module with respect to the action (\ref{qact}). Let $V$ be a ${\cal U}%
\left( {\frak g}_r\right) $-module.

\begin{proposition}
$V\simeq {\cal U}\left( {\frak g}\right) \bigotimes_{{\cal U}\left( {\frak g}%
_r\right) }V$ as a linear space.
\end{proposition}

Since half of ${\cal U}\left( {\frak g}_r\right) $ is acting on ${\cal U}%
\left( {\frak g}\right) $ on the left and the other half on the right, the
structure of a ${\cal U}\left( {\frak g}\right) $-module in ${\cal U}\left( 
{\frak g}\right) \bigotimes_{{\cal U}\left( {\frak g}_r\right) }V$ is
destroyed. However, the structure of ${\cal Z}$-module in ${\cal U}\left( 
{\frak g}\right) $ survives tensoring with $V$ over ${\cal U}\left( {\frak g}%
_r\right) $.

\begin{proposition}
For any $\zeta \in {\cal Z},$ $u\in {\cal U}\left( {\frak g}\right) ,v\in V,$
   $\zeta u\otimes v=u\otimes i_r\left( \zeta \right) v.$
\end{proposition}

This gives a new proof of theorem \ref{K}.

Fix $F\in {\frak g}_r^{*}$ and let ${\frak l}_F\subset {\frak g}_r$ be a
Lagrangian subalgebra subordinate to $F$. Let $V_F$ be the corresponding 1-dimensional ${\cal U}%
\left( {\frak l}_F\right) $ -module. Then ${\cal U}\left( {\frak g}\right)
\bigotimes_{{\cal U}\left( {\frak l}_F\right) }V_F$ may be identified with $%
V.$

Let us assume now that ${\frak g}={\frak g}_{+}\ \dot +\ {\frak g}_{-}$ and
the r-matrix is given by (\ref{standard}); in that case ${\cal U}\left( 
{\frak g}_r\right) \simeq $ ${\cal U}\left( {\frak g}_{+}\right) \otimes $ $%
{\cal U}\left( {\frak g}_{-}\right) .$ Let $W$ be a ${\cal U}\left( {\frak g}%
_{-}\right) $-module, ${\bf W}={\cal U}\left( {\frak g}\right) \bigotimes_{%
{\cal U}\left( {\frak g}_{-}\right) }W;$ we regard ${\bf W}$ as a left $%
{\cal U}\left( {\frak g}\right) $-module. There is a canonical embedding $%
W\hookrightarrow {\bf W}:w\mapsto 1\otimes w.$ Let ${\bf W}^{*}$ be the dual
module regarded as a right ${\cal U}\left( {\frak g}\right) $-module,  and let $%
W_0^{*}\subset {\bf W}$  be the subspace of ${\cal U}\left( {\frak g}_{+}\right) 
$-invariants.

\begin{lemma}
$W_0^{*}$ is isomorphic to the dual of $W.$
\end{lemma}

This is an easy corollary of the decomposition 
\[
{\cal U}\left( {\frak g}\right) \simeq {\cal U}\left( {\frak g}_{+}\right)
\otimes {\cal U}\left( {\frak g}_{-}\right) \simeq {\frak g}_{+}{\cal U}%
\left( {\frak g}\right) \oplus {\cal U}\left( {\frak g}_{-}\right) . 
\]
Fix a basis $\left\{ e_i\right\} $ in $W,$ and let $\left\{ e^i\right\} $ be
the dual basis in $W_0^{*}.$ Let $\Omega =\sum e^i\otimes e_i\in W_0^{*}{\bf %
\otimes }W$ be the canonical element; it defines a natural mapping 
\[
{\bf W}\rightarrow W:\varphi \mapsto \left\langle \varphi \right\rangle
_\Omega =\sum_i\left\langle e^i,\varphi \right\rangle e_i. 
\]

\begin{proposition}
For any $\varsigma \in {\cal Z},$ $\varphi \in {\bf W},$    $%
\left\langle \varsigma \varphi \right\rangle _\Omega =i_\varsigma $ $%
\left\langle \varphi \right\rangle _\Omega .$
\end{proposition}

Informally, we may say that quantum Hamiltonians acting in $W$ are {\em %
radial parts} of the Casimir operators. (The reasons for this terminology
will be obvious from our next example.)

\subsubsection{4.1.2. Toda lattice\label{quant toda}.}

Our first example is the generalized open Toda lattice. We retain the
notation of section 2.2.1. Let again ${\frak g}$ be a real split semisimple
Lie algebra, ${\frak g}={\frak k}\dot +{\frak a}\dot +{\frak n}$ its Iwasawa
decomposition, ${\frak b}={\frak a}\dot +{\frak n}$ the Borel subalgebra.
Let ${\frak d}_p\subset {\frak b}$ ($p\geq 0)$ be the generalized diagonals
defined in (\ref{diag}). The subalgebras ${\frak b}^{(p)}=\oplus _{r\geq p}%
{\frak d}_r$ define a decreasing ad$\ {\frak b}$-invariant filtration in $%
{\frak b.}$ Put ${\frak p}_k=\oplus _{0\leq p\leq k}s\left( {\frak d}%
_p\right) ,$ where $s:$ ${\frak g}\rightarrow {\frak p}:X\longmapsto \frac
12\left( id-\sigma \right) X$ is the projection onto the subspace of
anti-invariants of the Cartan involution $\sigma $. Put ${\frak s=b}/{\frak b%
}^{\left( 2\right) }$. Let $L$ $\in {\frak p}_1\otimes {\frak s}\subset 
{\frak p}_1\otimes {\cal U}\left( {\frak s}\right) $ be the canonical
element induced by the natural pairing ${\frak p}_1$ $\times $ ${\frak b}/%
{\frak b}^{\left( 2\right) }\rightarrow {\Bbb R}$. Let $P:{\cal U}({\frak g}%
)\rightarrow {\cal U}\left( {\frak b}\right) $ be the projection map onto $%
{\cal U}\left( {\frak b}\right) $ in the decomposition ${\cal U}({\frak g})=%
{\cal U}({\frak b})\oplus {\cal U}({\frak g}){\frak k}{\bf \ }$, $p:{\cal U}(%
{\frak b})\rightarrow U\left( {\frak s}\right) $ the specialization map. Let
us put $I_{{\frak s}}=p\circ P\left( {\cal Z}\right) ,$ where as usual $%
{\cal Z}$ is the center of ${\cal U}\left( {\frak g}\right) .$ Fix a
faithful representation $\left( \rho ,V\right) $ of ${\frak g}$ and put $%
L_V=\left( \rho \otimes id\right) L.$

\begin{proposition}
(i) $I_{{\frak s}}$ is a commutative subalgebra in ${\cal U}\left( {\frak s}%
\right) ;$ moreover, ${\cal Z}\rightarrow I_{{\frak s}}:z\longmapsto H_z$ is
an algebra isomorphism. (ii) $I_{{\frak s}}$ is generated by $%
H_k=tr_VL_V^k\in {\cal U}\left( {\frak s}\right) ,k=1,2,...$
\end{proposition}

Recall that ${\frak s}={\frak a}\ltimes {\frak u,}$ where ${\frak u}={\frak %
n/}\left[ {\frak n},{\frak n}\right] \ ,$ and the corresponding Lie group is 
$S=A\ltimes U,$ $U=N/N^{\prime }.$ Let ${\cal O}_T$ be the coadjoint orbit
of $S$ described in (\ref{orbtoda}), $f=\sum_{\alpha \in P}\left( e_\alpha
+e_{-\alpha }\right) \in $ ${\cal O}_T$ the marked point, $L_f=U$ the
corresponding Lagrangian subgroup, $V_f$ $={\Bbb C}$ the 1-dimensional $U$%
-module which corresponds to $f;$ we may also regard it as a ${\cal U}\left( 
{\frak u}\right) $-module.

\begin{proposition}
(i) The irreducible ${\cal U}\left( {\frak s}\right) $-module associated
with ${\cal O}_T$ is the induced module ${\cal U}\left( {\frak s}\right)
\otimes _{{\cal U}\left( {\frak u}\right) }V_f.$ (ii) The corresponding
unitary representation space ${\frak H}$ may be identified with $L_2\left( 
{\frak a}\right) .$ (iii) Let $\Delta $ $\in {\cal Z}$ be the quadratic
Casimir operator which corresponds to the Killing form. Then $H_\Delta $ is
precisely the Toda Hamiltonian $H_T\ $ acting in ${\frak H}$; it is given by 
\[
H_T=-\left( \partial ,\partial \right) +\sum_{\alpha \in P}\exp 2\alpha
\left( q\right) ,\partial =\left( \frac \partial {\partial q_1},...,\frac
\partial {\partial q_n}\right) .
\]
\end{proposition}

In other words, geometric quantization of the Toda Hamiltonian agrees with
its `naive' Schroedinger quantization. What is nontrivial, of course, is the
consistent definition of the quantum integrals of motion which correspond to
higher Casimirs.

Let us now turn to the description of the Toda lattice based on the
reduction procedure. In section 4.1.1 we discussed the reduction using the
language of ${\cal U}\left( {\frak g}\right) $-modules; in the present
setting we may assume that all these modules are actually integrable, i.e.,
come from representations of the corresponding Lie groups. In particular,
the action of ${\cal U}\left( {\frak g}\right) $ on itself by left (right)
multiplications corresponds to the regular representation of $G;$ we may
regard ${\frak H}=L_2\left( G\right) $ as a result of quantization of $%
T^{*}G.$

We are led to the following construction. Let $\chi _f$ be the character of $%
N$ defined by 
\[
\chi _f\left( \exp X\right) =\exp if\left( X\right) ,X\in {\frak n.} 
\]
Let ${\frak H}_f$ be the space of smooth functions on $G$ satisfying the
functional equation 
\begin{equation}
\psi \left( kxn\right) =\chi _f\left( n\right) \psi \left( x\right) ,k\in
K,n\in N  \label{whit}
\end{equation}
By Iwasawa decomposition, any such function is uniquely determined by its
restriction to $A\subset G.$ Thus there is an isomorphism $i:C^\infty \left( 
{\frak a}\right) \rightarrow {\frak H}_f:i\left( \psi \right) \left(
kan\right) =\chi _f\left( n\right) \psi \left( \log a\right) .$

\begin{lemma}
${\frak H}_f$ is invariant with respect to the action of ${\cal Z}\subset 
{\cal U}\left( {\frak g}\right) .$
\end{lemma}

For $\zeta \in {\cal Z}$ let $\delta \left( \zeta \right) \in EndC^\infty
\left( {\frak a}\right) $ be its radial part defined by 
\[
\zeta i\left( \psi \right) =i\left( \delta \left( \zeta \right) \psi \right)
,\psi \in C^\infty \left( {\frak a}\right) . 
\]

\begin{proposition}
\label{rad}We have $\delta \left( \zeta \right) =H_\zeta ,$ i.e., the radial
parts of the Casimir operators coincide with the quantum integrals of the
Toda lattice.
\end{proposition}

We are interested in eigenfunctions of the Toda Hamiltonian and higher Toda
integrals. By proposition \ref{rad} this is equivalent to description of the
eigenfunctions of Casimir operators on \ \ $G$ which satisfy the functional
equation $\left( \ref{whit}\right) .$ At the formal level this problem may
be solved as follows. Let $\left( \pi ,{\frak V}\right) $ be an
infinitesimally irreducible representation of $G,$ ${\frak V}^{*}$, the dual
representation. Assume that $w\in {\frak V}$ satisfies 
\begin{equation}
\pi \left( n\right) w=\chi _f\left( n\right) w\text{ for }n\in N.
\label{whittaker}
\end{equation}
(In that case $w$ is called a {\em Whittaker vector}.) Assume, moreover,
that in the dual space ${\frak V}^{*}$ there is a $K$-invariant vector $v\in 
{\frak V}^{*}.$ Put 
\begin{equation}
\psi \left( x\right) =\left\langle v,\pi \left( x\right) w\right\rangle .
\label{wavef}
\end{equation}

\begin{proposition}
$\psi $ is an eigenfunction of ${\cal Z}$ and satisfies the functional
equation $\left( \ref{whit}\right). $
\end{proposition}

Thus $\psi $ is essentially a Toda lattice eigenfunction.

This formal argument can be made rigorous. Since the Toda Hamiltonian
describes scattering in a repulsive potential, its spectrum is continuous;
so there are no chances that $\psi $ (which is called a {\em Whittaker
function}) is in $L_2.$ The problem is to find representations of $G$ such
that $\psi $ has the usual properties of a continuous spectrum eigenfunction
(i.e., the wave packets smoothed down with appropriate amplitudes are in $%
L_2).$ As it appears, the correct class are {\em spherical principal series
representations.}

\begin{definition}
Let $B=MAN$ be the Borel subgroup in $G,$ $\lambda \in {\frak a}^{*};$ let $%
\chi _\lambda :B\rightarrow {\Bbb C}$ be the 1-dimensional representation
defined by 
\[
\chi _\lambda \left( man\right) =\exp \left\langle \lambda -\rho ,\log
a\right\rangle . 
\]
The spherical principal series $\pi _\lambda $ is the induced
representation, $\pi _\lambda =ind_B^G\ \chi _\lambda .$
\end{definition}

Principal series representations are infinitesimally irreducible and hence
define a homomorphism $\pi :{\cal Z}\rightarrow P\left( {\frak a}^{*}\right) 
$ into the algebra of polynomials on ${\frak a}^{*}$ 
\[
\pi :z\longmapsto \pi _z,\text{ where }\pi _\lambda \left( z\right) :=\pi
_z\ \left( \lambda \right) \cdot id. 
\]
By a classical Harish Chandra theorem (\cite{dixmier}) $\pi $ is actually an
isomorphism onto the subalgebra $P\left( {\frak a}^{*}\right) ^W\subset
P\left( {\frak a}^{*}\right) $ of the Weyl group invariants. Hence the set
of principal series representations is sufficiently ample to separate points
of the algebraic spectrum of ${\cal Z}.$

The algebraic theory of Whittaker vectors is exposed in \cite{kost-whit},
the analytic theory leading to the Plancherel theorem for the Toda lattice
is outlined in \cite{STSToda}. A remarkable point is that the eigenfunctions
of the quantum system are expressed as matrix coefficients of appropriate
representations of the 'hidden symmetry group'. This is a special case of a
very general situation. Below we shall discuss another example of this kind,
for which the hidden symmetry algebra is infinite-dimensional.

\subsubsection{\label{quant gaudin}4.1.2. Gaudin model.}

The treatment of the Gaudin model is considerably more difficult, since in
this case the hidden symmetry group is infinite-dimensional. Remarkably, the
general pattern described in section 4.1.1 is fully preserved in this case
as well.

The algebra of observables for the Gaudin model is simply ${\cal U}({\frak g}%
^N).$ Let $V_\lambda $ be a finite-dimensional highest weight representation
of ${\frak g}$ with dominant integral highest weight $\lambda .$ Let ${\bf %
\lambda }=(\lambda _1,...\lambda _N)$ be the set of such weights; put ${\bf V%
}_{{\bf \lambda }}=\otimes _iV_{\lambda _i}.$ The space ${\bf V}_{{\bf %
\lambda }}$ is a natural Hilbert space associated with the Gaudin model; so
the `kinematical' part of quantization problem in this case is fairly
simple. Let us fix also an auxiliary representation ($\rho ,V).$ By analogy
with the classical case, we may introduce the {\em quantum Lax operator}, 
\begin{equation}
L_V(z)\in \ EndV(z)\ \otimes {\frak g}^N\subset EndV(z)\otimes {\cal U}(%
{\frak g}^N);\   \label{Lax matr}
\end{equation}
the definition remains exactly the same, but now we embed ${\frak g}^N$ into
the universal enveloping algebra ${\cal U}({\frak g}^N)$ instead of the
symmetric algebra ${\cal S}({\frak g}^N).$ The commutation relations for $%
L_V(z)$ essentially reproduce the Poisson brackets relations $\left( \ref
{pbr}\right) $, but this time the l.h.s. is a matrix of true commutators: 
\begin{equation}
\left[ L_V(u)\otimes _{,}L_V(v)\right] =\left[ r_V(u,v),L_V(u)\otimes
1+1\otimes L_V(v)\right].  \label{comm}
\end{equation}
The key point in (\ref{comm}) is the interplay of commutation relations in
the quantum algebra ${\cal U}({\frak g}^N)$ and the auxiliary matrix algebra 
$EndV(z)$. Formula (\ref{comm}) was the starting point of QISM (as applied
to models with {\em linear} commutation relations). Put

\begin{equation}
S(u)=\frac 12tr_V\left( L_V(u)\right) ^2;  \label{tr}
\end{equation}
using (\ref{comm}), it is easy to check that $S(u)$ form a commutative family
of Hamiltonians (called {\em Gaudin Hamiltonians}) in ${\cal U}({\frak g}%
^N)( $see, e.g., \cite{Jurco}). An important property of this commuting
family is that it possesses at least one `obvious' eigenvector $\mid 0>\in H$%
, the tensor product of highest weight vectors in $V_{\lambda _i};$ it is
usually called {\em the vacuum vector .}

One of the key ideas of QISM is to construct other eigenvectors by applying
to the vacuum creation operators which are themselves rational functions of $%
z$. This construction is called{\em the algebraic Bethe Ansatz.} Assume that $%
{\frak g}={\frak sl}_2$ and let $\left\{ E,F,H\right\} $ be its standard
basis. Put 
\begin{equation}
F(z)=\sum_{z_i\in D}\frac{F^{(i)}}{z-z_i},  \label{creat}
\end{equation}
where $F^{(i)}$ acts as $F$ in the i-th copy of ${\frak sl}_2$ and as $id$
in other places$.$ The Bethe vector is, by definition, 
\begin{equation}
\mid w_1,w_2,...w_m>=F(w_1)F(w_2)...F(w_m)\mid 0>.  \label{bethe}
\end{equation}

The Lax matrix (\ref{Lax matr}) applied to $\mid w_1,w_2,...w_m>$ becomes
triangular, i.e., 
\[
L(u)\mid w_1,w_2,...w_m>=\left( 
\begin{array}{cc}
a(u,w_1,w_2,...w_m) & * \\ 
0 & d(u,w_1,w_2,...w_m)
\end{array}
\right) ; 
\]
after a short computation this yields 
\[
\begin{array}{l}
S(u)\mid w_1,w_2,...w_m>= \\ 
\\ 
\qquad s_m(u)\mid w_1,w_2,...w_m>+\sum_{j=1}^N\frac{f_j}{u-w_j}\mid
w_1,w_2,...,w_{j-1},u,w_{j+1},...,w_m>,
\end{array}
\]
where 
\[
f_j=\sum_{i=1}^N\frac{\lambda _i}{w_j-z_i}-\sum_{s\neq j}\frac 2{w_j-w_s} 
\]
and $s_m(u)$ is a rational function, 
\begin{equation}
\begin{array}{c}
s_m(u)=\frac{c_V}2\chi _m(u)^2-c_V\partial _u\chi _m(u), \\ 
\chi _m(u)=\sum_{i=1}^N\frac{\lambda _i}{u-z_i}-\sum_{i=1}^m\frac 2{u-w_j\ }.
\end{array}
\label{miura}
\end{equation}
(The constant $c_V$ depends on the choice of $V.$) If all $f_j$ vanish, $%
\mid w_1,w_2,...w_m>$ is an eigenvector of $S(u)$ with the eigenvalue $%
s_m(u);$ equations 
\begin{equation}
\sum_{i=1}^N\frac{\lambda _i}{w_j-z_i}-\sum_{s\neq j}\frac 2{w_j-w_s}=0
\label{betheq}
\end{equation}
are called the {\em Bethe Ansatz equations}. (Notice that (\ref{betheq}) is
precisely the condition that $s_m(u)$ be nonsingular at $u=w_i.)$

For general simple Lie algebras the study of spectra of the Gaudin
Hamiltonians becomes rather complicated; one way to solve this problem is to
treat it inductively by choosing in ${\frak g}$ a sequence of embedded Lie
subalgebras of lower rank and applying the algebraic Bethe Ansatz to these
subalgebras. An alternative idea (which is completely parallel to the
treatment of the Toda lattice in section \ref{quant toda}) is to interpret
the Hamiltonians as radial parts of (infinite dimensional) Casimir operator
of the `global' Lie algebra ${\frak g}_D.$

It is impossible to apply theorem \ref{K} immediately in the affine case,
since the center of ${\cal U}({\frak g}_D)$ is trivial. The point is that
the invariants in the (suitably completed) symmetric algebra ${\cal S}(%
{\frak g}_D)$ are infinite series; an attempt to quantize these expression
leads to divergent expressions, unless some kind of ordering prescription is
introduced. The commutation relations for the normally ordered operators are
already nontrivial and they do not lie in the center of ${\cal U}({\frak g}%
_D).$ However, the situation can be amended by first passing to the central
extension of ${\cal U}({\frak g}_D)$ and then considering the quotient
algebra ${\cal U}_k(\ {\frak g}_D)=U({\frak \hat g}_D)/(c-k).$ It is known
that for the{\em \ critical value }of the central charge $k=-h^{\vee }$
(here $h^{\vee }$ is the {\em dual Coxeter number} of ${\frak g}$) the
(appropriately completed) algebra ${\cal U}_{-h^{\vee }}({\frak g}_D)$
possesses an ample center.

Let us recall first of all the construction of Sugawara operators. Let $L%
{\frak g}={\frak g}\left( \left( z\right) \right) $ be the `local' algebra
of formal Laurent series with coefficients in ${\frak g,}$Let $\omega $ be
the 2-cocycle on $L{\frak g}$ defined by 
\[
\omega (X,Y)=Res_{z=0}\left\langle X,\partial _zY\right\rangle dz 
\]
It is well known that highest weight representations of $L{\frak g}$ are
actually projective and correspond to the central extension $\widehat{L%
{\frak g}}{\cal \simeq }$ $L{\frak g}\oplus {\Bbb C}$ $c$ of $L{\frak g}$
defined by this cocycle. Let us consider formally the canonical element $%
J\in L{\frak g}\otimes L{\frak g}$ which corresponds to the inner product 
\[
(X,Y)=Res_{z=0}\left\langle X,Y\right\rangle dz 
\]
in $L{\frak g}$ . Fix a finite-dimensional representation $\left( \rho
,V\right) $ of ${\frak g;}$ it extends to the `evaluation representation' of 
$L{\frak g}$ in $V\otimes _{{\Bbb C}}{\Bbb C}\left( \left( z\right) \right)
. $ Let $\left( \pi _k,{\frak V}\right) $ be any level $k$ highest weight
representation of $\widehat{L{\frak g}}.$  (This means that the central
element $c\in \widehat{L{\frak g}}$ acts by $\pi _k\left( c\right) =k\cdot
id.)$ Put $J(z)=\left( \pi _k\ \otimes \rho \right) J.$ Since $J$ is not a
proper element in $L{\frak g}\otimes L{\frak g}$ $,$ $J(z)$ is a formal
series infinite in both directions; however, its coefficients are well
defined. Namely, let $\left\{ e_a\right\} $ be a basis in ${\frak g}$ ; the
algebra $L{\frak g}$ is spanned by the vectors $e_a\left( n\right)
=e_a\otimes z^n,n\in {\Bbb Z}$. It is easy to see that the coefficients of $%
J(z)$ are {\em finite} linear combinations of the vectors $\pi _k\left(
e_a\left( n\right) .\right) $

Put 
\begin{equation}
T\left( u\right) =\frac 12:trJ\left( u\right) ^2:,\;T\left( u\right)
=\sum_{n=-\infty }^\infty T_nu^n,
\end{equation}
where the normal ordering {\bf ::} is defined in the following way: 
\begin{eqnarray*}
\ &:&e_a\left( n\right) e_b\left( m\right) :=e_{\ b}\left( m\right)
e_a\left( n\right) ,\text{ if }n<0,m\geq 0, \\
\ &:&e_a\left( n\right) e_b\left( m\right) :=e_a\left( n\right) e_b\left(
m\right) \ \text{otherwise.}
\end{eqnarray*}
Due to normal ordering the coefficients of $T\left( u\right) $ are well
defined operators in $End{\frak V.}$ It is well known that 
\begin{equation}
\begin{array}{l}
\left[ T_m,T_n\right] =\left( k+h^{\vee }\right) \left[ \left( m-n\right)
T_{m+n}+\frac{k\dim {\frak g}}{12}\left( n^3-n\right) \delta _{n,-m}\right] ,
\\ 
\lbrack T_m,\pi _k\left( e_a\left( n\right) \right) ]=-\left( k+h^{\vee
}\right) \pi _k\left( e_a\left( n+m\right) \right) .
\end{array}
\end{equation}
Therefore, if $k\neq -h^{\vee },$ $T_n$ generate the Virasoro algebra; for $%
k=-h^{\vee }$ the elements $T_n,$ $n\in {\Bbb Z},$ are central.

One can show that for ${\frak g}={\frak sl}_2$ the center of ${\cal U}%
_{-h^{\vee }}(\widehat{L{\frak g}})$ is generated by $T_n,$ $n\in {\Bbb Z},;$
for higher rank algebras there are other Casimirs. An attempt to construct
these higher Casimir elements by considering $:trJ\left( u\right) ^n:$ for
arbitrary $n\geq 2$ runs into trouble. However, by applying a different
technique \cite{FeiFr} have proved that for the critical central charge the
algebra of Casimirs is very ample: essentially, there exists a Casimir
element with prescribed symbol and hence there is an isomorphism ${\cal S}(L%
{\cal {\frak g}})^{L{\cal {\frak g}}}\rightarrow {\cal Z}\left( U_{-h^{\vee
}}\left( \widehat{L{\frak g}}\right) \right) .$

The situation with the Lie algebra ${\frak g}_D$ is basically the same as
described above. This allows to realize the Gaudin Hamiltonians as radial
parts of appropriate Casimir operators. Let us describe this construction
(due to \cite{FFR}) more precisely. Let the Lie algebras ${\frak g}_D,\;%
{\frak g}_D^{+},\;{\frak g}_D^{++},\;{\frak g}(D)$ be as above. It will be
convenient to add one more point $\left\{ u\right\} $ to the divisor $D$ and
to attach to it the trivial representation $V_0$ of ${\frak g}.$ (This will
not affect the Hilbert space of our model, since $\otimes _{z_i\in
D}V_{\lambda _i}\otimes V_0=\left( \otimes _{z_i\in D}V_{\lambda _i}\right)
\otimes {\Bbb C\simeq }\otimes _{z_i\in D}V_{\lambda _i}$.) Thus we write $%
D^{\prime }=D\cup \left\{ u\right\} ,$ etc.${\frak \ }$Let $\omega $ be the
2-cocycle on ${\frak g}_{D^{\prime }}$ defined by 
\begin{equation}
\omega (X,Y)=\sum_{z_i\in D^{\prime }}\ Res(X_i,dY_i).  \label{cocycl}
\end{equation}
Let ${\frak \hat g}_{D^{\prime }}={\frak g}_{D^{\prime }}\oplus {\Bbb C}c$
be the central extension of ${\frak g}_{D^{\prime }}$ defined by this
cocycle. Note that since the restriction of $\omega $ to the subalgebra $%
{\frak g}(D^{\prime })\subset {\frak g}_{D^{\prime }}$ is zero, the algebra $%
{\frak g}(D^{\prime })$ is canonically embedded into ${\frak \hat g}%
_{D^{\prime }}$. Put ${\frak \hat g}_{D^{\prime }}^{+}={\frak g}_{D^{\prime
}}^{+}\oplus {\Bbb C}c$. Fix a highest weight representation $V_{\left( {\bf %
\lambda },0\right) }=\otimes _{z_i\in D}V_{\lambda _i}\otimes V_0$ of the
Lie algebra ${\frak g}^{N+1}={\frak g}_{D^{\prime }}^{+}/{\frak g}%
_{D^{\prime }}^{++}$ as above and let $V_{({\bf \lambda ,}0{\bf )}}^k$ be
the associated representation of ${\frak \hat g}_{D^{\prime }}^{+}$ on which
the center ${\bf C}c$ acts by multiplication by $k\in {\Bbb Z}{\bf ..}$ Let $%
{\bf V}_{({\bf \lambda ,}0{\bf )}}^k$ be the induced representation of $%
{\frak \hat g}_{D^{\prime }},$ 
\begin{equation}
{\bf V}_{({\bf \lambda ,}0{\bf )}}^k={\cal U}({\frak \hat g}_D)\otimes _{%
{\cal U}({\frak \hat g}_{D^{\prime }}^{+})}V_{({\bf \lambda ,}0{\bf )}}^k.
\label{Fock}
\end{equation}
There is a canonical embedding $V_{({\bf \lambda },0)}^k\hookrightarrow {\bf %
V}_{({\bf \lambda ,}0)}^k:v\mapsto 1\otimes v.$ Let $\left( {\bf V}_{({\bf %
\lambda },0)}^k\right) ^{*}$ be the dual of ${\bf V}_{({\bf \lambda ,}0)}$, $%
{\bf H}_{({\bf \lambda },0)}^k\subset \left( {\bf V}_{({\bf \lambda ,}%
0)}^k\right) ^{*}$ the subspace of ${\frak g}(D^{\prime })-$invariants.
Decomposition (\ref{split}) together with the obvious isomorphism $V_{\left( 
{\bf \lambda },0\right) }\simeq V_{{\bf \lambda }}$ immediately implies that 
${\bf H}_{({\bf \lambda },0)}$ is canonically isomorphic to $V_{{\bf \lambda 
}}^{*}.$ Let $\Omega \in {\bf H}_{{\bf \lambda }}\otimes V_{{\bf \lambda }}$
be the canonical element; it defines a natural mapping ${\bf V}_{({\bf %
\lambda },0)}^{k,}\rightarrow V_{{\bf \lambda }}$:$\varphi \mapsto
\left\langle \varphi \right\rangle _\Omega .$ (In Conformal Field Theory $%
\left\langle \varphi \right\rangle _\Omega $ are usually called {\em %
correlation functions}.) \ For any $x\in U(\hat {{\frak g}}_D)$ let $\gamma
_{{\bf \lambda }}\left( x\right) \in EndV_{{\bf \lambda }}$ be the linear
operator defined by the composition mapping $v\mapsto \left\langle
x(1\otimes v)\right\rangle _\Omega ,v\in V_{{\bf \lambda .}}$

Now suppose that $k=-h^{\vee }$ ; let again $J$ be the canonical element
which corresponds to the inner product (\ref{res}) in ${\frak g}((z))$, $%
J(z)=(id\otimes \rho )T$ , $T(z)=tr:J(z)^2:$ Let us embed ${\cal U}%
_{-h^{\vee }}({\frak g}{\bf (}(z)))$ into ${\cal U}_{-h^{\vee }}({\frak g}%
_D) $ sending it to the extra place $\left\{ u\right\} \subset D^{\prime }$
to which we attached the trivial representation of ${\frak g}.$

\begin{proposition}
$S(u)=\gamma _\lambda (J_V(-2))$ coincides with the Gaudin Hamiltonian.
\end{proposition}

\begin{remark}
Besides quadratic Hamiltonians associated with the Sugawara current there
are also higher commuting Hamiltonians which may be obtained using the
methods of \cite{FeiFr}.
\end{remark}

In the calculation above we started with a generalized Verma module ${\bf V}%
_{({\bf \lambda ,}0{\bf )}}^k;$ however, any representation of the critical
level will do. The problem is to find a sufficiently ample class of critical
level representations which will account for the spectrum of the Gaudin
model. Note that this construction is exactly similar to the description of
the Toda eigenfunctions in the finite-dimensional case. In the Toda case the
correct class consisted of the spherical principal series representations;
an important point is that principal series representations separate points
of the algebraic spectrum of ${\cal Z}.$ Now, at the critical level the
algebra of the Casimir operators is extraordinary rich; hence generalized
Verma modules which are parametrized by the dual of the extended Cartan
subalgebra $\widehat{{\frak a}}={\frak a}\oplus {\Bbb C}$ (i.e. depend on a
finite number of parameters) can not be used. Remarkably, there is another
class of representations of ${\frak g}_D,$ the {\em Wakimoto modules}, which
play the role of the principal series representations. Roughly, the idea is
to use the loop algebra $L{\frak a}={\frak a}\otimes {\Bbb C}\left( \left(
z\right) \right) $ as the substitute of the Cartan subalgebra, to extend its
characters trivially to $L{\frak n}$ and to take the induced $L{\frak g}$
-module. Due to normal ordering (which is necessary, since after central
extension $L{\frak a}$ becomes a Heisenberg algebra), this construction
should be modified (among other things, this leads to a shift of the level,
i.e., of the value of the central charge, which becomes $-h^{\vee }$ ).

For ${\frak g}={\frak sl}_2$ the Wakimoto module of the critical level has
the following explicit realization (\cite{EFrenkel}). Let ${\sf H}$ be the
Heisenberg algebra with generators $a_n,a_n^{*},n\in {\Bbb Z},$ and
relations 
\begin{equation}
\left[ a_n,a_m^{*}\right] =\delta _{n,-m}.  \label{heisenberg}
\end{equation}
Let ${\frak F}$ be the Fock representation of \ ${\sf H}$ with vacuum vector 
$v$ satisfying $a_nv=0,n\geq 0,a_n^{*}v=0,n>0.$ Put 
\begin{equation}
a\left( u\right) =\sum_{n\in {\Bbb Z}}a_nu^{-n-1},a^{*}\left( u\right)
=\sum_{n\in {\Bbb Z}}a_n^{*}u^{-n}.  \label{local}
\end{equation}
Let $\left\{ E,F,H\right\} $ be the standard basis of ${\frak sl}_2.$
Put $E\left( n\right) =E\otimes z^n,F\left( n\right) $ $=F\otimes
z^n,H\left( n\right) =H\otimes z^n;$ let 
\[
E\left[ u\right] =\sum_{n\in {\Bbb Z}}\ E\left( n\right) u^n,F\left[
u\right] =\sum_{n\in {\Bbb Z}}\ F\left( n\right) u^n,H\left[ u\right]
=\sum_{n\in {\Bbb Z}}\ F\left( n\right) u^n 
\]
be the corresponding `generating functions'. For any formal power series $%
\chi \left( u\right) =\sum_{n\in {\Bbb Z}}\chi _nu^{-n-1}$ define the
(projective) action of $\ L{\frak g}$ on $F$ by 
\begin{eqnarray}
E\left[ u\right] &=&a\left( u\right) ,H\left[ u\right] =-2:a\left( u\right)
a^{*}\left( u\right) :+\chi \left( u\right) ,  \label{module} \\
F\left[ u\right] &=&-:a\left( u\right) a^{*}\left( u\right) a^{*}\left(
u\right) :-2\partial _ua^{*}\left( u\right) +\chi \left( u\right)
a^{*}\left( u\right) .  \nonumber
\end{eqnarray}
Thus we get an $L{\frak g}$-module structure in ${\frak F}$ which depends on 
$\chi \left( u\right) ;$ this is the Wakimoto module of the critical level $%
k=-2$ (usually denoted by $W_{\chi \left( u\right) }).$

The Virasoro generators $T_n,n\in {\Bbb Z},$ are scalar in $W_{\chi \left(
u\right) },T_n=q_n\cdot id.$ Put $q\left( u\right) =\sum_{n\in {\Bbb Z}}\
q_nu^n.$ One can show that 
\begin{equation}
q\left( u\right) =\frac 14\chi \left( u\right) ^2-\frac 12\partial _u\chi
\left( u\right) ,  \label{Miura}
\end{equation}
i.e., $q$ is the Miura transform of $\chi .$ This means of course that 
\[
\partial _u^2-q\left( u\right) =\left( \partial _u-\frac 12\partial _u\chi
\right) \left( \partial _u+\frac 12\partial _u\chi \right) ; 
\]
as a matter of fact, due to normal ordering, if we make a change of
coordinates on the (formal) disk, $q$ transforms as a projective connection.
Now a comparison with formula (\ref{miura}) suggests that the Bethe Ansatz
is related to a Wakimoto module with an appropriate choice of $\chi .$
Namely, we consider the following global rational function 
\[
\lambda \left( u\right) =\sum_{i=1}^N\frac{\lambda _i}{u-z_i}%
-\sum_{j=1}^m\frac 2{u-w_j}. 
\]
Denote by $\lambda _i\left( u-z_i\right) $ its expansion at the points $%
z_{i,}$ $i=1,...,N,$ and by $\mu _j\left( u-w_j\right) $ its expansion at
the points $w_j,$ $j=1,...,m.$ We have 
\[
\lambda _i\left( t\right) =\frac{\lambda _i}t+...,\mu _j\left( t\right)
=-\frac 2t+\mu _j\left( 0\right) +..., 
\]
where 
\[
\mu _j\left( 0\right) =\sum_{i=1}^N\frac{\lambda _i}{w_j-z_i}-\sum_{s\neq
j}\frac 2{w_j-w_s}. 
\]
We want to attach the Wakimoto modules to the points on the Riemann sphere.
To this end let us observe that creation and annihilation operators $\left( 
\ref{local}\right) $ define a projective representation of the `local'
algebra $\Gamma ={\Bbb C}\left( \left( u\right) \right) $ $\oplus {\Bbb C}%
\left( \left( u\right) \right) du$ consisting of formal series and formal
differentials with the cocycle 
\[
\Omega _0\left( \left( \varphi _1,\alpha _1\right) ,\left( \varphi _2,\alpha
_2\right) \right) =Res_{u=0}\left( \varphi _1\alpha _2-\varphi _2\alpha
_1\right) . 
\]
`Global' algebra is the direct sum of local algebras with the cocycle 
\[
\Omega \left( \left( \varphi _1,\alpha _1\right) ,\left( \varphi _2,\alpha
_2\right) \right) =\sum_iRes_{z_i}\left( \varphi _1\alpha _2-\varphi
_2\alpha _1\right) . 
\]
Now let us consider the tensor product of the Wakimoto modules attached to
the points $z_i,w_j,$%
\[
{\Bbb W}=\bigotimes_{i=1}^NW_{\lambda _i\left( t\right)
}\bigotimes_{j=1}^mW_{\mu _j\left( t\right) }. 
\]
Let ${\Bbb H}$ be the corresponding `big' Heisenberg algebra. The
eigenfunctions of the Gaudin Hamiltonians may be constructed as the
appropriate correlation functions associated with ${\Bbb W}$ . To formulate
the exact statement we need one more definition. Let $V$ be a ${\frak g}%
\left( \left( t\right) \right) $ -module. For $X\in {\frak g},n\in {\Bbb Z}$
we put $X\left( n\right) =X\otimes t^n\in {\frak g}\left( \left( t\right)
\right) .$ A vector $v\in V$ is called a {\em singular vector of imaginary
weight} if $X\left( n\right) v=0$ for all $X\in {\frak g},n\in {\Bbb Z}$ $%
,n\geq 0.$ Let $v_j\in W_{\mu _j\left( t\right) }$ be the vacuum vector; put 
$w_j=a_{-1}v_j,$ where $a_{-1}$ is the creation operator introduced in $%
\left( \ref{heisenberg}\right) .$

\begin{lemma}
$.$ The vector $w_j$ is a singular vector of imaginary weight if and only if 
$\mu _j\left( 0\right) =0.$
\end{lemma}

Note that this condition coincides with $\left( \ref{betheq}\right) .$

Let $\tilde M_{\lambda _i}$ be the subspace of $W_{\lambda _i\left( t\right)
}\simeq {\cal F}$ generated from the vacuum vector by the creation operators 
$a_0^{*}.$

\begin{lemma}
$\tilde M_{\lambda _i}$ is stable with respect to the subalgebra ${\frak g}%
\subset L{\frak g}$ of constant loops and is isomorphic to the dual of the
Verma module $M_{\lambda _i}$ over ${\frak g}={\frak sl}_2$ with the highest
weight $\lambda _i=Res_{t=0}\lambda _i\left( t\right) .$
\end{lemma}

This assertion is immediate from the definition $\left( \ref{module}\right) $
of the $\widehat{{\frak sl}_2}$ -action on ${\cal F}.$

Let ${\Bbb W}^{*}$ be the dual of ${\Bbb W}.$ . Let ${\frak h}_{z,w}\ \ $be
the algebra of rational functions with values in the Cartan subalgebra $%
{\frak h}\subset {\frak sl}_2$ which are regular outside $\left\{
z_1,...,z_N,w_1,...,w_m\right\} .$ There is a natural embedding of ${\frak h}%
_{z,w}$ into the `big' Heisenberg algebra $H$ defined via expansion at each
point of the divisor $\left\{ z_1,...,z_N,w_1,...,w_m\right\} .$ Let ${\cal H%
}_{z,w}$ be the maximal abelian subalgebra of ${\Bbb H}$ which contains $%
{\cal H}_{z,w}.$

\begin{lemma}
The space of ${\cal H}_{z,w}$-invariants in ${\Bbb W}^{*}$ is 1-dimensional;
it is generated by a functional $\tau $ whose value on the tensor product of
the vacuum vectors of $W_{\lambda _i\left( t\right) },W_{\mu _j\left(
t\right) }$ is equal to 1; the restriction of $\tau $ to $\otimes
_{i=1}^N\tilde M_{\lambda _i}\otimes w_1\otimes ...\otimes w_m$ is
nontrivial.
\end{lemma}

Let $\psi \in \otimes _{i=1}^N\tilde M_{\lambda _i}$. We may write $\tau
\left( \psi \right) =\left\langle \psi ,\phi \right\rangle $ where $\phi \in
\otimes _{i=1}^NM_{\lambda _i}$ is a vector in the tensor product of Verma
modules over ${\frak sl}_2.$

\begin{theorem}
(\cite{FFR}) If the Bethe equations $\left( \ref{betheq}\right) $ are
satisfied, $\phi \in \otimes _{i=1}^NM_{\lambda _i}$ is an eigenvector of
the Gaudin Hamiltonians in $\otimes _{i=1}^NM_{\lambda _i}$ with the
eigenvalue \ref{miura}.
\end{theorem}

\begin{remark}
If the weights $\lambda _i$ are dominant integral, there is a natural
projection 
\[
\pi :\otimes _{i=1}^NM_{\lambda _i}\rightarrow \otimes _{i=1}^NV_{\lambda
_i}. 
\]
It is easy to see that $\pi $ maps the eigenvectors of the Gaudin
Hamiltonians in $\otimes _{i=1}^NM_{\lambda _i}$ onto those in $\otimes
_{i=1}^NV_{\lambda _i}.$
\end{remark}

The construction described above admits a generalization to arbitrary
semisimple Lie algebras. We shall only write down the generalized Bethe
equations. To parametrize a Bethe vector first choose a set $\left\{
w_1,...,w_m\right\} $, $w_j\in {\bf C},$ and assign to each $w_j$ a set of
simple roots $\left\{ \alpha _{i_j}\right\} _{i=1}^N$, one for each copy of $%
{\frak g}$ in ${\frak g}^N$. Let $F_{i_j}^i$ be the corresponding Chevalley
generator of ${\frak g}$ acting nontrivially in the i-th copy of ${\frak g}$%
. The straightforward generalization of the Bethe creation operator (\ref
{creat}) to the higher rank case is 
\begin{equation}
F(w_j)=\sum_{i=1}^N\frac{F_{i_j}^i}{w_j-z_i}  \label{crh}
\end{equation}
(\cite{BabF}). The problem with (\ref{crh}) is that these operators no
longer commute with each other. Hence one cannot use a string of creation
operators to produce an eigenvector as it is done in $\left( \ref{bethe}%
\right) .$ The correct way to decouple them is to use the Wakimoto modules $%
W_{\mu _j\left( t\right) }$ which correspond to the poles of the creation
operator; as above, the eigenvectors are generated by the singular vectors
of imaginary weight which exist if and only if the constant term in the
expansion of $\mu _j\left( t\right) $ satisfies certain orthogonality
conditions. This leads to the {\em generalized Bethe equations}: 
\begin{equation}
\sum_{i=1}^N\frac{(\lambda _i,\alpha _{i_j})}{w_j-z_i}-\sum_{s\neq j}\frac{%
\left( \alpha _{i_s},\alpha _{i_j}\right) }{w_j-w_s}=0,j=1,...,m.
\label{bethh}
\end{equation}

Bethe vectors are again computed as correlation functions.

Let us end this section with a brief remark on the Knizhnik-Zamolodchikov
equations. Recall that this is a system of equations satisfied by the
correlation functions for an {\em arbitrary} value of the central charge;
the critical value $c=-h^{\vee }$ corresponds to the semiclassical limit for
the KZ system (small parameter before the derivatives). The outcome of this
is two-fold: first, the Bethe vectors for the Gaudin model appear naturally
in the semiclassical asymptotics of the solutions of the KZ system (\cite{RV}%
). Moreover, the exact integral representation of the solutions (for any
value of the central charge) also involves the Bethe vectors (\cite{FFR}).

\subsection{ Quadratic case. Quasitriangular Hopf Algebras}

For an expert in Quantum Integrability the Gaudin model is certainly a sort of
a limiting special case.\ The real thing starts with the quantization of 
{\em quadratic} Poisson bracket relations (\ref{skl}). This is a much more
complicated problem which eventually requires the whole machinery of Quantum
Groups (and has led to their discovery). The substitute for the Poisson
bracket relations $\left( \ref{affskl}\right) $ is the famous relation 
\begin{equation}
R(u\ v^{-1})L^1(u)L^2(v)R(u\ v^{-1})^{-1}=L^2(v)L^1(u),  \label{RLL}
\end{equation}
where $R(u)$ is the quantum R-matrix satisfying the quantum Yang-Baxter
identity 
\begin{equation}
R_{12}(u)R_{13}(u\ v)R_{23}(v)=R_{23}(v)R_{13}(u\ v)R_{12}(u).  \label{YBE}
\end{equation}
To bring a quantum mechanical system into Lax form one has to arrange
quantum observables into a Lax matrix $L(u)$ (which is a rational function
of $u$) and to find an appropriate R-matrix satisfying (\ref{RLL}), (\ref
{YBE}). The first examples of quantum Lax operators were constructed by trial
and error method; in combination with the Bethe Ansatz technique this has
led to the explicit solution of important problems (\cite{STF}, \cite{FT},
Faddeev (\cite{Fadd}, \cite{Leshouches})).

The algebraic concept which brings order to the subject is that of {\em %
quasitriangular Hopf algebra} (\cite{QG}). Main examples of quasitriangular
Hopf algebras arise as q-deformations of universal enveloping algebras
associated with Manin triples. Remarkably, the general pattern represented
by Theorems \ref{AKS}, \ref{K} survives q-deformation. The standard way to
describe quantum deformations of simple finite-dimensional or affine Lie
algebras is by means of generators and relations generalizing the classical
Chevalley-Serre relations (\cite{QG}, \cite{Jimbo}). We shall recall this
definition below in section 4.2.2 in the finite-dimensional case, and in
section 4.2.3 for the quantized universal enveloping algebra of the loop
algebra $L\left( {\frak sl}_2\right) .$ A dual approach, due to \cite{FRT},
is to construct quantum universal enveloping algebras as deformations of
coordinate rings on Lie groups (regarded as linear algebraic groups). Of
course, construction of a quantum deformation of the Poisson algebra ${\cal F%
}\left( G\right) $ was one of the first results of the quantum group theory
and is, in fact, a direct generalization of the Baxter commutation relations 
$RT^1T^2=T^2T^1R.$ A nontrivial fact, first observed by Faddeev, Reshetikhin
and Takhtajan, is that the dual Hopf algebra (usually described as a
q-deformation of the universal enveloping algebra) may also be regarded as a
deformation of a Poisson algebra of functions on the dual group. More
generally, the FRT construction is related to the {\em quantum duality
principle}, which we shall now briefly discuss.

Let $\left( {\frak g},{\frak g}^{*}\right) $ be a factorizable Lie
bialgebra, $G,G^{*}$ the corresponding dual Poisson groups, ${\cal F}\left(
G\right) ,{\cal F}\left( G^{*}\right) $ the associated Poisson-Hopf algebras
of functions on $G,G^{*},$ and ${\cal F}_q\left( G\right) ,{\cal F}_q\left(
G^{*}\right) $ their quantum deformations. (For simplicity we choose the
deformation parameter $q$ to be the same for both algebras.) The quantum
duality principle asserts that these algebras are dual to each other as Hopf
algebras. More precisely, there exists a nondegenerate bilinear pairing 
\[
{\cal F}_q\left( G\right) \otimes {\cal F}_q\left( G^{*}\right) \rightarrow 
{\Bbb C}\left[ \left[ q\right] \right] 
\]
which sets the algebras ${\cal F}_q\left( G\right) ,{\cal F}_q\left(
G^{*}\right) $ into duality as Hopf algebras. Hence, in particular, we have,
up to an appropriate completion, 
\[
{\cal F}_q\left( G^{*}\right) \simeq {\cal U}_q\left( {\frak g}\right) . 
\]
In the dual way, we have also 
\[
{\cal F}_q\left( G\right) \simeq {\cal U}_q\left( {\frak g}^{*}\right) . 
\]
For factorizable Lie bialgebras the quantum deformations ${\cal F}_q\left(
G\right) ,{\cal F}_q\left( G^{*}\right) $ may easily be constructed once we
know the corresponding quantum R-matrices. An equivalence of this
formulation to the definition of Drinfeld and Jimbo is not immediate and
requires the full theory of universal R-matrices. Namely, starting with
Drinfeld's definition of a quasitriangular Hopf algebra we may construct
`quantum Lax operators' whose matrix coefficients generate the quantized
algebras of functions ${\cal F}_q\left( G\right) ,{\cal F}_q\left(
G^{*}\right) .$ Using explicit formulae for universal R-matrices one can, in
principle, express these generators in terms of the Drinfeld-Jimbo
generators. In the context of integrable models the FRT formulation has
several important advantages: it allows to state the quantum counterpart of
the main commutativity theorem as well as a transparent correspondence
between classical and quantum integrable systems. The analogue of the FRT
realization in the affine case is nontrivial, the key point being the
correct treatment of the central element which corresponds to the central
extension; it was described by \cite{RS}.

\subsubsection{4.2.1. Factorizable Hopf algebras and the
Faddeev-Reshetikhin-Takhtajan realization of quantized universal enveloping
algebras.}

\begin{definition}
Let $A$ be a Hopf algebra with coproduct $\Delta $ and antipode $S$; let $%
\Delta ^{\prime }$ be the opposite coproduct in $A$; $A$ is called
quasitriangular if 
\begin{equation}
\Delta ^{\prime }(x)={\cal R}\Delta (x){\cal R}^{-1}  \label{R}
\end{equation}
for all $x\in A$ and for some distinguished invertible element ${\cal R}\in
A\otimes A$ ({\em the universal R-matrix}) and, moreover, 
\begin{equation}
\left( \Delta \otimes id\right) {\cal R}={\cal R}_{13}{\cal R}_{23},\left(
id\otimes \Delta \right) {\cal R}={\cal R}_{13}{\cal R}_{12}.  \label{qt}
\end{equation}
Identities (\ref{qt}) imply that ${\cal R}$ satisfies the Yang-Baxter
identity 
\begin{equation}
{\cal R}_{12}{\cal R}_{13}{\cal R}_{23}={\cal R}_{23}{\cal R}_{13}{\cal R}%
_{12}.  \label{YB}
\end{equation}
\end{definition}

(We use standard tensor notation to denote different copies of the spaces
concerned.)

Let $A^0$ be the dual Hopf algebra equipped with the opposite coproduct (in
other words, its coproduct is dual to the opposite product in $A$). Put $%
{\cal R}^{+}={\cal R},{\cal R}^{-}=\sigma ({\cal R}^{-1})$ (here $\sigma $
is the permutation map in $A\otimes A,$ $\sigma (x\otimes y)=y\otimes x$)
and define the mappings $R_{\pm }:A^0\rightarrow A:f\mapsto \left\langle
f\otimes id,{\cal R}^{\pm }\right\rangle ;$ axioms (\ref{qt}) imply that $%
R_{\pm }$ are Hopf algebra homomorphisms. Define the action $A^0\otimes
A\rightarrow A$ by 
\begin{equation}
f\circ x=\sum_iR^{+}(f_i^{(1)})\;x\;S(R^{-}(f_i^{(2)}),\text{\ where}%
\;\Delta ^0f=\sum_if_i^{(1)}\otimes f_i^{(2)}.  \label{fact}
\end{equation}

\begin{definition}
$A$ is called {\em factorizable} if $A$ is a free $A^0$-module generated by $%
1\in A$.
\end{definition}

(Let us denote the corresponding linear isomorphism $A^0\rightarrow A$ by $F$
for future reference.)

There are several important examples of factorizable Hopf algebras:

\begin{itemize}
\item  Let\hfill $A$ \hfill be\hfill an \hfill arbitrary \hfill Hopf \hfill %
algebra; \hfill then \hfill its \hfill Drinfeld \hfill double \hfill$D\left(
A\right) $ \hfill is \hfill factorizable \\(\cite{RSfact}).

\item  Let ${\frak g}$ be a finite-dimensional semisimple Lie algebra; then $%
A={\cal U}_q\left( {\frak g}\right) $ is factorizable.

\item  Let ${\frak \hat g}$ be an affine Lie algebra; then $A={\cal U}%
_q\left( {\frak \hat g}\right) $ is factorizable.
\end{itemize}

(Observe that the last two cases are actually special cases of the
first one; up to now, the double remains the principal (if not the only)
source of factorizable Hopf algebras.)

Let ${\frak g}$ be a simple Lie algebra. Let $A={\cal U}_q({\frak g})$ be
the corresponding quantized enveloping algebra, $A^0$ be the dual of $A$
with the opposite coproduct. Let ${\cal R}$ be the universal R-matrix of $A.$
Let $V,W$ be finite-dimensional irreducible representations of $U_q({\frak g}%
).$ Let ${\cal L}\in A\otimes A^0$ be the canonical element. Put 
\begin{equation}
L_V=(\rho _{V\ }\ \otimes id){\cal L},R_{\pm }^{VW}\ =(\rho _{V\ }\otimes
\rho _W\ ){\cal R}_{\pm }.  \label{qlax}
\end{equation}
We may call $L_V\in EndV\ \otimes A^0$ the {\em universal quantum Lax
operator} (with auxiliary space $V$). Property (\ref{R}) immediately implies
that 
\begin{equation}
L_W^2\ L_V^1\ =R_{+}^{VW}\ L_V^1\ L_W^2\ \left( R_{+}^{VW}\ \right) ^{-1}.
\label{absRLL}
\end{equation}

\begin{proposition}
The associative algebra ${\cal F}_q\left( G\right) $ generated by the matrix
coefficients of $L_V$ satisfying the commutation relations $\left( \ref
{absRLL}\right) $ and with the matrix coproduct 
\[
\Delta L_V=L_V\stackrel{.}{\otimes }L_V
\]
is a quantization of the Poisson-Hopf algebra ${\cal F}\left( G\right) $
with the Poisson bracket $\left( \ref{ASkl}\right) ,\left( \ref{skl12}%
\right) .$
\end{proposition}

The dual algebra ${\cal F}_q\left( G^{*}\right) $ is described in the
following way. Let again ${\cal L}\in A\otimes A^0$ be the canonical
element. Put 
\begin{equation}
\begin{array}{c}
T_V\ ^{\pm }=\left( \ \rho _V\otimes R^{\pm }\right) {\cal L}\ \in EndV\
\otimes A, \\ 
T_V\ =T_V{}^{+}\;\left( id\otimes S\right) T_V\ ^{-}.
\end{array}
\label{Ldual}
\end{equation}

\begin{proposition}
(i) The matrix coefficients of $T_V\ ^{\pm }$ satisfy the following
commutation relations: 
\begin{equation}
\begin{tabular}{l}
$T_W^{\pm 2}\ \ T_V^{\pm 1}=R^{VW}\ T_V^{\pm 1}\ T_W^{\pm 2}\ \left( R^{VW}\
\right) ^{-1},$ \\ 
$T_W^{-2}T_V^{+1}=R^{VW}T_V^{+1}T_W^{-2}\left( R^{VW}\right) ^{-1}.j$ \\ 
\end{tabular}
\label{RTTdual}
\end{equation}
(ii) The associative algebra ${\cal F}_q\left( G^{*}\right) $ generated by
the matrix coefficients of $T^{\pm V}$ satisfying the commutation relations $%
\left( \ref{RTTdual}\right) $ and with the matrix coproduct 
\[
\Delta T^{\pm V}=T^{\pm V}\stackrel{.}{\otimes }T^{\pm V}
\]
is a quantization of the Poisson-Hopf algebra ${\cal F}\left( G^{*}\right) $
with Poisson bracket ($\ref{dual}),\left( \ref{dual12}\right) .$ \\ %
(iii) $T_V=(id\otimes F)L_V.$ The matrix coefficients of $T_V$ satisfy the
following commutation relations: 
\begin{equation}
\left( R_{+}^{VW}\right) ^{-1}T_W^2R_{+}^{VW}T_V^1=T_V^1\left(
R_{-}^{VW}\right) ^{-1}T_W^2R_{-}^{VW}.  \label{RTRT}
\end{equation}
(iv) The associative algebra generated by the matrix coefficients of $T_V$
satisfying the commutation relations $\left( \ref{RTTdual}\right) $ is a
quantization of the Poisson-Hopf algebra ${\cal F}\left( G^{*}\right) $ with
Poisson bracket ( \ref{dual12}). \\%
 (v) The pairing 
\[
\left\langle T_V^{\pm 1},L_W^2\right\rangle =R_{+}^{VW}
\]
sets the algebras ${\cal F}_q\left( G\right) ,$ ${\cal F}_q\left(
G^{*}\right) $ into duality as Hopf algebras.
\end{proposition}

Formulae (\ref{absRLL}, \ref{RTTdual}, \ref{RTRT}) are the exact quantum
analogues of the Poisson bracket relations (\ref{skl12}, \ref{dual}, \ref
{dual12}), respectively.

\subsubsection{4.2.2. Quantum commutativity theorem and quantum Casimirs.}

In complete analogy with the linear case the appropriate algebra of
observables which is associated with quantum Lax equations is not the
quasitriangular Hopf algebra $A$ but rather its dual \ $A^0$; the
Hamiltonians arise from the Casimir elements of $A$. We may summarize this
picture in the following heuristic correspondence principle. 
\[
\begin{tabular}{|cc|}
\hline
\multicolumn{1}{|c|}{\bf Quadratic Classical Case} & {\bf Quadratic Quantum
case} \\ \hline
\multicolumn{1}{|c|}{${\cal F}(G)$} & $A^0={\cal U}_q\left( {\frak g}%
^{*}\right) \simeq {\cal F}_q\left( G\right) $ \\ \hline
\multicolumn{1}{|c|}{${\cal F}\left( G^{*}\right) $} & $A={\cal U}_q\left( 
{\frak g}\right) \simeq {\cal F}_q\left( G^{*}\right) $ \\ \hline
\multicolumn{1}{|c|}{\em Classical r-matrices} & {\em Quantum R-matrices} \\ 
\hline
\multicolumn{1}{|c|}{{\em Casimir functions }{\em in} ${\cal F}\left(
G^{*}\right) $} & {\em Casimir operators} ${\cal Z}\subset {\cal U}_q\left( 
{\frak g}\right) $ \\ \hline
\multicolumn{1}{|c|}{{\em Symplectic leaves in} $G$} & {\em Irreducible
representations of} $A^0$ \\ \hline
\multicolumn{1}{|c|}{{\em Symplectic leaves in} $G^{*}$} & {\em Irreducible
representations of }$A$ \\ \hline
\end{tabular}
\]

According to this correspondence principle, in order to get a `quantum Lax
system' associated with a given factorizable Hopf algebra $A$ we may proceed
as follows:

\begin{itemize}
\item  {\em Choose a representation }$\left( \pi ,{\frak V}\right) ${\em \
of }$A^0.$

\item  {\em Choose a Casimir element }$\zeta \in {\cal Z}\left( A\right) $%
{\em \ and compute its inverse image }$F^{-1}\left( \zeta \right) \in A^0$%
{\em \ with respect to the factorization map }$F:A^0\rightarrow A.$

\item  {\em Put }${\cal H}_\zeta =\pi \left( F^{-1}\left( \zeta \right)
\right) .$
\end{itemize}

Moreover, we expect that the eigenvectors of the quantum Hamiltonian ${\cal H%
}_\zeta $  can be expressed as matrix coefficients of appropriate
representations of $A.$

In order to describe quantum Casimirs explicitly we have to take into
account the effect of twisting. We shall start with their description in the
finite-dimensional case. The results stated below are the exact quantum
counterparts of those described in section 3.2.2.

Let ${\frak g}$ be a finite-dimensional simple Lie algebra, $A={\cal U}_q(%
{\frak g})$ the quantized universal enveloping algebra of ${\frak g}.$ We
shall recall its standard definition in terms of the Chevalley generators (%
\cite{QG}, \cite{Jimbo}). Let $P$ be the set of simple roots of ${\frak g}.$
For $\alpha _i\in P$ we set $q_i=q^{\left\langle \alpha _i,\alpha
_i\right\rangle }.$ We denote the Cartan matrix of ${\frak g}.$
by $A_{ij}$
\begin{definition}
\label{UQ}${\cal U}_q({\frak g})$ is a free associative algebra with the
generators $\left\{ k_i,k_i^{-1},e_i,f_i\right\} _{i\in P}$ and relations 
\begin{eqnarray*}
k_i\cdot k_i^{-1}=k_i^{-1}k_i=1,\;\left[ k_i,k_j\right] =0, \\
\left[ a_i,e_j\right] =q_i^{A_{ij}}e_j,\;\left[ a_i,f_j\right]
=q_i^{-A_{ij}}e_j, \\
\left[ e_i,f_j\right] =\delta _{ij}\frac{\left( k_i^2-k_i^{-2}\right) }{%
q_i^2-q_i^{-2}}
\end{eqnarray*}
and the following q-Serre relations: 
\begin{eqnarray*}
\sum_{n=0}^{1-A_{ij}}(-1)^n\QATOPD[ ]
{1-A_{ij}}{n}_{q_i^2}e_i^{1-A_{ij}-n}e_je_i^n=0, \\
\sum_{n=0}^{1-A_{ij}}(-1)^n\QATOPD[ ]
{1-A_{ij}}{n}_{q_i^2}f_i^{1-A_{ij}-n}f_jf_i^n=0.
\end{eqnarray*}
(Here $\QATOPD[ ] {m}{n}_q$ are the q-binomial coefficients.)
\end{definition}

\begin{theorem}
(\cite{QG}) ${\cal U}_q({\frak g})$ is a quasitriangular Hopf algebra.
\end{theorem}

As usual, we denote by ${\cal R}\in {\cal U}_q({\frak g})\otimes {\cal U}_q(%
{\frak g})$ its universal R-matrix and put $R_{VW}=\left( \rho _{V\otimes
}\rho _W\right) {\cal R}$

We recall the following well known corollary of definition \label{Uq}\ref{UQ}%
. Let $H\ $ be the `Cartan subgroup' generated by the elements $k_i,i\in
P, $ ${\Bbb C}\left[ H\right] $ its group algebra. We denote by ${\cal U}_q(%
{\frak n}_{\pm }),$ ${\cal U}_q({\frak b}_{\pm })\subset {\cal U}_q({\frak g}%
)\ $the subalgebras generated by $e_i,f_i$ (respectively, by $k_i,e_i$ and
by $k_i,f_i).$

\begin{proposition}
(i) ${\Bbb C}\left[ H\right] ,$ ${\cal U}_q({\frak n}_{\pm }),{\cal U}_q(%
{\frak b}_{\pm })$ are Hopf subalgebras in ${\cal U}_q({\frak g}).$(ii) $%
{\cal U}_q({\frak n}_{\pm })$ is a two-sided Hopf ideal in ${\cal U}_q(%
{\frak b}_{\pm })$ and ${\cal U}_q({\frak b}_{\pm })/{\cal U}_q({\frak n}%
_{\pm })\simeq {\Bbb C}\left[ H\right] .$
\end{proposition}

Let $\pi _{\pm }:{\cal U}_q({\frak b}_{\pm })\rightarrow {\Bbb C}\left[
H\right] $ be the canonical projection. Consider the embedding 
\[
i:{\Bbb C}\left[ H\right] \rightarrow {\Bbb C}\left[ H\right] \otimes {\Bbb C%
}\left[ H\right] :a\longmapsto \left( id\otimes S\right) \Delta a.
\]
(Here $\Delta $ is the coproduct in commutative and cocommutative Hopf
algebra ${\Bbb C}\left[ H\right] $ and $S$ is its antipode.)

\begin{proposition}
The dual of ${\cal U}_q({\frak g})$ may be identified with the subalgebra 
\[
A^0=\left\{ x\in {\cal U}_q({\frak b}_{+})\otimes {\cal U}_q({\frak b}%
_{-});\pi _{+}\otimes \pi _{-}(x)\in i\left( {\Bbb C}\left[ H\right] \right)
\right\} .
\]
\end{proposition}

\begin{lemma}
For any $h\in H$,  $h\otimes h{\cal R}={\cal R}h\otimes h.$
\end{lemma}

The algebra $A$ admits a (trivial) family of deformations $A_h={\cal U}_q(%
{\frak g})_h,h\in H,$ defined by the following prescription: 
\begin{equation}
e_i\longmapsto he_i,f_i\longmapsto hf_i,k_i\longmapsto hk_i.  \label{deform}
\end{equation}
It is instructive to look at the commutation relations of $A_h$ in the FRT
realization. Let $V,W\in {\sf Rep}A;$ we define the `generating matrices' $%
T_V^{\pm }\in A\otimes EndV,T_W^{\pm }\in A\otimes EndW$, as in (\ref{Ldual}).

\begin{proposition}
(i) The deformation (\ref{deform}) transforms $T_V^{\pm }$ into $\rho
_V(h^{\mp 1})T_V^{\pm };$ the commutation relations in the deformed algebra $%
A_h$ amount to 
\begin{equation}
\begin{tabular}{l}
$T_W^{\pm 2}\ \ T_V^{\pm 1}=R^{VW}\ T_V^{\pm 1}\ T_V^{\pm 2}\ \left( R^{VW}\
\right)^{-1} ,$ \\
$T_W^{-2}T_V^{+1}=R_h^{VW}T_V^{+1}T_W^{-2}\left( R^{VW}\right) ^{-1},$%
\end{tabular}
\label{RTThdual}
\end{equation}
where \ $R_h^{VW}=\rho _V\otimes \rho _W\left( {\cal R}_h\right) $ and $%
{\cal R}_h=h\otimes h^{-1}{\cal R}h^{-1}\otimes h.$ \\%
 (ii) Put $T_h=\rho
_V\left( h^2\right) T_V^{+}\left( T_V^{-}\right) ^{-1};$ the
`operator-valued matrix' $T_h$ satisfies the following commutation relations 
\begin{equation}
\left( R_{+}^{VW}\right) ^{-1}T_W^2R_{+}^{VW}T_V^1=T_V^1\left(
R_{-}^{VW}\right) ^{-1}T_W^2R_{-}^{VW}.  \label{hRTRT}
\end{equation}
\end{proposition}

Let $L_V\in A^{*}\otimes EndV$ be the `universal Lax operator' introduced in
(\ref{qlax}). Fix $h\in H$ and put 
\begin{equation}
l_V^h\ =tr_V\rho _V(h^2)L^V\ .  \label{transfer}
\end{equation}
Elements $l_V^h\ \in A^0$ are usually called (twisted) {\em transfer
matrices. }The following theorem is one of the key results of QISM (as
applied to the finite-dimensional case).

\begin{theorem}
For any representations $V,W$ and for any$\ h\in H\ $, 
\begin{equation}
l_V^h\ l_W^h\ =l_W^h\ l_V^h\ .  \label{Comm}
\end{equation}
\end{theorem}

The mutual commutativity of transfer matrices is a direct corollary of the
commutation relations ($\ref{absRLL});$ twisting by elements $h\in H$ is
compatible with these relations, since 
\[
\left( h\otimes h\right) {\cal R}={\cal R}\left( h\otimes h\right) . 
\]

We want to establish a relation between the transfer matrices and the
Casimirs of ${\cal U}_q({\frak g})_{h^{\prime }}$ (where $h^{\prime }\in H$
may be different from $h).$ To account for twisting we must slightly modify
the factorization map.

Let $H\times A\rightarrow A$ be the natural action of $H$ on $A\ $ by right
translations, and let $H\ \times A^0\rightarrow A^0$ be the contragredient action. The
algebra $A={\cal U}_q({\frak g})$ is factorizable; let $F:A^0\rightarrow A$
be the isomorphism induced by the action (\ref{fact}); for $h\in H$ put $%
F^h=F\circ h.$ Let us compute the image of the universal Lax operator $L_V\
\in EndV\ \otimes A^0$ under the `factorization mapping' $\left( id\otimes
F^h\right) ;$ it is easy to see that 
\[
\left( id\otimes F^h\right) (L_V):=T_V^{h\pm }=\left( id\otimes h^{\pm
1}\right) T_V^{1\pm } 
\]
Put 
\begin{equation}
t_V^h=tr_VT_V^{h+}\left( T_V^{h-}\right) ^{-1}=tr_V\rho _V(h)^2T_V^1;
\label{qtr}
\end{equation}
clearly, we have $t_V^h=F(l_V^h).$

\begin{theorem}
\label{Qk}(i) (\cite{FRT}) Suppose that $h=q^{-\rho }$ , where $2\rho \in 
{\frak h}$ is the sum of positive roots of ${\frak g}$ . Then all
coefficients of $t_V^h$ are central in $A={\cal U}_q\left( {\frak g}\right) .
$ (ii) The center of $A\ $ is generated by $\left\{ t_V^h\right\} $ (with $V$
ranging over all irreducible finite-dimensional representations of ${\cal U}%
_q\left( {\frak g}\right) )$.(iii) For any $s\in H$ we have $%
l_V^s=F^{sh^{-1}}(t_V^h).$
\end{theorem}

Theorem \ref{Qk} is an exact analogue of theorem \ref{K} for
finite-dimensional factorizable quantum groups. Its use is twofold: first,
it provides us with commuting quantum Hamiltonians; second, the
eigenfunctions of these Hamiltonians may be constructed as appropriate
matrix coefficients (`correlation functions') of irreducible
representations of ${\cal U}_q\left( {\frak g}\right) .$ Again, to apply the
theorem to more interesting and realistic examples we must generalize it to
the affine case.

\subsubsection{4.2.3. Quantized affine Lie algebras.}

Two important classes of infinite-dimensional quasitriangular Hopf algebras
are quantum affine Lie algebras and the full Yangians (that is, the doubles
of the Yangians defined in (\cite{QG})). For concreteness I shall consider
only the first class. To avoid technical difficulties we shall consider the
simplest nontrivial case, that of the affine ${\frak sl}_2.$

The standard definition of ${\cal U}_q\widehat{({\frak sl}_2)}$ is by means
of generators and relations:

\begin{definition}
${\cal U}_q\widehat{({\frak sl}_2)}$ is a free associative algebra over $%
{\Bbb C}\left[ q,q^{-1}\right] $ with generators $e_i,f_i,K_i,$ $K_i^{-1},$ $%
i=0,1,$ which satisfy the following relations: 
\begin{eqnarray}
K_iK_j &=&K_jK_i,  \nonumber \\
K_ie_j &=&q^{A_{ij}}e_jK_i,\;K_if_j=q^{-A_{ij}}f_jK_i,  \label{rel} \\
\left[ e_i,f_j\right]  &=&\delta _{ij}\left( q-q^{-1}\right) \left(
K_i-K_i^{-1}\right) ,  \nonumber
\end{eqnarray}
where $A_{ij}$ is the Cartan matrix of the affine ${\frak sl}%
_2,\;A_{00}=A_{11}=-A_{01}=-A_{10}=2.$ In addition to (\ref{rel}) the
following {\em q-Serre} relations are imposed: 
\begin{eqnarray}
e_i^3e_j-\left( 1+q+q^{-1}\right) \left( e_i^2e_je_i-e_ie_je_i^2\right)
-e_je_i^3 &=&0,  \label{serre} \\
f_i^3f_j-\left( 1+q+q^{-1}\right) \left( f_i^2f_jf_i-f_if_jf_i^2\right)
-f_jf_i^3 &=&0.  \nonumber
\end{eqnarray}
The Hopf structure in ${\cal U}_q\widehat{({\frak sl}_2)}$ is defined by 
\begin{eqnarray}
\Delta K_i &=&K_i\otimes K_i,  \nonumber \\
\Delta e_i &=&e_i\otimes 1+K_i\otimes e_i, \\
\Delta f_i &=&f_i\otimes 1+K_i\otimes e_i.  \nonumber
\end{eqnarray}
The element $K_0K_1\ is$ central in ${\cal U}_q\widehat{({\frak sl}_2)};$
the quotient of ${\cal U}_q\widehat{({\frak sl}_2)}$ over the relation 
\[
K_0K_1=q^k,k\in {\Bbb Z},
\]
is called the level $k$ quantized universal enveloping algebra ${\cal U}_q%
\widehat{({\frak sl}_2)}_k.$
\end{definition}

An alternative realization of ${\cal U}_q\widehat{({\frak sl}_2)}_k$ is
based on the explicit use of the quantum R-matrix; in agreement with the
quantum duality principle, this realization may be regarded as an explicit
quantization of the Poisson algebra of functions ${\cal F}_q\left(
G^{*}\right) $, where $G^{*}$ is the {\em Poisson dual} of the loop group $%
L(SL_2).$ The relevant Poisson structure on $G^{*}$ has to be {\em twisted}
in so as to take into account the central charge $k;$ moreover, in complete
analogy with the non-deformed case, the center of ${\cal U}_q\widehat{(%
{\frak sl}_2)}_k$ is nontrivial only for the critical value of the central
charge $k_{crtit}=-2.$ (The critical value $k=-2$ is a specialization, for $%
{\frak g}={\frak sl}_2,$ of the general formula $k=-h^{\vee }$ ; thus the
value of $k_{crit}$ is not affected by q-deformation.)

Introduce the R-matrix 
\begin{equation}
R(z)=\left( 
\begin{array}{cccc}
1 & 0 & 0 & 0 \\ 
0 & \frac{1-z}{q-zq^{-1}} & \frac{z\left( q-q^{-1}\right) }{q-zq^{-1}} & 0
\\ 
0 & \frac{q-q^{-1}}{q-zq^{-1}} & \frac{1-z}{q-zq^{-1}} & 0 \\ 
0 & 0 & 0 & 1
\end{array}
\right) .  \label{Rq}
\end{equation}
Let ${\cal U}_q^{\prime }\widehat{({\frak gl}_2)}_k$ be an associative
algebra with generators $l_{ij}^{\pm }\left[ n\right] ,i,j=1,2,n\in \pm 
{\Bbb Z}$ $_{+}$. In order to describe commutation relations in ${\cal U}%
_q^{\prime }\widehat{({\frak gl}_2)}_k$ we introduce the generating series $%
L^{\pm }(z)=\left\| l_{ij}^{\pm }\left( z\right) \right\| ,$%
\begin{equation}
l_{ij}^{\pm }\left( z\right) =\sum_{n=0}^\infty l_{ij}^{\pm }\left[ \pm
n\right] z^{\pm n};  \label{loper}
\end{equation}
we assume, moreover, that $l_{ij}^{\pm }\left[ 0\right] $ are upper (lower)
triangular: $l_{ij}^{+}\left[ 0\right] =l_{ji}^{-}\left[ 0\right] =0$ for $%
i<j.$ The defining relations in ${\cal U}_q^{\prime }\widehat{({\frak gl}_2)}%
_k$ are 
\begin{eqnarray}
l_{ii}^{+}\left[ 0\right] l_{ii}^{-}\left[ 0\right] &=&l_{ii}^{-}\left[
0\right] l_{ii}^{+}\left[ 0\right] =1,i=1,2,  \nonumber \\
R\left( \frac zw\right) L_1^{\pm }\left( z\right) L_2^{\pm }\left( w\right)
&=&L_2^{\pm }\left( w\right) L_1^{\pm }\left( z\right) R\left( \frac
zw\right) ,  \label{krel} \\
R\left( \frac zwq^{-k}\right) L_1^{+}\left( z\right) L_2^{-}\left( w\right)
&=&L_2^{-}\left( w\right) L_1^{+}\left( z\right) R\left( \frac zwq^k\right) 
\nonumber
\end{eqnarray}
(we use again the standard tensor notation). Relations (\ref{krel}) are
understood as relations between formal power series in $\frac zw.$ Observe the
obvious parallels between (\ref{krel}) and (\ref{RTThdual}): the shift in
the argument of $R$ reflects the effects of the central extension; in the
finite-dimensional case this extension is of course trivial and the shift
may be eliminated by a change of variables. The commutation relations (\ref
{krel}) are the exact quantum analogue of the Poisson bracket relations (\ref
{tp12}) with $p=q^k$.

\begin{lemma}
The coefficients of the formal power series 
\[
q\det_{}(L^{\pm }(z)):=l_{11}^{\pm }\left( zq^2\right) \left( l_{22}^{\pm
}\left( z\right) -l_{21}^{\pm }\left( z\right) l_{11}^{\pm }\left( z\right)
^{-1}l_{12}^{\pm }\left( z\right) \right) 
\]
are central in ${\cal U}_q^{\prime }\widehat{({\frak gl}_2)}_k.$
\end{lemma}

\begin{theorem}
\label{current}The quotient of ${\cal U}_q^{\prime }\widehat{({\frak gl}_2)}%
_k$ over the relations $q\det (L^{\pm }(z))=1$ is isomorphic to ${\cal U}_q%
\widehat{({\frak sl}_2)}_k.$
\end{theorem}

The accurate proof of this theorem (and of its generalization to arbitrary
quantized affine Lie algebras) is far from trivial (see \cite{ding}); it
relies on a still another realization of quantized affine Lie algebras, the
so-called {\em Drinfeld's new realization} (\cite{NR}); the important point
for us is that the realization described in theorem \ref{current} is adapted
to the explicit description of the Casimir elements.

Let us put 
\begin{equation}
L\left( z\right) =L^{+}\left( q^{-k/2}z\right) L^{-}\left( zq^{k/2}\right)
^{-1}  \label{curr}
\end{equation}
(this is the so-called {\em full quantum current})

\begin{theorem}
(\cite{RS}) \label{qcenter}(i) The quantum current (\ref{curr}) satisfies
the commutation relations 
\begin{eqnarray}
R\left( \frac zw\right) L_1(z)R^{-1}(\frac{q^{-k/2}z}w)L_2(w)=  \label{commk}
\\
\ \ \ \ \ \ L_2(w)R(\frac{q^{k/2}z}w)L_1(z)R\left( \frac zw\right) ^{-1}. 
\nonumber
\end{eqnarray}
(ii) Suppose that $k=-2;$ then the coefficients of the formal series 
\begin{equation}
t(z)=tr\ q^{2\rho }L\left( z\right)   \label{gene}
\end{equation}
are central elements in the algebra ${\cal U}_q\widehat{({\frak sl}_2)}_k.$
\end{theorem}

As usual, $\rho $ stands for half the sum of positive roots. In the present
setting we have simply $q^{2\rho }=\left( 
\begin{array}{cc}
q & 0 \\ 
0 & q^{-1}
\end{array}
\right) .$ The critical value $k=-2$ is a special case of the general
formula $k=-h^{\vee }$ (the dual Coxeter number); thus the critical value
remains the same as in the non-deformed case.

The definition of the quantum current and formulae (\ref{commk},\ref{gene})
are in fact quite general and apply to arbitrary quantized loop algebras.
There are all reasons to expect that the center of the quotient algebra $%
A_{-h^{\vee }}=A/(c+h^{\vee })$ for the critical value of the central charge
is generated by (\ref{gene}) (with $V$ ranging over all irreducible
finite-dimensional representations of ${\cal U}_q({\frak \hat g}))$ (cf. (%
\cite{etingof}, \cite{EFrenkel})).There are two main subtle points:

\begin{enumerate}
\item  Our \hfill construction\hfill involved\hfill the\hfill 2-dimensional %
\hfill `evaluation\hfill representation' \hfill ${\cal U}_q({\frak sl}_2)%
\rightarrow \allowbreak \break  End{\Bbb C}((z)).$ 
In the general case we must describe finite-dimensional
representations $V(z)$ of ${\cal U}_q({\frak g})$ with spectral parameter;
this requires additional technical efforts.

\item  In the ${\frak sl}_2$ case it was possible to use the generating
series (\ref{loper}) as an exact alternative to the Drinfeld-Jimbo
definition; in the general case the analogues of relations (\ref{krel})
are still valid but there may be others; anyway, the generating series $%
L_{V(z)}^{\pm }$ may be defined as the image of the canonical element: 
\[
L_{V(z)}^{\pm }=\left( id\otimes \rho _{V\left( z\right) }\right) \circ
(R_{\pm }\otimes id)L; 
\]
the definition of the quantum current and theorem \ref{qcenter} are still valid.
\[
\]
\end{enumerate}

We may now use theorem \ref{qcenter} to relate the elements of the center to
quantum Hamiltonians. Our main pattern remains the same. For concreteness
let us again assume that $A={\cal U}_q\widehat{({\frak sl}_2)}.$ Let us take
its dual $A^0$ as the algebra of observables. Let ${\cal L}\in A\otimes A^0$
be the canonical element. For any finite-dimensional ${\frak sl}_2$-module $%
V $ there is a natural representation of $A$ in $V((z))$ (evaluation
representation) Let \ $W$ be another finite-dimensional ${\frak sl}_2$%
-module. Put $L^V(z)=(\rho _{V(z)}M\otimes id){\cal L},L^W(z)=(\rho
_{W(z)}\otimes id){\cal L},R^{VW}(z)=(\rho _{V(z)}\otimes \rho _{W(z)}){\cal %
R}.$ We may call $L^V(z)\in EndV((z))\ \otimes A^0$ the {\em universal
quantum Lax operator} (with auxiliary space $V$). Fix a finite-dimensional
representation $\pi $ of $A^0;$ one can show that $\left( id\otimes \pi
\right) L_V(z)$ is{\em \ rational} in $z.$ (Moreover, \cite{Tar} showed how
to use this dependence on $z$ to classify finite-dimensional representations
of $A^0.)$ Property (\ref{R}) immediately implies that 
\begin{equation}
L_2^W(w)L_1^V(v)=R^{VW}(vw^{-1})L_1^V(v)L_2^W(w)\left(
R^{VW}(vw^{-1})\right) ^{-1}.  \label{Rtt}
\end{equation}
Moreover, $R^{VW}(z)$ satisfies the Yang-Baxter identity (\ref{YBE}).

Let $H$ be the Cartan subgroup of $G$. Fix $h\in H$ and put $%
l_V^h(v)=tr_V\rho _V(h^2)L^V(v)$, then $l_V^h(v)l_W^h(w)=l_W^h(w)l_V^h(v).$
Elements $l_V^h(v)$ are called (twisted) {\em transfer matrices. }

It is convenient to extend $A$ by adjoining to it group-like elements which
correspond to the extended Cartan subalgebra ${\frak \hat h}={\frak h}\oplus 
{\Bbb C}d$ in ${\frak \hat g}{\bf .}$ Let $\hat H=H\times {\Bbb C}^{*}$ be
the corresponding `extended Cartan subgroup' in $A\ $generated by elements 
$h=q^\lambda ,\lambda \in {\frak h},t=q^{kd}.$ Let $\hat H\ \times
A\rightarrow A$ be its natural action on $A\ $ by right translations, $\hat
H\ \times A^0\rightarrow A^0$ the contragredient action. The algebra $A=%
{\cal U}_q({\frak \hat g})$ is factorizable; let $F:A^0\rightarrow A$ be the
isomorphism induced by the action (\ref{fact}); for $h\in \hat H$ put $%
F^h=F\circ h.$ 
\begin{equation}
\begin{array}{c}
L_V^h(z)^{\pm }=\left( \ id\otimes R^{\pm }\circ h\right) L_V(z)\in
EndV\left[ \left[ z^{\pm 1}\right] \right] \otimes A, \\ 
L_V^h(z)=L_V^h(z)^{+}\;\left( id\otimes S\right) L_V^h(z)^{-},h\in \hat H;
\end{array}
\label{TTT}
\end{equation}
it is easy to see that 
\[
L_V^h(z)^{\pm }=\left( id\otimes h^{\pm 1}\right) L_V^1(z)^{\pm } 
\]
and hence 
\[
T_V^h(z)=\left( id\otimes h\right) T_V^1(z)\left( id\otimes h\right) ; 
\]
moreover, if $h=s^{-1}t,$ where $s\in H$ and $t=q^{kd},$ we have also 
\[
T_V^h(z)=\left( \rho _V\left( s\right) \otimes id\right) T_V^t(z)\left( \rho
_V\left( s\right) \otimes id\right) . 
\]
Put $t_V^h=tr_VT_V^h(z)=tr_V\rho _V(s)^2T_V^t(z);$ clearly, we have $%
t_V^h(z)=F^t(l_V^s(z)).$ Let $2\rho $ be the sum of positive roots of $%
{\frak g}={\frak sl}_2,$ $h^{\vee }=2$ its dual Coxeter number (our notation
again hints at the general case).

\begin{theorem}
\label{Cas} (i) (\cite{RS})\hfill Suppose\hfill that\hfill $h=$\\ $q^{-\rho
}q^{-h^{\vee }d}$. Then all coefficients of $t_V^h(z)$ are central in $U.$
(ii) For any $s\in H$ we have $l_V^s(z)=F^{sh^{-1}}(t_V^h(z)).$
\end{theorem}

Thus the duality between Hamiltonians and Casimir operators holds for
quantum affine algebras as well. This allows to anticipate connections
between the generalized Bethe Ansatz, representation theory of quantum
affine algebras at the critical level and the q-KZ equation (\cite{FR}, \cite
{SM}). The results bearing on these connections are already abundant (\cite
{etingof}, \cite{TV}, \cite{EFrenkel}), although they are still not in their
final form.

As we have already observed in the classical context, the Hopf structure on $A^0$ is
perfectly suited to the study of lattice systems. Let $\Delta
^{(N)}:A^0\rightarrow \otimes ^NA^0$ be the iterated coproduct map. $\otimes
^NA^0$ may be interpreted as the algebra of observables associated with a
multiparticle system. Put $\hat t_V^h(z)=\Delta ^{(N)}t_V^h(z);$ Laurent
coefficients of $\hat t_V^h(z)$ provide a commutative family of Hamiltonians
in $\otimes ^NA^0.$ Let $i_n:A^0\rightarrow \otimes ^NA^0$ be the natural
embedding, $i_n:x\mapsto 1\otimes ...\otimes x\otimes ...\otimes 1;$ put $%
L_V^n=\left( id\otimes i_n\right) L_V.$ Then 
\begin{equation}
\hat t_V^h(z)=tr_V\left( \rho _V(h)\overleftarrow{\prod_n}L_V^n\right) .
\label{mon}
\end{equation}
Formula (\ref{mon}) has a natural interpretation in terms of lattice
systems: $L_V^n$ may be regarded as `local' Lax operators attached to the
points of a periodic lattice $\Gamma ={\Bbb Z}/N{\Bbb Z};$ commuting
Hamiltonians for the big system arise from the monodromy matrix $M_V=%
\overleftarrow{\prod }L_V^n$ associated with the lattice. Finally, the twist 
$h\in H$ defines a quasiperiodic boundary condition on the lattice. The
study of the lattice system again breaks into two parts: (a) {\em Find the
joint spectrum of} $\hat l_V(z).$ (b) {\em Reconstruct the Heisenberg
operators corresponding to 'local' observables and compute their correlation
functions.} This is the {\em Quantum Inverse Problem} (profound results on
it are due to \cite{SM}.)

\end{document}